\documentclass[pdflatex,sn-mathphys-num]{sn-jnl}% Math and Physical Sciences Numbered Reference Style

\usepackage{graphicx}%
\usepackage{multirow}%
\usepackage{amsmath,amssymb,amsfonts}%
\usepackage{amsthm}%
\usepackage{mathrsfs}%
\usepackage[title]{appendix}%
\usepackage{xcolor}%
\usepackage{textcomp}%
\usepackage{manyfoot}%
\usepackage{booktabs}%
\usepackage{algorithm}%
\usepackage{algorithmicx}%
\usepackage{algpseudocode}%
\usepackage{listings}%
\usepackage{upgreek}
\usepackage{epsfig}
\usepackage{subcaption}
\graphicspath{ {./fig/} }
\usepackage{cleveref}
\def\eu{\ensuremath{\mathrm{e}}}
\def\iu{\ensuremath{\mathrm{i}}}
\def\du{\ensuremath{\mathrm{d}}}
\def\circv{\ensuremath{\mathring{\mathbf{v}}}}
\def\bfv{\ensuremath{\mathbf{v}}}
\def\hatHone{\ensuremath{\widehat{H}_{\mathrm{I}}}}
\def\hatHtwo{\ensuremath{\widehat{H}_{\mathrm{II}}}}
\def\hatHbulk{\ensuremath{\widehat{H}_{\mathrm{bulk}}}}
\def\rmI{\ensuremath{\mathrm{I}}}
\def\rmII{\ensuremath{\mathrm{II}}}
\def\ehalf{\ensuremath{\mathrm{e}^{\iu\frac{nk}{2}}}}

\def\Inter{\ensuremath{\mathrm{inter}}}
\def\Intra{\ensuremath{\mathrm{intra}}}
\newcommand{\expik}{\eu^{\iu k }}
\newcommand{\expmik}{\eu^{-\iu k }}
\theoremstyle{thmstyleone}%
\newtheorem{theorem}{Theorem}%  meant for continuous numbers
\newtheorem{proposition}[theorem]{Proposition}% 
\newtheorem{lemma}{Lemma}

\theoremstyle{thmstyletwo}%
\newtheorem{remark}{Remark}%

\theoremstyle{thmstylethree}%

\begin{document}
	
	\title[Zero-energy Edge States]{Zero-energy Edge States of Tight-Binding Models for Generalized Honeycomb-Structured  Materials}
	
	
	\author[1]{\fnm{Borui} \sur{Miao}}\email{mbr@mail.tsinghua.edu.cn}
	
	\author[1,2]{\fnm{Yi} \sur{Zhu}}\email{yizhu@tsinghua.edu.cn}

	%\affil*[1]{\orgdiv{Dapartment of Mathematical Sciences}, \orgname{Tsinghua University}, \orgaddress{\city{Beijing}, \postcode{100084}, \country{China}}}
	
	\affil[1]{\orgdiv{Yau Mathematical Sciences Center}, \orgname{Tsinghua University}, \orgaddress{ \city{Beijing}, \postcode{100084}, \country{China}}}
	
	\affil[2]{ \orgname{Beijing Institute of Mathematical Sciences and Applications}, \orgaddress{ \city{Beijing}, \postcode{101408} \country{China}}}
	
	\abstract{
		Generalized honeycomb-structured materials have received increasing attention due to their novel topological properties. 
		In this article, we investigate zero-energy edge states in tight-binding models for such materials with two different interface configurations: type-I and type-II, which are analogous to zigzag and armchair interfaces in the honeycomb structure. We obtain the necessary and sufficient conditions for the existence of such edge states and rigorously prove the existence of paired zero-energy edge states. More specifically, type-II interfaces support two zero-energy states exclusively between topologically distinct materials. For type-I interfaces, zero-energy edge states exist between both topologically distinct and identical materials when hopping coefficients satisfy specific constraints. We further prove that the two energy curves for edge states exhibit a strict crossing. Numerical simulations are presented to demonstrate that wave packets associated with these edge states remain localized near the interface and can be chosen to propagate in either direction, reflecting the pairing of the two edge modes.
	}


	\maketitle	
	\section{Introduction and Summary of Results}
	\subsection{Introduction}
	Owing to their novel physical properties and potential applications, topological materials have attracted considerable attention from researchers in condensed matter physics and quantum mechanics. One of their key features is topologically protected edge states localized near the interface between two topologically distinct materials \cite{Shockley1939,Kane2005,Hasan2010,Rechtsman2013}. Specifically, the interface between two gapped materials with different topological indices gives rise to localized edge states in the band gap. Topological concepts can also be generalized to describe edge states in various physical systems, such as magnetic systems \cite{Ablowitz2020,Shapiro2022}, non-linear systems \cite{Ablowitz2017,Ablowitz2019} and non-Hermitian systems \cite{Ammari2024a, Ammari2024}. \par 	
	A typical system that supports edge states is the 2D honeycomb structure. Its novel properties are related to the conical intersection in the dispersion surfaces, known as the Dirac point \cite{Wallace1947,Neto2009,Fefferman2012,Ammari2020}. A rigorous characterization of these properties was first provided in the pioneering work of Fefferman and Weinstein \cite{Fefferman2012}. Following the framework, the existence of a Dirac point can be extended to general elliptic operators with honeycomb structures \cite{LeeThorp2018}. Ammari and his collaborators also proved similar results for materials in the subwavelength regime by applying layer potential methods to Helmholtz systems \cite{Ammari2018,Ammari2020,Ammari_2020}. In honeycomb structures, a band gap can be created by breaking the inversion symmetry \cite{Neto2009,Schaibley2016,Xiao2007}. 
	Recently, generalized honeycomb-structured materials \cite{Wu2015,Yves2017} have been shown to support the existence of a double Dirac point at the center of the Brillouin zone $ \Gamma $ \cite{Miao2024,Cao2023}. A band gap can also be opened by contracting or dilating the inclusions. The band inversion at $\Gamma$ is commonly regarded as a signature of a topological transition in the local model obtained from perturbations near $\Gamma$; see \cite{Wu2015}. However, the precise definition of a bulk topological invariant is more subtle and remains under discussion \cite{Palmer2021,Xu2023}. In this work, we therefore focus on a direct spectral analysis of the corresponding tight-binding Hamiltonians.

	\par 
	To generate edge states in honeycomb structures, domain wall modulation or sharp termination at the interface between two topologically distinct materials can be employed \cite{LeeThorp2016,Xia2023}. The rigorous characterization of the existence of edge states for domain wall modulation was first provided by \cite{Fefferman_2016}. It has also been verified that the edge or interface configuration determines the existence of edge states \cite{Drouot2020,Fefferman2022,Fefferman2023}. However, analysis of continuous systems remains challenging, especially for sharply-terminated materials that are frequently used in photonics. Therefore, physicists adopt tight-binding approximations that reduce the continuum models to discrete models to obtain insights into the edge states \cite{Fefferman_2017,Fefferman2020,Ammari2022}. In \cite{Fefferman2023}, the authors proved mathematically that flat-band edge states arise for a zigzag-type edge but disappear for an armchair-type edge in tight-binding models. 
	The existence of edge states in generalized honeycomb materials remains largely unexplored, with only a few existing studies addressing this issue \cite{Cao2024}. \par
	
	In this article, we investigate the existence of zero-energy edge states in tight-binding models for generalized honeycomb-structured materials with different interface configurations. Specifically, we study the \emph{zigzag} and \emph{armchair} interfaces, which will be called \emph{type-I} and \emph{type-II} interfaces, respectively. With the Floquet-Bloch transform, we seek to determine the point spectrum of an infinite discrete Hamiltonian operator with quasi-momentum $ k $. We obtain results on:
	\begin{itemize}
		\item[1.] The existence of a pair of zero-energy edge states at $ k=0 $. For the type-I interface, zero-energy edge states may occur between either topologically distinct or identical materials when the parameters satisfy specific conditions. For the type-II interface, they occur only between topologically distinct materials.
		\item[2.] The local behavior near $ k=0 $ given the existence of zero-energy edge states. The point spectrum of the corresponding Hamiltonian operators strictly crosses at $ k=0 $. This result is valid for both types of interfaces. 
	\end{itemize}
	These results provide insights into the dynamical behavior of wave packets; see \cite{Fefferman2013,Watson2018,LeeThorp2018}. Additionally, we numerically simulate the wave packet dynamics at the interfaces.
	
	\subsection{Main Results}
	We give a brief summary of the results presented in this article. Specifically, we establish the following results for type-I and type-II interfaces, respectively. 
	\begin{itemize}
		\item[(1)] \emph{Existence and characterization of zero-energy edge states}: In \Cref{thm:ZeroEigenvalueI} and \Cref{thm:ZeroEigenvalueII}, we provide the necessary and sufficient conditions for the existence of a two-fold eigenvalue for the Hamiltonian operators $ \hatHone(0),\hatHtwo(0) $, respectively. For $ \hatHone(0) $, \Cref{thm:ZeroEigenvalueI} shows that the zero-energy edge states do not generally exist. However, it also demonstrates that zero-energy edge states can exist at interfaces between both topologically distinct and identical materials. Given the existence of zero-energy edge states, the operator $ \hatHone(0) $ has two linearly independent edge states, as characterized in \Cref{thm:ZeroEigenvalueI}. \par 
		For $ \hatHtwo(0) $, \Cref{thm:ZeroEigenvalueII} establishes that zero-energy edge states always exist at the interface between two topologically distinct materials but do not arise at the interface between two topologically identical materials. When zero-energy edge states exist, the operator $ \hatHtwo(0) $ likewise possesses two linearly independent edge states; see \Cref{thm:ZeroEigenvalueII}. We numerically compute the spectra of $ \hatHone(k)$ and $\hatHtwo(k) $; see \Cref{fig:tbedge_zero}, \Cref{fig:tbedge_nonzero}, and \Cref{fig:tbedge_2}.
		\item[(2)] \emph{Local behavior of $ E_{\mathrm{I}}(k)$ and $ E_{\rmII}(k)$ near $ k=0 $}: In \Cref{thm:CrossingCurveI} and \Cref{thm:CrossingCurveII}, we prove that, whenever zero-energy edge states exist, the corresponding eigenvalue branches cross non-tangentially at $k=0$. More specifically, the eigenvalues admit the following asymptotic expansion near $ 0 $
		\begin{gather*}
			\begin{split}
				E_{\rmI,\pm}(k) &= \pm E_{\rmI}^{(1)}|k| + \mathcal{O}(k^2),\\
				E_{\rmII,\pm}(k) &= \pm E_{\rmII}^{(1)}|k| + \mathcal{O}(k^2),
			\end{split}
		\end{gather*}
		for small but nonzero $ k $. Here $ E_{\rmI}^{(1)},E_{\rmII}^{(1)} $ are nonzero constants. These results indicate that a wave packet localized at the interface can propagate in both directions. As proved for continuous systems \cite{Fefferman2013,Watson2018}, the group velocity of wave packets is given by $ \partial_{k}E(0) $. Thus, one may expect that a wave packet localized at either a type-I or type-II interface can propagate in either direction if the initial value is properly chosen.
		\item[(3)] \emph{Numerical simulation of wave packets that concentrate near $ E_{\rmI}(0)=E_{\rmII}(0)=0 $}: We numerically simulate the time dynamics of wave packets propagating through a local defect in discrete Schr\"odinger systems  
		\[ \iu \partial_{t} \Phi = \widehat{H}\Phi, \]
		 to illustrate the theoretical results and present some intuitions for further studies. The results are shown in \Cref{fig:edgeIdef_diff}, \Cref{fig:edgeIdef_same} and \Cref{fig:edgeIIdef}.
		Since the wave packets are localized in the band gap, the numerical simulations show no significant radiation
		into the bulk material. These numerical observations suggest potential applications in designing bidirectional waveguides. 
		
	\end{itemize}
	\subsection{Structure of the Article}
	This paper is arranged as follows. In \Cref{sec:ProbForm}, we introduce the overall problem formulation. define the Hamiltonian operator $ \widehat{H}_{\mathrm{full}} $ for tight-binding models. By appropriately defining the coefficients of the Hamiltonian operator and applying Floquet-Bloch reduction, we obtain the Hamiltonians $ \hatHone(k) $ and $ \hatHtwo(k) $ for type-I and type-II interfaces, respectively, where $ k\in [-\pi,\pi) $. In \Cref{sec:TypeI} and \Cref{sec:TypeII}, we first analyze the existence of zero-energy edge states in type-I and type-II interfaces, respectively. We then derive the asymptotic expansion of eigenvalues near $ k=0 $ and $ E_{\rmI}(0) = E_{\rmII}(0) = 0 $. In \Cref{sec:Numerics}, we present some numerical illustrations of previous results. We also simulate the time evolution of wave packets as they propagate through a local defect. In the appendix, we discuss the properties of the bulk spectrum in \Cref{apsec:bulk}.

	\section{Problem Formulation}\label{sec:ProbForm}
	
	In this section, we introduce the tight-binding model for generalized honeycomb structures with nearest neighborhood approximation. The tight-binding model arises from reducing eigenvalue problems on the whole space $ \mathbb{R}^{d} $ to problems on lattice sites. 
	\subsection{Lattice Points}
\begin{figure}[htbp]\centering
		\includegraphics[width=1\linewidth]{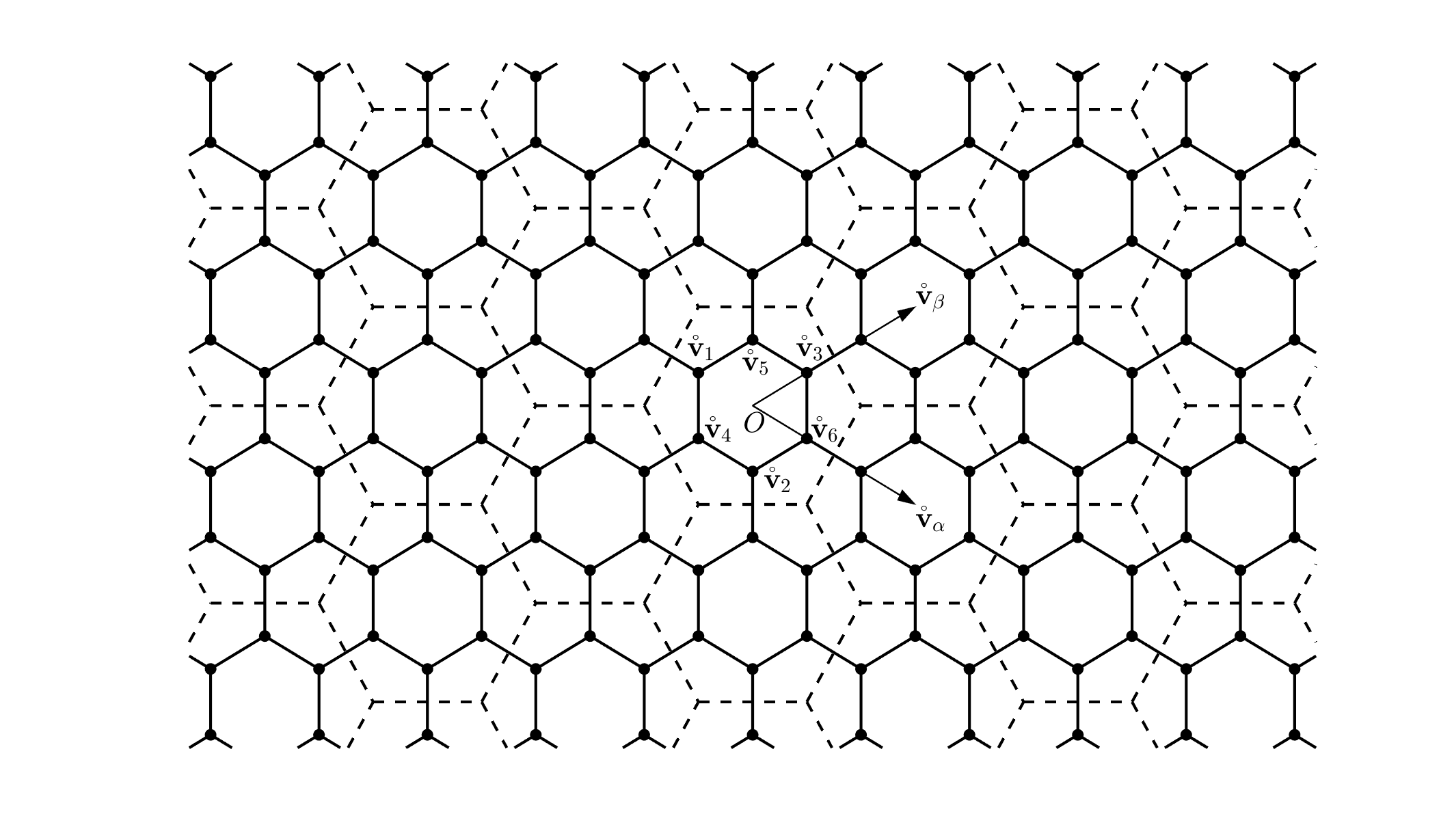}
		\caption{Generalized Honeycomb Structure. The dashed lines denote the edge of a unit hexagonal cell. }
		\label{fig:Material}
	\end{figure}

	To define the geometry, we first define the \emph{triangular lattice} $ \Lambda \triangleq \mathbb{Z} \circv_\alpha \oplus \mathbb{Z}\circv_\beta $, where its basis vectors are given by 
	\[ \circv_\alpha \triangleq (\sqrt{3}/2,-1/2)^{T},\quad \circv_\beta \triangleq (\sqrt{3}/2,1/2)^{T}. \]
	The \emph{generalized honeycomb sites} are given by 
	\[ \mathbb{GH} \triangleq \bigcup_{j=1}^6{\circv_j+\Lambda}, \]
	where the vectors $ \{ \circv_j \}_{j=1}^6 $ are given by 
	\begin{gather*}
		\begin{aligned}
			&\circv_1 \triangleq -\circv_\alpha/3,\quad \circv_2 \triangleq\circv_\alpha/3-\circv_\beta/3 ,\quad\circv_3 \triangleq\circv_\beta/3 ,\\
			&\circv_4 \triangleq -\circv_3,\quad \circv_5 \triangleq -\circv_2,\quad \circv_6 \triangleq -\circv_1.
		\end{aligned}
	\end{gather*}
	See \Cref{fig:Material} for an example.
	\par 
	For any $ \omega = \circv_j + \lambda \in \circv_j+ \Lambda\subset \mathbb{GH} $, the \emph{nearest neighborhoods $ \mathcal{N}_\omega $} are given by
	\begin{gather}
		\mathcal{N}_\omega\triangleq\left\{\begin{split}
			&\{\omega + \circv_k\}_{k=1,2,3},\quad j = 1,2,3,\\
			&\{\omega + \circv_k\}_{k = 4,5,6},\quad j = 4,5,6,
		\end{split}\right.
	\end{gather}
	where $ \lambda\in \Lambda $ denotes the lattice point in the triangular lattice $ \Lambda $.
	The sets of nearest neighbors for $ \omega = \circv_j + \lambda $ are given by 
	\begin{gather}\label{eqn:NearestNeighbour}
		\mathcal{N}_\omega=\left\{\begin{split}
			&\{\circv_4+\lambda,\circv_5+\lambda,\circv_6+(\lambda-\circv_{\alpha})\},\quad j=1,\\
			&\{\circv_4+\lambda,\circv_5+(\lambda+\circv_{\alpha}-\circv_{\beta}),\circv_6+\lambda\},\quad  j=2,\\
			&\{\circv_4+(\lambda+\circv_{\beta}),\circv_5+\lambda,\circv_6+\lambda\},\quad  j=3,\\
			&\{\circv_1+\lambda,\circv_2+\lambda,\circv_3+(\lambda-\circv_{\beta})\},\quad j=4,\\
			&\{\circv_1+\lambda,\circv_2+(\lambda-\circv_{\alpha}+\circv_{\beta}),\circv_3+\lambda\},\quad  j=5,\\
			&\{\circv_1+(\lambda+\circv_{\alpha}),\circv_2+\lambda,\circv_3+\lambda\},\quad  j=6.
		\end{split}\right.
	\end{gather} 
	Note that the definition of the set of nearest neighbors $ \mathcal{N}_{\omega} $ depends on the lattice point $ \lambda $. One also notes that the nearest neighborhoods can be categorized into two types, namely \emph{intracell and intercell sets of nearest neighbors} $ \mathcal{N}_{\omega,\Intra}, \mathcal{N}_{\omega,\Inter} $, 
	\begin{gather}\label{eqn:NNint}
		\begin{aligned}
			\mathcal{N}_{\omega,\Intra}&\triangleq\left\{\begin{split}
				&\{\circv_4+\lambda,\circv_5+\lambda\},\quad j=1,\\
				&\{\circv_4+\lambda,\circv_6+\lambda\},\quad  j=2,\\
				&\{\circv_5+\lambda,\circv_6+\lambda\},\quad  j=3,\\
				&\{\circv_1+\lambda,\circv_2+\lambda\},\quad j=4,\\
				&\{\circv_1+\lambda,\circv_3+\lambda\},\quad  j=5,\\
				&\{\circv_2+\lambda,\circv_3+\lambda\},\quad  j=6.
			\end{split}\right.\\
			\mathcal{N}_{\omega,\Inter}&\triangleq\left\{\begin{split}
				&\{\circv_6+(\lambda-\circv_{\alpha})\},\quad j=1,\\
				&\{\circv_5+(\lambda+\circv_{\alpha}-\circv_{\beta})\},\quad  j=2,\\
				&\{\circv_4+(\lambda+\circv_{\beta})\},\quad  j=3,\\
				&\{\circv_3+(\lambda-\circv_{\beta})\},\quad j=4,\\
				&\{\circv_2+(\lambda-\circv_{\alpha}+\circv_{\beta})\},\quad  j=5,\\
				&\{\circv_1+(\lambda+\circv_{\alpha})\},\quad  j=6.
			\end{split}\right.
		\end{aligned}
	\end{gather}
	\subsection{The Full Hamiltonian Operator}
	Given the configuration of the discrete system, we define the wave function on each site $\Psi\in  l^2(\mathbb{GH}) \simeq l^2(\Lambda;\mathbb{C}^6)$. 
	\[ \Psi \triangleq \{( \psi_{_{\circv_j+\lambda}} )_{j=1}^6\}_{\lambda \in \Lambda}=\left\{( \psi_{_{j,\lambda}} )_{j=1}^6\right\}_{\lambda \in \Lambda}. \]
	Here the Hilbert space $ l^2(\Lambda;\mathbb{C}^6) $ is equipped with an inner product $ (\cdot,\cdot) $ given by
	\[ (\Phi,\Psi) \triangleq \sum_{\lambda\in \Lambda}\sum_{j=1}^6 \overline{\phi}_{_{j,\lambda}} \psi_{_{j,\lambda}},\quad \Phi,\Psi\in l^2(\Lambda;\mathbb{C}^6).   \]
	From \eqref{eqn:NearestNeighbour}, the \emph{full Hamiltonian operator} $ \widehat{H}_{\mathrm{full}} $ acting on $ \Psi\in l^2(\Lambda;\mathbb{C}^6) $, with nearest neighborhood approximation, is given by
	\begin{gather}\label{eqn:HamiltonFull}
		\begin{aligned}
			(\widehat{H}_{\mathrm{full}}\Psi)_{_{1,\lambda}} &\triangleq -b_{_\lambda} \psi_{_{4,\lambda}}  -b_{_{\lambda}} \psi_{_{5,\lambda}} -a_{_{\lambda}}\psi_{_{6,\lambda-\circv_\alpha}},\\
			(\widehat{H}_{\mathrm{full}}\Psi)_{_{2,\lambda}} &\triangleq -b_{_{\lambda}} \psi_{_{4,\lambda}}  -d_{_{\lambda+\circv_\alpha-\circv_\beta}} \psi_{_{5,\lambda+\circv_\alpha-\circv_\beta}} -b_{_{\lambda}}\psi_{_{6,\lambda}},\\
			(\widehat{H}_{\mathrm{full}}\Psi)_{_{3,\lambda}} &\triangleq -c_{_\lambda} \psi_{_{4,\lambda+\circv_\beta}} - b_{_{\lambda}} \psi_{_{5,\lambda}} -b_{_\lambda}\psi_{_{6,\lambda}},\\
			(\widehat{H}_{\mathrm{full}}\Psi)_{_{4,\lambda}} &\triangleq -b_{_{\lambda}} \psi_{_{1,\lambda}}-b_{_\lambda}\psi_{_{2,\lambda}} -c_{_{\lambda-\circv_{\beta}}} \psi_{_{3,\lambda-\circv_{\beta}}},\\
			(\widehat{H}_{\mathrm{full}}\Psi)_{_{5,\lambda}} &\triangleq -b_{_{\lambda}} \psi_{_{1,\lambda}}-d_{_\lambda}\psi_{_{2,\lambda-\circv_{\alpha}+\circv_{\beta}}} -b_{_\lambda} \psi_{_{3,\lambda}},\\
			(\widehat{H}_{\mathrm{full}}\Psi)_{_{6,\lambda}} &\triangleq -a_{_{\lambda+\circv_{\alpha}}} \psi_{_{1,\lambda+\circv_{\alpha}}}-b_{_\lambda}\psi_{_{2,\lambda}} -b_{_\lambda} \psi_{_{3,\lambda}}.
		\end{aligned}
	\end{gather}
	The \emph{hopping coefficients} between nearest neighborhoods $ \{ a_{_\lambda},b_{_\lambda},c_{_\lambda},d_{_\lambda} \}_{\lambda\in \Lambda} $ are strictly positive real numbers with uniform lower and upper bounds. 
	\begin{equation}\label{eqn:HoppingCoef}
		\max_{\lambda\in \Lambda}\{a_{_\lambda},b_{_\lambda},c_{_\lambda},d_{_\lambda}\}\le C <+\infty,\quad \min_{\lambda\in \Lambda}\{a_{_\lambda},b_{_\lambda},c_{_\lambda},d_{_\lambda}\}\ge c>0.
	\end{equation}
	Here $ \{b_{_{\lambda}}\}_{\lambda\in \Lambda} $ denotes the \emph{intracell hopping coefficients}, while $ \{ a_{_\lambda},c_{_\lambda},d_{_\lambda} \}_{\lambda\in \Lambda} $ denote the \emph{intercell hopping coefficients}. It is clear that the full Hamiltonian operator $ \widehat{H}_{\mathrm{full}} $ is self-adjoint and bounded from the Hilbert space $ l^2(\Lambda;\mathbb{C}^6) $ to itself. Therefore, its spectrum satisfies $ \sigma(\widehat{H}_{\mathrm{full}})\subset \mathbb{R} $.

	\subsection{Interface Configuration and Floquet-Bloch Decomposition}
	In this article, we investigate two types of interface configurations, namely the type-I interface and the type-II interface. \par
\begin{figure}[htbp] 
		\centering
		\includegraphics[width=0.5\textwidth]{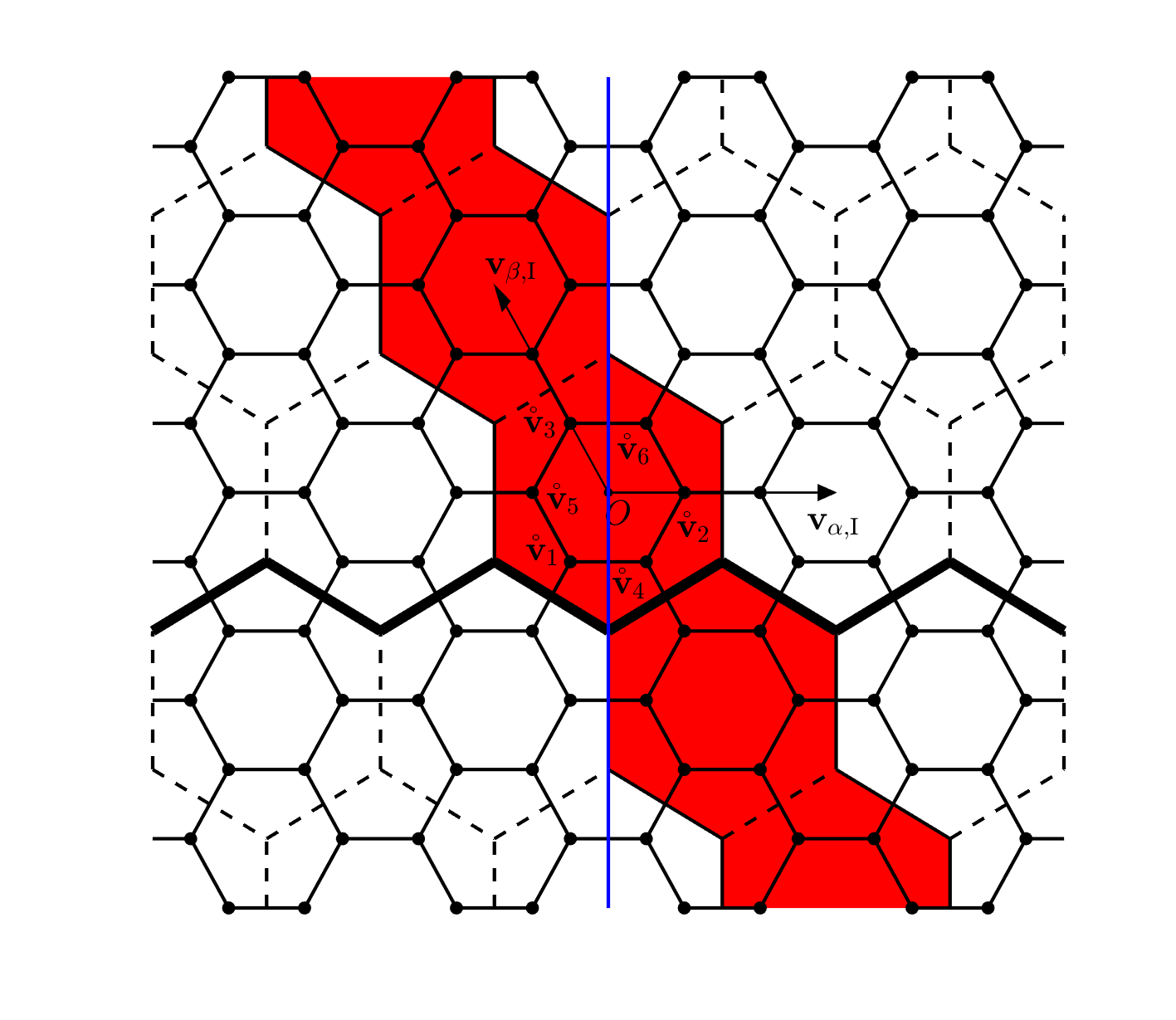}
		\hspace{-0.3cm}
		\includegraphics[width=0.5\textwidth]{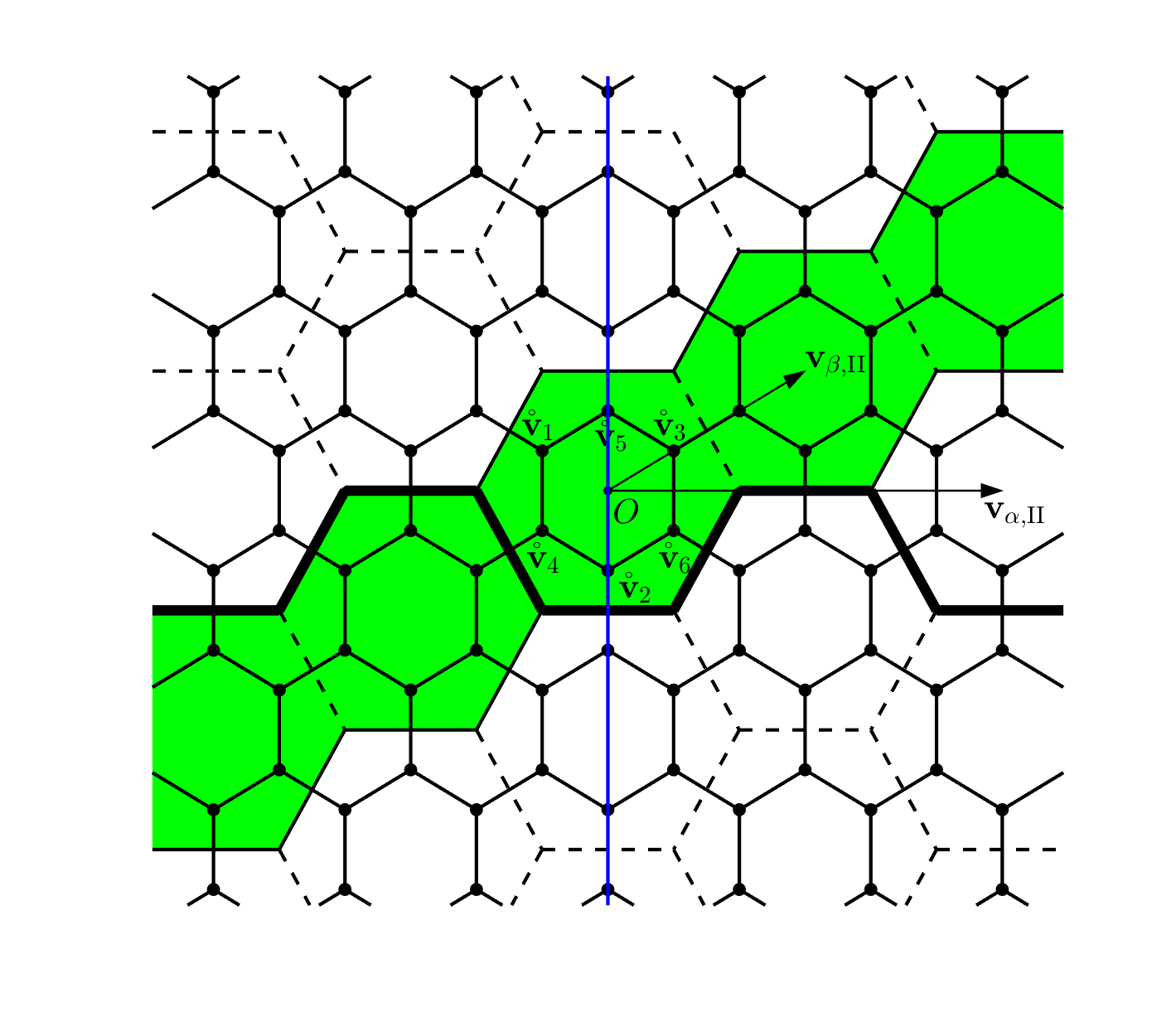}
		\caption{Type-I and type-II interfaces and the unit cell. Left panel: An example of a type-I interface. The cells in red denote the unit cell of the type-I interface. Right panel: An example of a type-II interface. The cells in green denote the unit cell of the type-II interface. In both cases, the unit cell is periodic in the $ \bfv_{\alpha,\rmI}(\bfv_{\alpha,\rmII}) $ direction and extends in the $ \bfv_{\beta,\rmI}(\bfv_{\beta,\rmII}) $ direction. The bold black lines denote the interface between the two materials. The blue lines denote the reflection axes. }
		\label{fig:CellExample}
	\end{figure}

 
	We now give a rigorous definition for both types of interfaces. First, we define two infinite subsets of the generalized honeycomb structure $ \mathbb{GH}_{\pm} $, which satisfy $ \mathbb{GH}_{+}\cup\mathbb{GH}_{-}=\mathbb{GH} $ and $ \mathbb{GH}_{+}\cap\mathbb{GH}_{-}=\emptyset $. Following the previous definitions, we define the \emph{interior sites} and \emph{interface sites} of $ \mathbb{GH}_{\pm} $ by 
	\begin{gather}
		\begin{split}
			\mathrm{int}\big( \mathbb{GH}_{\pm} \big) &\triangleq \big\{ \omega: \mathcal{N}_{\omega} \subseteq \mathbb{GH}_{\pm}, \omega\in \mathbb{GH}_{\pm} \big\},\\
			\partial \big(\mathbb{GH}_{\pm}\big) &\triangleq \big\{ \omega: \mathcal{N}_{\omega} \nsubseteq \mathbb{GH}_{\pm}, \omega\in \mathbb{GH}_{\pm} \big\}.
		\end{split} 
	\end{gather}
	It follows easily that $ \mathbb{GH}_{\pm} = \mathrm{int}( \mathbb{GH}_{\pm} )\cup \partial (\mathbb{GH}_{\pm}) $.
	As shown in \Cref{fig:CellExample}, both physical systems are $ \bfv_{\alpha} $-periodic and extend in the $ \bfv_{\beta}$ direction. 
	Hence, we define 
	\begin{gather}
		\begin{aligned}
			\bfv_{\alpha,\rmI} &= \circv_{\alpha}-\circv_{\beta},\quad \bfv_{\beta,\rmI} = \circv_{\beta},\quad \text{for type I edge},\\
			\bfv_{\alpha,\rmII} &= \circv_{\alpha}+\circv_{\beta},\quad \bfv_{\beta,\rmII} = \circv_{\beta},\quad \text{for type II edge}.
		\end{aligned}
	\end{gather} 
	For both cases, one easily has 
	\[ \Lambda = \{ m\bfv_{\alpha,\rmI}+n\bfv_{\beta,\rmI}:m,n\in \mathbb{Z} \}  = \{ m\bfv_{\alpha,\rmII}+n\bfv_{\beta,\rmII}:m,n\in \mathbb{Z} \}. \]
	Therefore the components of the wave function $ \{( \psi_{_{j,\lambda}} )_{j=1}^6\}_{\lambda \in \Lambda} $ and hopping coefficients $ \{ a_{_\lambda},b_{_\lambda},c_{_\lambda},d_{_\lambda} \}_{\lambda\in \Lambda} $ can be relabeled as $ \Psi = \{ (\psi_{_{j,m,n}})_{j=1}^{6} \}_{m,n\in\mathbb{Z}} $ and $ \{ a_{_{m,n}},b_{_{m,n}},c_{_{m,n}},d_{_{m,n}} \}_{m,n\in \mathbb{Z}} $, respectively. Here we emphasize that for different types of interface, the hopping coefficients $ \{ a_{_{m,n}},b_{_{m,n}},c_{_{m,n}},d_{_{m,n}} \}_{m,n\in \mathbb{Z}} $ correspond to different hopping coefficients $ \{ a_{_\lambda},b_{_\lambda},c_{_\lambda},d_{_\lambda} \}_{\lambda\in \Lambda} $ defined in \eqref{eqn:HamiltonFull}. \par 
	With a slight abuse of notation, for the type-I interface, we let
	
	\[ \mathbb{GH}_{+} = \bigcup_{m\in\mathbb{Z},n\ge 0}\bigg[\bigcup_{j=1}^{6}\{\circv_{j}+m\mathbf{v}_{\alpha,\rmI}+n\mathbf{v}_{\beta,\rmI}\} \bigg], \quad \mathbb{GH}_{-} = \mathbb{GH}\backslash \mathbb{GH}_{+}.\]
	Similarly, for the type-II interface, we let 
	\[ \mathbb{GH}_{+} = \bigcup_{m\in\mathbb{Z},n\ge 0}\bigg[\bigcup_{j=1}^{6}\{\circv_{j}+m\mathbf{v}_{\alpha,\rmII}+n\mathbf{v}_{\beta,\rmII}\} \bigg], \quad \mathbb{GH}_{-} = \mathbb{GH}\backslash \mathbb{GH}_{+}.\]
	
	In this article, we specialize to the case where the intracell and intercell hopping coefficients at the interior sites $\mathrm{int}(\mathbb{GH}_{\pm})$ are given by
		\[ b_{_{\lambda}} = b_\pm, a_{_\lambda} = c_{_\lambda} = d_{_\lambda} = b_\pm + \delta_\pm,\] 
		respectively. At the interface sites $\partial(\mathbb{GH}_{\pm})$, they are given by $ b_{_\pm} $ and $ c $, respectively. Therefore, the Hamiltonian operator can be written out as follows. For $ \omega\in \mathrm{int}(\mathbb{GH}_{\pm}) $, one has
	\begin{gather}
		\begin{aligned}
			(\widehat{H}_{\mathrm{full}}\Psi)_{_{1,\lambda}} &\triangleq -b_{_\pm} \psi_{_{4,\lambda}}  -b_{_{\pm}} \psi_{_{5,\lambda}} -(b_{_{\pm}}+\delta_{_{\pm}})\psi_{_{6,\lambda-\circv_\alpha}},\\
			(\widehat{H}_{\mathrm{full}}\Psi)_{_{2,\lambda}} &\triangleq -b_{_{\pm}} \psi_{_{4,\lambda}}  - (b_{_{\pm}}+\delta_{_{\pm}}) \psi_{_{5,\lambda+\circv_\alpha-\circv_\beta}} -b_{_{\pm}}\psi_{_{6,\lambda}},\\
			(\widehat{H}_{\mathrm{full}}\Psi)_{_{3,\lambda}} &\triangleq -(b_{_{\pm}}+\delta_{_{\pm}}) \psi_{_{4,\lambda+\circv_\beta}} - b_{_{\pm}} \psi_{_{5,\lambda}} -b_{_\pm}\psi_{_{6,\lambda}},\\
			(\widehat{H}_{\mathrm{full}}\Psi)_{_{4,\lambda}} &\triangleq -b_{_{\pm}} \psi_{_{1,\lambda}}-b_{_\pm}\psi_{_{2,\lambda}} -(b_{_{\pm}}+\delta_{_{\pm}})\psi_{_{3,\lambda-\circv_{\beta}}},\\
			(\widehat{H}_{\mathrm{full}}\Psi)_{_{5,\lambda}} &\triangleq -b_{_{\pm}} \psi_{_{1,\lambda}}-(b_{_{\pm}}+\delta_{_{\pm}})\psi_{_{2,\lambda-\circv_{\alpha}+\circv_{\beta}}} -b_{_\pm} \psi_{_{3,\lambda}},\\
			(\widehat{H}_{\mathrm{full}}\Psi)_{_{6,\lambda}} &\triangleq -(b_{_{\pm}}+\delta_{_{\pm}}) \psi_{_{1,\lambda+\circv_{\alpha}}}-b_{_\pm}\psi_{_{2,\lambda}} -b_{_\pm} \psi_{_{3,\lambda}}.
		\end{aligned}
	\end{gather}
	For $\omega\in \partial (\mathbb{GH}_{\pm})$, one has 
	\begin{gather}
		\begin{aligned}
			(\widehat{H}_{\mathrm{full}}\Psi)_{_{1,\lambda}} &\triangleq -b_{_\pm} \psi_{_{4,\lambda}}  -b_{_{\pm}} \psi_{_{5,\lambda}} -c\psi_{_{6,\lambda-\circv_\alpha}},\\
			(\widehat{H}_{\mathrm{full}}\Psi)_{_{2,\lambda}} &\triangleq -b_{_{\pm}} \psi_{_{4,\lambda}}  - c \psi_{_{5,\lambda+\circv_\alpha-\circv_\beta}} -b_{_{\pm}}\psi_{_{6,\lambda}},\\
			(\widehat{H}_{\mathrm{full}}\Psi)_{_{3,\lambda}} &\triangleq -c \psi_{_{4,\lambda+\circv_\beta}} - b_{_{\pm}} \psi_{_{5,\lambda}} -b_{_\pm}\psi_{_{6,\lambda}},\\
			(\widehat{H}_{\mathrm{full}}\Psi)_{_{4,\lambda}} &\triangleq -b_{_{\pm}} \psi_{_{1,\lambda}}-b_{_\pm}\psi_{_{2,\lambda}} -c\psi_{_{3,\lambda-\circv_{\beta}}},\\
			(\widehat{H}_{\mathrm{full}}\Psi)_{_{5,\lambda}} &\triangleq -b_{_{\pm}} \psi_{_{1,\lambda}}-c\psi_{_{2,\lambda-\circv_{\alpha}+\circv_{\beta}}} -b_{_\pm} \psi_{_{3,\lambda}},\\
			(\widehat{H}_{\mathrm{full}}\Psi)_{_{6,\lambda}} &\triangleq -c \psi_{_{1,\lambda+\circv_{\alpha}}}-b_{_\pm}\psi_{_{2,\lambda}} -b_{_\pm} \psi_{_{3,\lambda}}.
		\end{aligned}
	\end{gather}

	From a physical point of view, the hopping coefficient between two sites is larger when the corresponding atoms or inclusions are closer, and vice versa; see \cite{Ammari2022} for example. Thus, the situation in which $ \delta_{+}\delta_{-}<0 $ models the connection of two materials with different configurations, for instance the contracted and dilated generalized honeycomb structures \cite{Miao2024}. In this article, we also discuss the situation in which $ \delta_{+}\delta_{-}>0 $, which corresponds to an interface between two materials of the same type. It is also worth noting that the parameter $ c $ denotes the hopping coefficient across the interface, which is a positive constant determined by the underlying physical system.

	Before presenting the interface problems, we briefly review the bulk properties of the problem, namely $b_{_\lambda} = b$ and $a_{_\lambda}=c_{_\lambda}=d_{_\lambda}= b+\delta$. By the Floquet-Bloch decomposition, we can prove that there is a double Dirac point at the center of the Brillouin zone when $\delta=0$, while a band gap opens when $ \delta\neq 0 $. We refer to \Cref{apsec:bulk} for a detailed review of bulk properties of the generalized honeycomb structure. 
	
	From the definition, it is clear that for both types of interfaces, the operators depend on the parameters $ (b_{\pm},\delta_{\pm},c) $. We denote the corresponding operators by  $ \widehat{H}_{\mathrm{I}},\widehat{H}_{\mathrm{II}}:l^2(\Lambda;\mathbb{C}^6)\to l^2(\Lambda;\mathbb{C}^6) $. We describe the intrinsic symmetries of these two interface operators below.  \par 
	Again, by the Floquet-Bloch decomposition, one can decompose the Hilbert space $ l^2(\Lambda;\mathbb{C}^6) $ into a direct integral over the space $ l^2_{k}(\mathbb{Z};\mathbb{C}^6) $, with quasi-momentum $ k\in [-\pi,\pi) $. Here $ l^2_{k}(\mathbb{Z};\mathbb{C}^6) $ denotes the space of $ k $-quasiperiodic wave functions $ \boldsymbol{\upphi} = \{ (\phi_{_{j,n}})_{j=1}^{6} \}_{n\in\mathbb{Z}} $ with 
	\[ \Psi = \{(\psi_{_{j,m,n}}=\eu^{\iu k m}\phi_{_{j,n}})_{j=1}^6\}_{m,n\in\mathbb{Z}},\; j=1,2,\ldots,6,\]
	such that 
	\[\sum_{n\in \mathbb{Z}}\sum_{j=1}^6|\phi_{_{j,n}}|^2<+\infty.\]
	The Hamiltonian operators $ \widehat{H}_{\mathrm{I}},\widehat{H}_{\mathrm{II}} $ can also be decomposed into a direct integral of operators $ \hatHone(k),\hatHtwo(k) $ for $ k\in [-\pi,\pi) $, acting on $ \mathbf{u},\mathbf{w}\in l^2_{k}(\mathbb{Z};\mathbb{C}^6) $. For the type-I interface, one has  
	\begin{gather}\label{eqn:OneHamiltonian}
		\begin{aligned}
			(\hatHone(k)\mathbf{u})_{_{1,n}} &= -b_{_{n,\rmI}}u_{_{4,n}}-b_{_{n,\rmI}}u_{_{5,n}}-c_{_{n-1,\rmI}}\expmik u_{_{6,n-1}},\\
			(\hatHone(k)\mathbf{u})_{_{2,n}} &= -b_{_{n,\rmI}} u_{_{4,n}}-d_{_{n,\rmI}}\expik u_{_{5,n}}-b_{_{n,\rmI}}u_{_{6,n}},\\
			(\hatHone(k)\mathbf{u})_{_{3,n}} &= -c_{_{n,\rmI}}u_{_{4,n+1}} -b_{_{n,\rmI}}u_{_{5,n}}-b_{_{n,\rmI}} u_{_{6,n}},\\
			(\hatHone(k)\mathbf{u})_{_{4,n}} &= -b_{_{n,\rmI}} u_{_{1,n}}-b_{_{n,\rmI}}u_{_{2,n}}-c_{_{n-1,\rmI}}u_{_{3,n-1}},\\
			(\hatHone(k)\mathbf{u})_{_{5,n}} &= -b_{_{n,\rmI}} u_{_{1,n}}-d_{_{n,\rmI}}\expmik u_{_{2,n}}-b_{_{n,\rmI}}u_{_{3,n}},\\
			(\hatHone(k)\mathbf{u})_{_{6,n}} &= -c_{_{n,\rmI}}\expik u_{_{1,n+1}}-b_{_{n,\rmI}}u_{_{2,n}}-b_{_{n,\rmI}}u_{_{3,n}}.
		\end{aligned}
	\end{gather}
	The coefficients $ \{ a_{_{n,\rmI}},b_{_{n,\rmI}},c_{_{n,\rmI}},d_{_{n,\rmI}} \} $ are given by 
	\begin{gather}\label{eqn:CoefI}
		\begin{split}		
			a_{_{n,\rmI}} = \left\{ 
			\begin{aligned}
				&b_{+}+\delta_{+},\quad n\ge 1,\\
				&c ,\quad n = 0,\\
				&b_{-}+\delta_{-},\quad n\le -1,
			\end{aligned}\right.\quad 
			b_{_{n,\rmI}} = \left\{
			\begin{aligned}
				&b_{+},\quad n\ge 0,\\
				&b_{-},\quad n\le -1,
			\end{aligned}\right.\\
			c_{_{n,\rmI}} = a_{_{n+1,\rmI}},\quad
			d_{_{n,\rmI}} = \left\{
			\begin{aligned}
				&b_{+}+\delta_{+},\quad n\ge 0,\\
				&b_{-}+\delta_{-},\quad n\le -1.
			\end{aligned}\right.
		\end{split}
	\end{gather}
	For the type-II interface, one has 
	\begin{gather}\label{eqn:TwoHamiltonian}
		\begin{aligned}
			(\hatHtwo(k)\mathbf{w})_{_{1,n}} &= -b_{_{n,\rmII}}w_{_{4,n}}-b_{_{n,\rmII}}w_{_{5,n}}-c_{_{n,\rmII}}\expmik w_{_{6,n+1}},\\
			(\hatHtwo(k)\mathbf{w})_{_{2,n}} &= -b_{_{n,\rmII}} w_{_{4,n}}-d_{_{n-2,\rmII}}\expik w_{_{5,n-2}}-b_{_{n,\rmII}}w_{_{6,n}},\\
			(\hatHtwo(k)\mathbf{w})_{_{3,n}} &=-c_{_{n,\rmII}}w_{_{4,n+1}} -b_{_{n,\rmII}}w_{_{5,n}}-b_{_{n,\rmII}} w_{_{6,n}},\\
			(\hatHtwo(k)\mathbf{w})_{_{4,n}} &= -b_{_{n,\rmII}} w_{_{1,n}}-b_{_{n,\rmII}}w_{_{2,n}}-c_{_{n-1,\rmII}}w_{_{3,n-1}},\\
			(\hatHtwo(k)\mathbf{w})_{_{5,n}} &= -b_{_{n,\rmII}} w_{_{1,n}}-d_{_{n,\rmII}}\expmik w_{_{2,n+2}}-b_{_{n,\rmII}}w_{_{3,n}},\\
			(\hatHtwo(k)\mathbf{w})_{_{6,n}} &= -c_{_{n-1,\rmII}}\expik w_{_{1,n-1}}-b_{_{n,\rmII}}w_{_{2,n}}-b_{_{n,\rmII}}w_{_{3,n}}.
		\end{aligned}
	\end{gather}
	The coefficients $ \{ a_{_{n,\rmII}},b_{_{n,\rmII}},c_{_{n,\rmII}},d_{_{n,\rmII}} \} $ are given by
	\begin{gather}\label{eqn:CoefII}
		\begin{split}		
			a_{_{n,\rmII}} = \left\{ 
			\begin{aligned}
				&b_{+}+\delta_{+},\quad n\ge 0,\\
				&c ,\quad n = -1,\\
				&b_{-}+\delta_{-},\quad n\le -2,
			\end{aligned}\right.\quad 
			b_{_{n,\rmII}} = \left\{
			\begin{aligned}
				&b_{+},\quad n\ge 0,\\
				&b_{-},\quad n\le -1,
			\end{aligned}\right.\\
			c_{_{n,\rmII}} = a_{_{n,\rmII}},\quad
			d_{_{n,\rmII}} = \left\{ 
			\begin{aligned}
				&b_{+}+\delta_{+},\quad n\ge 0,\\
				&c ,\quad n = -1,-2,\\
				&b_{-}+\delta_{-},\quad n\le -3.
			\end{aligned}\right.
		\end{split}
	\end{gather}
	The operators are self-adjoint and bounded from $ l^2_k(\mathbb{Z};\mathbb{C}^6) $ to itself. The point spectra near $ 0 $ are denoted by $ E_{\rmI}(k),E_{\rmII}(k) $, respectively. Here we adopt different notations ($ \mathbf{u}$ and $\mathbf{w} $) for wave functions to signify the difference between the two types of interfaces after Floquet-Bloch decomposition.\par 
	
	\section{Edge States for Type I Interface}\label{sec:TypeI}
	In this section, we first investigate the point spectrum at $ E_{\rmI}(k)=0 $ for the type-I interface near $ k=0 $. More specifically, we look for nonzero solutions to the following equations that decay as $ n\to \pm \infty $:
	\begin{gather}\label{eqn:ZeroIEdge}
		\begin{split}
			&\left\{\begin{aligned}
				& -b_{_{n,\rmI}}u_{_{4,n}}-b_{_{n,\rmI}}u_{_{5,n}}-c_{_{n-1,\rmI}}\expmik u_{_{6,n-1}}=0,\\
				& -b_{_{n,\rmI}} u_{_{4,n}}-d_{_{n,\rmI}}\expik u_{_{5,n}}-b_{_{n,\rmI}}u_{_{6,n}}=0,\\
				& -c_{_{n,\rmI}}u_{_{4,n+1}} -b_{_{n,\rmI}}u_{_{5,n}}-b_{_{n,\rmI}} u_{_{6,n}}=0,
			\end{aligned}\right.\\
			&\left\{\begin{aligned}
				& -b_{_{n,\rmI}} u_{_{1,n}}-b_{_{n,\rmI}}u_{_{2,n}}-c_{_{n-1,\rmI}}u_{_{3,n-1}}=0,\\
				& -b_{_{n,\rmI}} u_{_{1,n}}-d_{_{n,\rmI}}\expmik u_{_{2,n}}-b_{_{n,\rmI}}u_{_{3,n}}=0,\\
				& -c_{_{n,\rmI}}\expik u_{_{1,n+1}}-b_{_{n,\rmI}}u_{_{2,n}}-b_{_{n,\rmI}}u_{_{3,n}}=0.
			\end{aligned}
			\right.	\end{split}
	\end{gather}
	Before investigating the existence of zero eigenvalues, we first notice the intrinsic symmetry of the Hamiltonian operator $ \hatHone(k) $. Let $ \mathbf{u} = \{ \mathbf{u}(n) \}_{n\in\mathbb{Z}}\in l^2_k(\mathbb{Z};\mathbb{C}^6)$, where $ \mathbf{u}(n) = ( u_{_{j,n}} )_{j=1}^6\in \mathbb{C}^6 $ denotes its components. The bounded operators $\widehat{T}(k), \widehat{V}:l^2_k(\mathbb{Z};\mathbb{C}^6)\to l^2_k(\mathbb{Z};\mathbb{C}^6)$ can be defined componentwise as follows
	\begin{equation}
		\begin{split}
			[\widehat{T}(k) \mathbf{u}](n) &\triangleq 
			\eu^{-\iu n k}\begin{pmatrix}
				0_{3\times 3}           & \mathrm{Id}_{3\times 3} \\
				\mathrm{Id}_{3\times 3} & 0_{3\times 3}
			\end{pmatrix}\overline{\mathbf{u}(n)},\\
			(\widehat{V} \mathbf{u})(n) &\triangleq 
			\begin{pmatrix}
				\mathrm{Id}_{3\times 3} & 0_{3\times 3}            \\
				0_{3\times 3}           & -\mathrm{Id}_{3\times 3}
			\end{pmatrix}\mathbf{u}(n).
		\end{split}
	\end{equation}

	Geometrically, $\widehat{V}$ represents the sublattice symmetry of the Hamiltonian operators, while $\widehat{T}(k)$ is the antiunitary reflection operator with respect to the blue axis, as shown in the left panel of \Cref{fig:CellExample}.
	
	By direct calculation one has 
	\begin{proposition}
		For any $ k\in [-\pi,\pi) $ and $  \mathbf{u} \in l_k^2(\mathbb{Z};\mathbb{C}^6) $ satisfying $ \hatHone(k)\mathbf{u} = E_{\rmI}(k)\mathbf{u} $ for some $ E_{\rmI}(k)\in\mathbb{R} $, one has
		\begin{equation}
			\begin{aligned}
				\hatHone(k) \widehat{T}(k) \mathbf{u} = \widehat{T}(k)\hatHone(k) \mathbf{u} = E_{\rmI}(k) \widehat{T}(k)\mathbf{u},\\
				\widehat T(k)^2=\operatorname{Id},\quad \hatHone(k) \widehat{V} \mathbf{u} = -E_{\rmI}(k) \widehat{V} \mathbf{u}.
			\end{aligned}
		\end{equation}
	\end{proposition}
	
	When $E_{\rmI}(k)=0$, the two sets of equations in \eqref{eqn:ZeroIEdge} are decoupled. Therefore, to construct a zero-energy state, it is sufficient to first seek a nonzero solution that takes the form
	\begin{equation}
		\mathbf{u}(n) = (0,0,0,u_{_{4,n}},u_{_{5,n}},u_{_{6,n}})^{T}.
	\end{equation}
	Then $ \widehat{T}(k)\mathbf{u} $ is another nonzero solution, and $ \mathbf{u} $ and $ \widehat{T}(k)\mathbf{u} $ are linearly independent. It follows that:
	\begin{proposition}\label{prop:TwoFold}
		For any $ k\in [-\pi,\pi) $, the eigenspace associated with the zero eigenvalue $ \{ \mathbf{u}\in l^2_k(\mathbb{Z};\mathbb{C}^6): \hatHone(k) \mathbf{u}=0 \} $ is even dimensional. 
	\end{proposition}\par 
	%We will prove the following theorem characterizing the existence of zero eigenvalue at $ k=0 $.
	%\begin{theorem}  % This theorem should be 
	%	Suppose the hopping coefficients of the Hamiltonian operator $ \hatHone(k) $ are given by \eqref{eqn:CoefI}. Then there exists a two-fold zero eigenvalue at $ E_{\rmI}(0)=0 $ if and only if the relation \eqref{eqn:MatchingCond} holds.
	%\end{theorem}
	
	\subsection{Asymptotic Behavior at Infinity}
	From the intrinsic symmetry, we can first investigate the existence of nonzero solutions to the following equations for $ n\in \mathbb{Z} $ and $ k\in [-\pi,\pi) $, 
	\begin{gather}\label{eqn:DecoupleEqn}
		\left\{\begin{aligned}
			& -b_{_{n,\rmI}}u_{_{4,n}}-b_{_{n,\rmI}}u_{_{5,n}}-c_{_{n-1,\rmI}}\expmik u_{_{6,n-1}}=0,\\
			& -b_{_{n,\rmI}} u_{_{4,n}}-d_{_{n,\rmI}}\expik u_{_{5,n}}-b_{_{n,\rmI}}u_{_{6,n}}=0,\\
			& -c_{_{n,\rmI}}u_{_{4,n+1}} -b_{_{n,\rmI}}u_{_{5,n}}-b_{_{n,\rmI}} u_{_{6,n}}=0,\\
		\end{aligned}\right.
	\end{gather}
	such that the corresponding eigenfunction $ \mathbf{u} = \{ (u_{_{j,n}})_{j=1}^{6} \}_{n\in \mathbb{Z}}\in l^2_k(\mathbb{Z};\mathbb{C}^6) $. \par 
	We first define several matrices $ \{ A_j(b,\varepsilon,k) \}_{j=1}^{6} $ to describe the asymptotic behavior at infinity. 
	\begin{gather}\label{eqn:DefPropagationMat}
		\begin{split}
			A_1 \triangleq \begin{pmatrix}
				-b & 0\\
				-b  & -(b+\varepsilon)\expmik
			\end{pmatrix},\quad
			A_2 \triangleq \begin{pmatrix}
				-b &-(b+\varepsilon)\expik \\
				0 & -b
			\end{pmatrix},\\ 
			A_3 \triangleq \begin{pmatrix}
				-(b+\varepsilon)\expmik & 0\\
				-b & -b
			\end{pmatrix},\quad 
			A_4 \triangleq \begin{pmatrix}
				-b & -b\\
				0 & -(b+\varepsilon)
			\end{pmatrix},\\ 
			A_5 \triangleq \begin{pmatrix}
				-b & 0\\
				-(b+\varepsilon)\expik	 & -b
			\end{pmatrix},\quad 
			A_6 \triangleq \begin{pmatrix}
				-(b+\varepsilon) & -b\\
				0 & -b
			\end{pmatrix}.
		\end{split}
	\end{gather}
	These matrices $ \{A_{j}(b,\varepsilon,k)\}_{j=1}^6 $ depend on real parameters $ b, \varepsilon, k $ with $ b>0 $ and $ b+\varepsilon>0 $. 
	We further define $ \widetilde{A}_{1}(k),\widetilde{A}_{6} $
	\begin{equation}
		\widetilde{A}_1(k)\triangleq \begin{pmatrix}
			-b_{+}  & 0 \\
			-b_{+}  & -c\expmik 
		\end{pmatrix},\quad \widetilde{A}_6 \triangleq \begin{pmatrix}
			-c & -b_{-}\\
			0 & -b_{-} 
		\end{pmatrix}.
	\end{equation} 
	From these matrices, we can characterize the behavior as $ n\to \pm\infty $ by rewriting the above equations \eqref{eqn:DecoupleEqn}, which gives the following lemma
	\begin{lemma}\label{lem:PropagationMat}
		For a nonzero solution $\mathbf{u}$ to \eqref{eqn:DecoupleEqn}, the matrices $ -A_2^{-1}A_1,-A_4^{-1}A_3 $ and $ -A_6^{-1}A_5 $ satisfy
		\begin{gather}
			\begin{split}
				&-A_2^{-1}\widetilde{A}_1(k)(
				u_{_{4,0}},u_{_{6,-1}})^{T} = (
				u_{_{6,0}},u_{_{5,0}}
				)^{T},\\
				&-A_2^{-1}A_1(u_{_{4,2m}},u_{_{6,2m-1}})^{T} = (u_{_{6,2m}},u_{_{5,2m}})^{T} ,\quad m\ge 1,\\
				&-A_4^{-1}A_3(
				u_{_{6,2m}},
				u_{_{5,2m}}
				)^{T} = (
				u_{_{5,2m+1}},
				u_{_{4,2m+1}}
				)^{T}, \quad m \ge 0,\\
				&-A_6^{-1}A_5(
				u_{_{5,2m+1}},
				u_{_{4,2m+1}}
				)^{T} = (
				u_{_{4,2m+2}},
				u_{_{6,2m+1}}
				)^{T}, \quad m \ge 0.
			\end{split}
		\end{gather}
		when $ (b,\varepsilon,k) = (b_{+},\delta_{+},k) $. A similar result holds for $ (b,\varepsilon,k) = (b_{-},\delta_{-},k) $ 
		\begin{gather}
			\begin{split}
				&-A_5^{-1}\widetilde{A}_6
				(u_{_{4,0}},u_{_{6,-1}})^{T} = (u_{_{5,-1}},u_{_{4,-1}})^{T},\\
				&-A_5^{-1}A_6(
				u_{_{4,2m}},
				u_{_{6,2m-1}}
				)^{T} = (
				u_{_{5,2m-1}},
				u_{_{4,2m-1}}
				)^{T},\quad m\le -1,\\
				&-A_3^{-1}A_4(u_{_{5,2m-1}},u_{_{4,2m-1}})^{T} 
				= (u_{_{6,2m-2}},u_{_{5,2m-2}})^{T},\quad m\le 0,\\
				&-A_1^{-1}A_2(u_{_{6,2m-2}},u_{_{5,2m-2}})^{T} = 
				(u_{_{4,2m-2}},u_{_{6,2m-3}})^{T}, m \le 0.
			\end{split}
		\end{gather}
		Therefore, by defining the \emph{propagation matrix} $  P(b_+,\delta_{+},k) \triangleq  -A_6^{-1}A_5A_4^{-1}A_3A_2^{-1}A_1$, one has 
		\begin{gather}\label{eqn:DynamicToInf}
			\begin{split}
				&P(b_+,\delta_{+},k)A_1^{-1}\widetilde{A}_1(k)
				(u_{_{4,0}},u_{_{6,-1}})^{T} = (u_{_{4,2}},u_{_{6,1}}	)^{T},\\
				&P(b_+,\delta_{+},k)  
				(u_{_{4,2m}},u_{_{6,2m-1}})^{T} = 
				(u_{_{4,2m+2}},u_{_{6,2m+1}})^{T},\quad m\ge 1,\\
				&P^{-1}(b_{-},\delta_{-},k)A_6^{-1}\widetilde{A}_6
				(u_{_{4,0}},u_{_{6,-1}})^{T} = 
				(u_{_{4,-2}},u_{_{6,-3}})^{T},\\
				&P^{-1}(b_{-},\delta_{-},k) 
				(u_{_{4,2m}},u_{_{6,2m-1}})^{T} = 
				(u_{_{4,2m-2}},	u_{_{6,2m-3}})^{T},\quad m\le -1.			
			\end{split}
		\end{gather}
	\end{lemma}
	To complete the analysis, one has to study the eigenvalues of the matrix $ P(b_{+},\delta_{+},k) $ (respectively, $ P^{-1}(b_{-},\delta_{-},k) $) to characterize the asymptotic behavior at $ +\infty$ (respectively, $-\infty$). By direct calculation, we are able to explicitly derive the eigenvalues of the matrix $ P(b,\varepsilon,k) $.
	\begin{theorem}\label{thm:PropagationMat}
		And the matrix $ P(b,\varepsilon,k) $ has two distinct eigenvalues $ \lambda_1(b,\varepsilon,k)$, $\lambda_2(b,\varepsilon,k) $ when $ k\neq 0 $ or $\varepsilon\neq 0$. Moreover, when $ k\neq 0 $ or $\varepsilon\neq 0$, we have 
		\begin{equation}\label{eqn:Eigenvalues}
			0<|\lambda_1|<1,\quad |\lambda_2|>1.
		\end{equation}
	\end{theorem}
	\begin{proof}
		We begin by carrying out explicit calculations:
		\begin{gather*}
			\begin{split}
				A_2^{-1}A_1 = \frac{1}{b^2}\begin{pmatrix}
					b^2-b(b+\varepsilon)\expik & -(b+\varepsilon)^2\\
					b^2 & b(b+\varepsilon)\expmik
				\end{pmatrix}\\
				A_4^{-1}A_3 = \frac{1}{b(b+\varepsilon)}\begin{pmatrix}
					(b+\varepsilon)^2\expmik - b^2 & -b^2\\
					b^2 & b^2
				\end{pmatrix},\\
				A_6^{-1}A_5 = \frac{1}{b(b+\varepsilon)}\begin{pmatrix}
					b^2-b(b+\varepsilon)\expik & -b^2\\
					(b+\varepsilon)^2\eu^{\iu k} & b(b+\varepsilon)
				\end{pmatrix},
			\end{split}
		\end{gather*}
		From this, one can easily calculate the determinant of $ P(b,\varepsilon,k) $
		\[ \lambda_1\lambda_2=\det[P(b,\varepsilon,k)]  = \eu^{-2\iu k},  \]
		As for its trace, a direct but lengthy calculation yields
		\begin{equation}\label{eqn:MatrixForm}
			P(b,\varepsilon,k) = \begin{pmatrix}
				\alpha(b,\varepsilon,k) & \beta(b,\varepsilon,k)\\
				-\expik\beta(b,\varepsilon,k) & \gamma(b,\varepsilon,k)
			\end{pmatrix},
		\end{equation} 
		where the elements are given by 
		\begin{gather}\label{eqn:MatrixElement}
			\begin{split}
				\alpha(b,\varepsilon,k) &= -\expik\frac{(b+\varepsilon)^2}{b^2} + 2\frac{b+\varepsilon}{b}-\expmik + \eu^{2\iu k} -4\expik\frac{b}{b+\varepsilon}+ 4\frac{b^2}{(b+\varepsilon)^2},\\
				\beta(b,\varepsilon,k) &= -\frac{(b+\varepsilon)^3}{b^3} + \expmik\frac{(b+\varepsilon)^2}{b^2} + \expik\frac{b+\varepsilon}{b}-3+2\expmik\frac{b}{b+\varepsilon},\\
				\gamma(b,\varepsilon,k) &= \frac{(b+\varepsilon)^4}{b^4} - \expik\frac{(b+\varepsilon)^2}{b^2} + 2\frac{b+\varepsilon}{b} - \expmik.
			\end{split}
		\end{gather}
		Therefore, we have
		\begin{equation*}
			\operatorname{trace}[P(b,\varepsilon,k)] =\Big( \frac{2b}{b+\varepsilon}+\frac{(b+\varepsilon)^2}{b^2} - \expik \Big)^2 - 2\eu^{-\iu k}.
		\end{equation*}
		A minimizer $k^\ast$ of the real part of $ \operatorname{trace}(P(b,\varepsilon,k)) $ with respect to $k$ satisfies
		\[ \sin k^{\ast}\bigg(2+\frac{4b}{b+\varepsilon} + \frac{2(b+\varepsilon)^2}{b^2}\bigg) -2\sin2k^\ast = 0.  \]
		The above relation holds if and only if $ k^\ast=0 $ or $ k^\ast=\pi $. It is easy to check that the minimizer is $ k^\ast = 0 $. Now we have 
		\[\Re \{\operatorname{trace}[P(b,\varepsilon,k)]\}\ge\Big( \frac{2b}{b+\varepsilon}+\frac{(b+\varepsilon)^2}{b^2} - 1 \Big)^2 - 2.\]
		Minimizing over $ \varepsilon $, we conclude that the minimizer can be achieved if and only if $ \varepsilon = 0 $. Therefore, we have  
		\[ |\lambda_1 +\lambda_2|\ge \Re(\operatorname{trace}(P))>2. \]
		However, if $ \lambda_1 = \eu^{\iu \theta} $, then $ \lambda_2 = \eu^{-\iu(2k+\theta)} $ and one has
		\[ |\lambda_1+\lambda_2| = |\eu^{\iu \theta} + \eu^{-\iu(2k+\theta)}|\le 2. \]
		This contradiction finishes the proof.		
	\end{proof}
	Denote the eigenvector $ v_i(b_{+},\delta_{+},k) $ corresponding to eigenvalue $ \lambda_i $. 
	To ensure that the eigenvector $ \mathbf{u} $ decays at $ \pm \infty $, we impose the following two conditions
	\begin{align}
		A_{1}^{-1}(b_{+},\delta_{+},k)\widetilde{A}_1(k) \widetilde{v}_0&= g_{_1}v_1(b_{+},\delta_{+},k),\label{eqn:DecayPlusInf}\\
		A_{6}^{-1}(b_{-},\delta_{-},k)\widetilde{A}_6\widetilde{v}_0&=g_{_2}v_2(b_{-},\delta_{-},k).\label{eqn:DecayMinusInf}
	\end{align}
	Here $ g_{_1},g_{_2} $ are two complex nonzero constants, while $ \widetilde{v}_0  = (u_{_{4,0}},u_{_{6,-1}})^{T}$. Note that the first condition \eqref{eqn:DecayPlusInf} corresponds to exponential decay as $ n\to +\infty $. By \eqref{eqn:DynamicToInf} and \Cref{thm:PropagationMat}, one has for $ m\ge 1 $, 
	\begin{gather*}
		\begin{split}
			(u_{_{4,2m}},u_{_{6,2m-1}})^{T} &= (P(b_+,\delta_{+},k))^{m}A_{1}^{-1}(b_{+},\delta_{+},k)\widetilde{A}_1(k)
			\widetilde{v}_0\\
			&= g_{_1}[\lambda_1(b_{+},\delta_{+},k)]^{m}v_{1}(b_{+},\delta_{+},k).
		\end{split}
	\end{gather*}
	For $ m\le -1 $, one similarly has 
	\begin{gather*}
		\begin{split}
			(u_{_{4,2m}},u_{_{6,2m-1}})^{T} &= (P(b_-,\delta_{-},k))^{m}A_{6}^{-1}(b_{-},\delta_{-},k)\widetilde{A}_6
			\widetilde{v}_0\\
			&= g_{_2}[\lambda_2(b_{-},\delta_{-},k)]^{m}v_{2}(b_{-},\delta_{-},k).
		\end{split}
	\end{gather*}
	The above procedures give 
	\[ \Vert(u_{_{4,2m}},u_{_{6,2m-1}})^{T}\Vert \le C\rho^{|m|},\quad m\in\mathbb{Z},\]
	for given constants 
	\[ C = \max( |g_{_1}|\Vert v_1 \Vert,|g_{_2}|\Vert v_2 \Vert   ) ,\quad \rho = \max [ |\lambda_1|(b_{+},\delta_{+},k),|\lambda_2|^{-1}(b_{-},\delta_{-},k) ]<1,\] 
	provided that $ \delta_{+},\delta_{-}\neq 0 $ or $ k\neq 0 $. \par 
	To simplify the conditions \eqref{eqn:DecayPlusInf} and \eqref{eqn:DecayMinusInf}, we proceed by direct calculation, 
	\[ v_2(b_{-},\delta_{-},k) = g_{_3}A_{6}^{-1}\widetilde{A}_6\widetilde{v}_0 = g_{_4}A_{6}^{-1}\widetilde{A}_6 \widetilde{A}_{1}^{-1}(k)A_{1} v_1(b_{+},\delta_{+},k), \]
	where we let $ A_1 = A_1(b_{+},\delta_{+},k), A_{6} = A_6(b_{-},\delta_{-},k) $. 
	
	Therefore, we are led to the constraints by properly choosing a complex constant $ g_{_5} $,
	\begin{equation}\label{eqn:Constraints}
		\begin{pmatrix}
			1 & 0\\
			0 & \frac{(b_{+}+\delta_{+})(b_{-}+\delta_{-})}{c^2}
		\end{pmatrix}v_1(b_{+},\delta_{+},k)=g_{_5}v_2(b_{-},\delta_{-},k).
	\end{equation}
	Combining the above results and Proposition \ref{prop:TwoFold}, we have
	\begin{theorem}\label{thm:ZeroEigenvalueI}
		Let $k\in[-\pi,\pi)$ and assume that
		\[
		b_{\pm}>0,\quad
		b_{\pm}+\delta_{\pm}>0,\quad
		c>0,
		\]
		and
		\[
		(k,\delta_{+})\neq(0,0),
		\quad
		(k,\delta_{-})\neq(0,0).
		\] 
		The necessary and sufficient condition for the existence of a two-fold eigenvalue at $ 0 $ for the Hamiltonian operator $ \hatHone(k) $, as defined in \eqref{eqn:OneHamiltonian}, is that the following equality holds
		\begin{equation}\label{eqn:CondZeroEigValue}
			\begin{pmatrix}
				1 & 0\\
				0 & \frac{(b_{+}+\delta_{+})(b_{-}+\delta_{-})}{c^2}
			\end{pmatrix}v_1(b_{+},\delta_{+},k)=g_{_5}v_2(b_{-},\delta_{-},k).
		\end{equation} 
		Here $ v_{i}(b,\varepsilon,k) $ denotes the eigenvector of $ P(b,\varepsilon,k)  $ that corresponds to eigenvalue $ \lambda_{i} $, as defined in \Cref{thm:PropagationMat}.
	\end{theorem}
	
	\subsection{Existence of zero-energy Edge State at $ k=0 $}
	In this section, we discuss the existence of zero eigenvalues when $ k=0 $. Following \eqref{eqn:MatrixElement}, we first characterize the matrix elements of $ P(b,\varepsilon,0) $.
	\begin{lemma}\label{lem:MatElementForm}
		And the matrix elements of $ P(b,\varepsilon,0) $, which are given by \eqref{eqn:MatrixForm} and \eqref{eqn:MatrixElement}, satisfy the following properties:
		\begin{itemize}
			\item  $ \alpha, \beta$ are strictly decreasing with respect to $\varepsilon$;
			\item  $ \gamma,  \gamma-\alpha-2\beta, \gamma-\alpha+2\beta $ and $ \gamma-\alpha $ are strictly increasing with respect to $ \varepsilon $;
			\item $ \alpha+\gamma$ has a unique minimum at $ \varepsilon=0 $. 
		\end{itemize}
		Moreover, $ \alpha(b,0,0) = \gamma(b,0,0) = 1 $, while $ \beta(b,0,0)=0 $ for all $ b>0 $.
	\end{lemma}
	From \Cref{thm:PropagationMat} and \Cref{thm:ZeroEigenvalueI}, we explicitly calculate the eigenvectors of the matrix $ P(b,\varepsilon,0) $.
	\begin{theorem}\label{thm:PMatrixEigV}
		When $ k=0 $ and $ \varepsilon\neq 0 $, the matrix $ P(b,\varepsilon,0) $ has two eigenvalues $ \lambda_1(b,\varepsilon,0),\lambda_2(b,\varepsilon,0) $ such that $ 0<\lambda_1<1 $ and $ \lambda_2 = \lambda_1^{-1} >1 $. Their corresponding eigenvectors are given by 
		\begin{equation}\label{eqn:EigenvectorForm}
			v_1 = (1,f_1(b,\varepsilon,0))^{T},\quad v_2 = (1,f_1^{-1}(b,\varepsilon,0))^{T}.
		\end{equation}
		Here 
		\begin{gather}\label{eqn:EigenvecCoeff}
			f_1(b,\varepsilon,0)=\frac{-\alpha+\gamma-\sqrt{ (\alpha-\gamma)^2-4\beta^2 }}{2\beta}
			\left\{\begin{split}
				&<-1,\quad  \varepsilon<0, \\
				&\in (-1,0),\quad \varepsilon>0.
			\end{split}\right.
		\end{gather}
		The parameters $ \alpha,\beta,\gamma $ are given in \eqref{eqn:MatrixElement} for fixed $ (b,\varepsilon,0) $.
	\end{theorem}
	\begin{proof}
		Notice that the matrix $ P(b,\varepsilon,0) $ is given explicitly by
		\begin{equation}
			P(b,\varepsilon,0) = \begin{pmatrix}
				\alpha(b,\varepsilon,0) & \beta(b,\varepsilon,0) \\
				-\beta(b,\varepsilon,0) & \gamma(b,\varepsilon,0)
			\end{pmatrix}.
		\end{equation}
		One concludes that the eigenvector can be represented by $ v_1 = (1,f_1)^T, v_2 = (1,f_2)^T $, where $ f_1,f_2 $ are nonzero since $ \beta(b,\varepsilon,0) \neq 0 $ for $ \varepsilon\neq 0 $ by \Cref{lem:MatElementForm}. From the definition $ Pv_j = \lambda_j v_j $, one has
		\begin{gather*}
			\left\{\begin{split}
				\alpha+\beta f_j &= \lambda_j \\
				-\beta+\gamma f_j &= \lambda_j f_j
			\end{split}\right.,\quad j = 1,2.
		\end{gather*}
		By simple calculation, one can check that $ f_1,f_2 $ satisfy the following quadratic equation
		\[ \beta f_j^2 + (\alpha-\gamma)f_j+\beta = 0,\quad j = 1,2. \]
		Here the equation has two distinct real roots since the determinant $ (\alpha-\gamma)^2-4\beta^2> 0 $ by \Cref{lem:MatElementForm}. It is also obvious that $ f_2 = f_1^{-1} $.
		As a result, $ f_1,f_2 $ are given by
		\[ \frac{-\alpha+\gamma\pm\sqrt{ (\alpha-\gamma)^2-4\beta^2 }}{2\beta}. \]
		The corresponding eigenvalues $ \lambda_1,\lambda_{2} $ are given by 
		\[ \frac{\alpha+\gamma \pm \sqrt{(\alpha-\gamma)^2-4\beta^2}}{2}. \]\par 
		To identify the eigenvalue $ \lambda_1 $ satisfying $ |\lambda_{1}|\in (0,1) $ and associated coefficient $f_1$, \Cref{lem:MatElementForm} gives
		\[ 0<\frac{\alpha+\gamma - \sqrt{(\alpha-\gamma)^2-4\beta^2}}{2}<1,\quad \frac{\alpha+\gamma +\sqrt{(\alpha-\gamma)^2-4\beta^2}}{2}>1. \]
		Therefore 
		\[ f_1(b,\varepsilon,0) = \frac{-\alpha+\gamma-\sqrt{ (\alpha-\gamma)^2-4\beta^2 }}{2\beta}.  \]
		From \Cref{lem:MatElementForm}, when $ \varepsilon<0 $,
		\[ \alpha-\gamma+\sqrt{ (\alpha-\gamma)^2-4\beta^2}>2\beta,  \]
		while for $ \varepsilon>0 $, 
		\[ 2\beta<\alpha-\gamma
		+\sqrt{ (\alpha-\gamma)^2-4\beta^2 }<0. \] 
		As a result, one has
		\begin{gather*}
			f_1(b,\varepsilon,0)=\frac{-\alpha+\gamma-\sqrt{ (\alpha-\gamma)^2-4\beta^2 }}{2\beta}
			\left\{\begin{split}
				&<-1,\quad  \varepsilon<0, \\
				&\in (-1,0),\quad \varepsilon>0.
			\end{split}\right.
		\end{gather*}
	\end{proof}
	This fact means one can always choose $ c = c^{\ast}_{\rmI}>0 $ such that
	\begin{equation}
		\frac{(b_{+}+\delta_{+})(b_{-}+\delta_{-})}{c^2}f_{1}(b_{+},\delta_{+},0) = f_{2}(b_{-},\delta_{-},0) = f_{1}^{-1}(b_{-},\delta_{-},0), 
	\end{equation}
	since $ f_{1}(b,\varepsilon,0),f_{2}(b,\varepsilon,0) $ are always negative for $ \varepsilon\neq 0 $. It is obvious that $ c^{\ast}_{\rmI} $ is unique, if we suppose the hopping coefficients are positive. Combining the above arguments, we have
	\begin{theorem}\label{thm:MatchingZero}
		For $ k=0 $ and $ \delta_{+},\delta_{-}\neq 0 $, the necessary and sufficient condition for the existence of a two-fold eigenvalue for the Hamiltonian operator $ \hatHone(0) $, as defined in \eqref{eqn:OneHamiltonian}, is that the hopping coefficients satisfy the following condition
		\begin{equation}\label{eqn:MatchingCond}
			(b_{+}+\delta_{+})(b_{-}+\delta_{-})f_{1}(b_{+},\delta_{+},0) f_{1}(b_{-},\delta_{-},0)= c^2.
		\end{equation}
		Here $ f_1 $ is given in \eqref{eqn:EigenvecCoeff} and \eqref{eqn:MatrixElement}. 
	\end{theorem}
	\begin{remark}\label{rmk:TypeI}
		From the above existence result, one does not need to combine two topologically distinct materials to guarantee the existence of zero eigenvalues. Further, zero-energy edge states exist only when the hopping coefficients satisfy specific conditions, which means that they are not robust against perturbations. Nevertheless, they can be used to build propagating states. For example, choose $b_{+}=b_{-}=60$ and $\delta_{+}=\delta_{-}=30$ so that the two bulk materials are identical. The condition in \Cref{thm:MatchingZero} is satisfied when
		\[
		c=\frac{15\times(23-\sqrt{385})}{2}.
		\]
		Thus the type-I interface can support zero-energy edge states even in the absence of a topological transition.
	\end{remark}
	From Proposition \ref{prop:TwoFold} and previous discussions, we know that for the type-I interface, the eigenvectors corresponding to the zero eigenvalue at $ k=0 $ can be properly chosen such that 
	\begin{equation}\label{eqn:TwoFoldEigenvec}
		\begin{split}
			\mathbf{u}^{A} &= \{( 0,0,0,u_{_{4,n}},u_{_{5,n}},u_{_{6,n}} )^{T}\}_{n\in \mathbb{Z}},\\ 
			\mathbf{u}^{B} &=\widehat{T}(0)\mathbf{u}^{A} =  \{( u_{_{4,n}},u_{_{5,n}},u_{_{6,n}},0,0,0 )^{T}\}_{n\in \mathbb{Z}}.
		\end{split}
	\end{equation}
	In the following discussions, we suppose $ \mathbf{u}^A,\mathbf{u}^B $ have real components and are normalized so that $ ( \mathbf{u}^A,\mathbf{u}^A)= (\mathbf{u}^B,\mathbf{u}^B)=1$.
	Although $\widehat T(0)$ relates the two zero-energy modes, these
		modes do not form a Kramers pair, since $ \widehat T(0)^2\neq -\operatorname{Id}$. Consequently, this symmetry alone does not imply a Kramers-type prohibition of backscattering.
	\subsection{Local Behavior of $ E_{\rmI}(k) $ near $ k=0 $}
	In this section, we study the local behavior of $ E_{\rmI}(k) $ given the existence of zero eigenvalues at $ k=0 $.
	\begin{remark}\label{rmk:spectral_isolation}
	The two-fold zero eigenvalue of $\widehat H_{\rm I}(0)$ is isolated from the remaining spectrum. By Weyl's sequence argument, it can be verified that 
	\[
	\bigcup_{k\in[-\pi,\pi)}\sigma_{\rm ess}\bigl(\widehat H_{\rm I}(k)\bigr)
	=
	\sigma_{\rm ess}\bigl(\widehat H_{\mathrm{bulk}}(b_{+},\delta_{+})\bigr)
	\cup
	\sigma_{\rm ess}\bigl(\widehat H_{\mathrm{bulk}}(b_{-},\delta_{-})\bigr).
	\]
	Therefore,
	Proposition \ref{approp:bulk_gap} implies
	\[
		\sigma_{\rm ess}\bigl(\widehat H_{\rm I}(0)\bigr)
		\cap(-g_{\rm ess},g_{\rm ess})
		=
		\emptyset,
		\quad
		g_{\rm ess}:=\min\{|\delta_+|,|\delta_-|\}>0.	\]
	Moreover, zero cannot be an accumulation point of nonzero eigenvalues.
	Otherwise, there would exist distinct eigenvalues
	$E_j\neq0$, $E_j\to0$, and normalized eigenvectors $\mathbf{u}_j$ satisfying
	\[
	\widehat H_{\rmI}(0)\mathbf{u}_j=E_j\mathbf{u}_j.
	\]
	Since $\widehat H_{\rm I}(0)$ is self-adjoint, the vectors $\{\mathbf{u}_j\}_{j=1}^{\infty}$ are orthonormal and hence $\mathbf{u}_j\rightharpoonup0$. Furthermore,
	\[
		\bigl\|\widehat H_{\rm I}(0)\mathbf{u}_j\bigr\|
		=
		|E_j|\longrightarrow0.
	\]
	Thus $\{\mathbf{u}_j\}_{j=1}^{\infty}$ is a Weyl sequence for zero, which would imply $0\in\sigma_{\rm ess}(\widehat H_{\rm I}(0))$, contradicting the essential spectral gap above. Consequently,
	\[
	\operatorname{dist}
			\left(
			0,\,
			\sigma\bigl(\widehat H_{\rm I}(0)\bigr)\setminus\{0\}
			\right)>0.
	\]
	
	From this we can choose a contour $\Gamma$ in $\mathbb{C}$ around $0$ that contains no other point of $\sigma(\hatHone(0))$. Since $\hatHone(k)$ is analytic with respect to $k$, define 
	\[
	P(k)=\frac{1}{2\pi\iu}\oint_\Gamma(z-\hatHone(k))^{-1}\du z.
	\]
	For small $|k|$, $P(k)$ is analytic and has rank two. This justifies the perturbation arguments in the proof of \Cref{thm:CrossingCurveI}.
	\end{remark}
		
	We will prove the following theorem.
	\begin{theorem}\label{thm:CrossingCurveI}
		Suppose the parameters $ ( b_{\pm},\delta_{\pm},c) $ are properly chosen such that for the Hamiltonian operator $ \hatHone(k) $, there exists a two-fold zero eigenvalue at $ k=0 $. Then there exist two distinct point spectra $ E_{\rmI,-}(k),E_{\rmI,+}(k) $ near zero that admit the following asymptotic expansion.
		\begin{equation}
			E_{\rmI,\pm}(k) = \pm E^{(1)}_{\rmI} | k | + \mathcal{O}(k^2).	
		\end{equation}
		Here $ E^{(1)}_{\rmI} $ is a positive real constant determined from the parameters $ ( b_{\pm},\delta_{\pm},c) $.
	\end{theorem}
	\begin{proof}
		
		For the type-I interface, one has by direct Taylor expansion
		\[ \hatHone(k) = \hatHone(0)+k\widehat{H}^{(1)}_{\mathrm{I}}+\mathcal{O}(k^2). \]
		Here the first order part $ \widehat{H}^{(1)}_{\mathrm{I}}:l^2_k(\mathbb{Z};\mathbb{C}^6)\to l^2_k(\mathbb{Z};\mathbb{C}^6)$ is a bounded operator given by 
		\begin{gather}
			\left\{\begin{split}
				(\widehat{H}^{(1)}_{\rmI}\mathbf{u})_{_{1,n}}&=\iu c_{_{n-1,\rmI}}u_{_{6,n-1}},\\
				(\widehat{H}^{(1)}_{\rmI}\mathbf{u})_{_{2,n}}&=-\iu d_{_{n,\rmI}}u_{_{5,n}},\\
				(\widehat{H}^{(1)}_{\rmI}\mathbf{u})_{_{3,n}}&=0,
			\end{split}\right.\quad 
			\left\{
			\begin{split}
				(\widehat{H}^{(1)}_{\rmI}\mathbf{u})_{_{4,n}}&=0,\\
				(\widehat{H}^{(1)}_{\rmI}\mathbf{u})_{_{5,n}}&=\iu d_{_{n,\rmI}}u_{_{2,n}},\\
				(\widehat{H}^{(1)}_{\rmI}\mathbf{u})_{_{6,n}}&=-\iu c_{_{n,\rmI}}u_{_{1,n+1}}.
			\end{split}
			\right.
		\end{gather}
		Therefore, we suppose that the eigenvalue problem near $ k=0 $
		\begin{equation}\label{eqn:Eigval_Neark}
			\hatHone(k)\mathbf{u}_{k} = E_{\rmI,\pm}(k)\mathbf{u}_{k},
		\end{equation}
		admits the following asymptotic expansions
		\begin{equation}
			\begin{split}
				\hatHone(k) &= \hatHone(0)+k\hatHone^{(1)}+k^2\widehat{H}_{\rmI}^{(r)},\\ E_{\rmI}(k) &= kE^{(1)}_{\rmI,\pm}+\mathcal{O}(k^2),\\ 
				\mathbf{u}_{k} &= a_1\mathbf{u}^{A}+a_2\mathbf{u}^{B}+k\mathbf{u}^{(1)}_{\rmI}+\mathcal{O}(k^2).
			\end{split}
		\end{equation}
		Here $ a_1,a_2 $ are complex constants. By matching the first order terms in \eqref{eqn:Eigval_Neark}, we have 
		\[ \widehat{H}^{(1)}_{\rmI}(a_1\mathbf{u}^{A}+a_2\mathbf{u}^{B})+\hatHone(0)\mathbf{u}^{(1)}=E_{\rmI}^{(1)}(a_1\mathbf{u}^{A}+a_2\mathbf{u}^{B}),\]
		Taking the inner product with $ \mathbf{u}^{A},\mathbf{u}^{B} $ on both sides, one has
		\begin{equation}
			\left\{\begin{split}
				(  \mathbf{u}^{A},\widehat{H}^{(1)}_{\rmI}(a_1\mathbf{u}^{A}+a_2\mathbf{u}^{B}) )+ ( \mathbf{u}^{A},\hatHone(0)\mathbf{u}_{\rmI}^{(1)} ) &= ( \mathbf{u}^{A}, E_{\rmI}^{(1)}(a_1\mathbf{u}^{A}+a_2\mathbf{u}^{B}) ),\\
				(  \mathbf{u}^{B},\widehat{H}^{(1)}_{\rmI}(a_1\mathbf{u}^{A}+a_2\mathbf{u}^{B}) )+ ( \mathbf{u}^{B},\hatHone(0)\mathbf{u}_{\rmI}^{(1)} ) &= ( \mathbf{u}^{B}, E_{\rmI}^{(1)}(a_1\mathbf{u}^{A}+a_2\mathbf{u}^{B}) ).
			\end{split}\right.
		\end{equation}
		Since the Hamiltonian operator $ \hatHone(0) $ is self-adjoint, we are led to
		\begin{gather*}
			\begin{aligned}
				( \mathbf{u}^{A},\hatHone(0)\mathbf{u}_{\rmI}^{(1)} )&=( \hatHone(0)\mathbf{u}^{A},\mathbf{u}_{\rmI}^{(1)} )=0, \\
				( \mathbf{u}^{B},\hatHone(0)\mathbf{u}_{\rmI}^{(1)} )&= ( \hatHone(0)\mathbf{u}^{B},\mathbf{u}_{\rmI}^{(1)} )=0.
			\end{aligned}
		\end{gather*}
		Therefore, we conclude that the values $ E_{\rmI}^{(1)} $ are given by the eigenvalues of the matrix 
		\begin{equation}
			\mathbf{M}_{0,\rmI} \triangleq \begin{pmatrix}
				( \mathbf{u}^{A},\hatHone^{(1)}\mathbf{u}^{A})  &(\mathbf{u}^{A},\hatHone^{(1)}\mathbf{u}^{B}) \\ 
				(\mathbf{u}^{B},\hatHone^{(1)}\mathbf{u}^{A})& (\mathbf{u}^{B},\hatHone^{(1)}\mathbf{u}^{B} )
			\end{pmatrix}.
		\end{equation}
		From the definition \eqref{eqn:TwoFoldEigenvec}, we have 
		\begin{equation}
			\begin{split}
				( \mathbf{u}^{A},\hatHone^{(1)}\mathbf{u}^{A}) &= ( \mathbf{u}^{B},\hatHone^{(1)}\mathbf{u}^{B})=0,\\
				( \mathbf{u}^{A},\hatHone^{(1)}\mathbf{u}^{B}) &= \iu \sum_{n\in \mathbb{Z}}(d_{_{n,\rmI}}|u_{_{5,n}}|^2-c_{_{n,\rmI}}u_{_{6,n}}u_{_{4,n+1}}), \\
				( \mathbf{u}^{B},\hatHone^{(1)}\mathbf{u}^{A}) &= \iu \sum_{n\in \mathbb{Z}}(-d_{_{n,\rmI}}|u_{_{5,n}}|^2+c_{_{n-1,\rmI}}u_{_{6,n-1}}u_{_{4,n}})=-( \mathbf{u}^{A},\hatHone^{(1)}\mathbf{u}^{B}).
			\end{split}
		\end{equation}
		Thus, one only needs to calculate $ ( \mathbf{u}^{A},\hatHone^{(1)}\mathbf{u}^{B}) $. We recall \eqref{eqn:DecayPlusInf} and \eqref{eqn:DecayMinusInf}, which yield the following equality
		\begin{equation}\label{eqn:ConditionExist}
			\widetilde{v}_0 = g_{_6}\widetilde{A}_1^{-1}(0)A_1v_1(b_{+},\delta_{+},0)  = g_{_7}\widetilde{A}_6^{-1}A_6v_2(b_{-},\delta_{-},0), 
		\end{equation}
		Here $ g_{_6},g_{_7} $ are two real constants and $ \widetilde{v}_0 = (u_{_{4,0}},u_{_{6,-1}})^{T} $. By direct calculation, one has	
		\[ \widetilde{v}_0 = g_{_6}\begin{pmatrix}
			1 & 0 \\
			0 & \frac{b_{+}+\delta_{+}}{c}
		\end{pmatrix} v_1(b_{+},\delta_{+},0). \]
		From previous results, we have, for $ m\ge 1 $,
		\begin{equation}\label{eqn:VecFormulaEven1}
			(u_{_{4,2m}},u_{_{6,2m-1}})^{T}
			=(P(b_+,\delta_{+},0))^{m}A_{1}^{-1}(b_{+},\delta_{+},0)\widetilde{A}_1(0)
			\widetilde{v}_0 = g_{_6}\lambda_{1}^{m}v_1.
		\end{equation}
		For $ m\le -1 $, we have
		\begin{equation}\label{eqn:VecFormulaEven2}
			(u_{_{4,2m}},u_{_{6,2m-1}})^{T}=(P(b_-,\delta_{-},0))^{m}A_{6}^{-1}(b_{-},\delta_{-},0)\widetilde{A}_6
			\widetilde{v}_0 = g_{_7} \lambda_2^{m}v_2.
		\end{equation}
		Thus we have, from \Cref{thm:PMatrixEigV},
		\[ u_{_{6,2m-1}}u_{_{4,2m}}<0,\quad m\in \mathbb{Z}. \]\par 
		Now we continue to calculate $ u_{_{6,2m}}u_{_{4,2m+1}} $ for $ m\ge 0 $. By \Cref{lem:PropagationMat}, we have
		\begin{gather*}
			\begin{split}
				-A_2^{-1}(b_{+},\delta_{+},0)\widetilde{A}_1(0)\widetilde{v}_0 &= (
				u_{_{6,0}},u_{_{5,0}}
				)^{T},\\
				-A_2^{-1}A_1(b_{+},\delta_{+},0)(u_{_{4,2m}},u_{_{6,2m-1}})^{T} &= (u_{_{6,2m}},u_{_{5,2m}})^{T} ,\quad m\ge 1,\\
				-A_5^{-1}A_6(b_{+},\delta_{+},0)
				(u_{_{4,2m+2}},u_{_{6,2m+1}})^{T} &= (u_{_{5,2m+1}},	u_{_{4,2m+1}})^{T}
				, \quad m \ge 0.\\
			\end{split}
		\end{gather*}
		It is clear from \eqref{eqn:ConditionExist}, \eqref{eqn:VecFormulaEven1} and \eqref{eqn:VecFormulaEven2} that the following relations hold
		\begin{gather*}
			\begin{split}
				(u_{_{6,0}},u_{_{5,0}})^{T}&= -g_{_6}A_2^{-1}A_1v_1(b_{+},\delta_{+},0),\\
				(u_{_{6,2m}},u_{_{5,2m}})^{T} &= -g_{_{6}}A_2^{-1}A_1\lambda_{1}^{m}v_{1}(b_{+},\delta_{+},0) ,\quad m\ge 1,\\
				(u_{_{5,2m+1}},u_{_{4,2m+1}})^{T} &= -g_{_6} A_5^{-1}A_6
				\lambda^{m+1}_{1}v_1(b_{+},\delta_{+},0),\quad m\ge 0.
			\end{split}
		\end{gather*}
		Therefore, by \eqref{eqn:DefPropagationMat}, one has
		\begin{gather*}
			\begin{split}
				u_{_{6,0}} &= g_{_6}\big[ (t_{+}-1)+t_{+}^{2}f_1\big],\\
				u_{_{6,2m}} &= g_{_6}\lambda^{m}_{1}\big[ (t_{+}-1)+t_{+}^2f_1\big],\\
				u_{_{4,2m+1}} &= g_{_6}\lambda_1^{m+1}\big[ t_{+}^{2} +(t_{+}-1)f_1\big].
			\end{split}
		\end{gather*}
		Here $ t_{+} = (b_{+}+\delta_{+})/b_{+} $. It is obvious that 
		\[ u_{_{6,2m}}u_{_{4,2m+1}} = g^2_{_6}\lambda^{2m+1}_{1}\big[ (t_{+}-1)+t_{+}^2f_1\big]\big[ t_{+}^{2} +(t_{+}-1)f_1\big].\]
		From \Cref{thm:PMatrixEigV}, one has 
		\[ f_1+\frac{1}{f_1} = \frac{-\alpha+\gamma}{\beta} = -\frac{t_{+}^2+t_{+}+2}{t_{+}}. \]
		By tedious calculations one has
		\[ u_{_{6,2m}}u_{_{4,2m+1}} = g^2_{_6}\lambda^{2m+1}_{1}f_{1}<0.\]
		Similar procedures hold for $ m\le -1 $, leading to 
		\[ u_{_{6,n}}u_{_{4,n+1}} <0,\quad n\in\mathbb{Z}. \]
		Therefore, we conclude that 
		\[ ( \mathbf{u}^{A},\hatHone^{(1)}\mathbf{u}^{B})\neq 0, \]
		which means that the matrix $ \mathbf{M}_{0,\rmI} $ is nonzero. This matrix has two distinct real eigenvalues $ \pm \iu( \mathbf{u}^{A},\hatHone^{(1)}\mathbf{u}^{B}) $, which means that the crossing at $ k=0 $ is non-tangential. Therefore, if the two eigenvalues are arranged in increasing order, we have 
		\[E_{\rmI,\pm} (k)= \pm |( \mathbf{u}^{A},\hatHone^{(1)}\mathbf{u}^{B})|  \,|k|+\mathcal{O}(k^2) \triangleq \pm E_{\rmI}^{(1)}|k|+\mathcal{O}(k^2).\]
	\end{proof}
	
	\begin{remark}
		The above theorem shows that a wave packet concentrated near $E_{\rmI}(0)=0$ will propagate in both directions with nonzero group speed. A rigorous characterization for continuous systems is proved in \cite{Fefferman2013}.
	\end{remark}
	\section{Edge States for Type II Interface}\label{sec:TypeII}
	In this section, we investigate the existence of zero-energy edge states for the type-II interface and the behavior of the associated eigenvalue branches near $ k=0 $. Specifically, we study the existence of nonzero solutions to the following equations, where the coefficients $ \{ b_{_{n,\rmII}},c_{_{n,\rmII}},d_{_{n,\rmII}} \}_{n\in \mathbb{Z}} $ are given by \eqref{eqn:CoefII}:
	\begin{gather}\label{eqn:ZeroIIEdge}
		\begin{split}
			&\left\{
			\begin{aligned}
				& -b_{_{n,\rmII}}w_{_{4,n}}-b_{_{n,\rmII}}w_{_{5,n}}-c_{_{n,\rmII}}\expmik w_{_{6,n+1}}=0,\\
				& -b_{_{n,\rmII}} w_{_{4,n}}-d_{_{n-2,\rmII}}\expik w_{_{5,n-2}}-b_{_{n,\rmII}}w_{_{6,n}}=0,\\
				&-c_{_{n,\rmII}}w_{_{4,n+1}} -b_{_{n,\rmII}}w_{_{5,n}}-b_{_{n,\rmII}} w_{_{6,n}}=0,
			\end{aligned}\right.\\
			&\left\{
			\begin{aligned}
				& -b_{_{n,\rmII}} w_{_{1,n}}-b_{_{n,\rmII}}w_{_{2,n}}-c_{_{n-1,\rmII}}w_{_{3,n-1}}=0,\\
				& -b_{_{n,\rmII}} w_{_{1,n}}-d_{_{n,\rmII}}\expmik w_{_{2,n+2}}-b_{_{n,\rmII}}w_{_{3,n}}=0,\\
				& -c_{_{n-1,\rmII}}\expik w_{_{1,n-1}}-b_{_{n,\rmII}}w_{_{2,n}}-b_{_{n,\rmII}}w_{_{3,n}}=0.
			\end{aligned}
			\right.\end{split}
	\end{gather}
	We will discuss two situations in this section: connecting topologically distinct materials and connecting topologically identical materials. In the former situation, one has $ \delta_{+}\delta_{-}<0 $, while in the latter situation, one has $ \delta_{+}\delta_{-}>0 $. \par 
	First and foremost, we also notice the intrinsic symmetry of the Hamiltonian operator $ \widehat{H}_{\rmII} $. Let $ \mathbf{w} = \{ \mathbf{w}(n) \}_{n\in\mathbb{Z}} $, where $ \mathbf{w}(n) = ( w_{_{j,n}} )_{j=1}^6\in \mathbb{C}^6 $ denotes its components. Define the antiunitary reflection operator $ \widehat{R}(k):l^2_{k}(\mathbb{Z};\mathbb{C}^6)\to l^2_{k}(\mathbb{Z};\mathbb{C}^6) $ as follows 
	\begin{equation}
		[\widehat{R}(k)\mathbf{w}] \triangleq \eu^{\iu n k} 
		\begin{pmatrix}
			P_{3\times3}	& 0_{3\times3}\\
			0_{3\times3}    & P_{3\times3}	
		\end{pmatrix}
		\overline{\mathbf{w}(n)},\quad P_{3\times3} \triangleq
		\begin{pmatrix}
			0 & 0 & 1\\
			0 & 1 & 0\\
			1 & 0 & 0
		\end{pmatrix},
	\end{equation}
	which is an analogue of $\widehat{T}(k)$. 
	Similarly, $\widehat{R}(k)$ represents the reflection symmetry across the blue axis in type-II interface, as shown in the right panel of \Cref{fig:CellExample}. By direct calculation, one also has
	\begin{proposition}
		For any $ k\in [-\pi,\pi) $ and $  \mathbf{w} \in l_k^2(\mathbb{Z};\mathbb{C}^6) $ satisfying $ \hatHtwo(k)\mathbf{w} = E_{\rmII}(k)\mathbf{w} $ for some $ E_{\rmII}(k)\in\mathbb{R} $, one has
		\begin{gather}
			\begin{split}
				\hatHtwo(k) \widehat{R}(k) \mathbf{w} &= \widehat{R}(k)\hatHtwo(k) \mathbf{w} = E_{\rmII}(k) \widehat{R}(k)\mathbf{w},\\
				\widehat{R}(k)^2 &= \operatorname{Id},\;  \hatHtwo(k) \widehat{V} \mathbf{w} = -E_{\rmII}(k)\widehat{V} \mathbf{w}.
			\end{split}
		\end{gather}
	\end{proposition}
	One can also deduce from the above proposition that for $ E_{\rmII}=0 $, one can seek two solutions $ \mathbf{w}^{A} = \{w^A(n)\}_{n\in \mathbb{Z}},\mathbf{w}^{B} = \{w^{B}(n) \}_{n\in \mathbb{Z}} $ that take the form
	\begin{equation}
		w^{A}(n) = ( 0,0,0,w_{_{4,n}},w_{_{5,n}},w_{_{6,n}} )^{T},\quad w^{B}(n) = (w_{_{1,n}},w_{_{2,n}},w_{_{3,n}},0,0,0 )^{T},
	\end{equation}
	respectively. It is worth noting that the eigenspace associated with the zero eigenvalue $ \{ \mathbf{w}\in l^{2}_{k}(\mathbb{Z};\mathbb{C}^6):\hatHtwo(k)\mathbf{w}=0 \} $ need not be even dimensional, which is different from the type-I interface. \par 
	In the following section we seek solutions $ \mathbf{w}^{A},\mathbf{w}^{B} $ that satisfy 
	\begin{equation}\label{eqn:Inversion}
		\widehat{R}(k)\mathbf{w}^{A} = \mathbf{w}^{A},\quad \widehat{R}(k)\mathbf{w}^{B} = \mathbf{w}^{B}.
	\end{equation}
	The reason is that, if $ \widehat{R}(k)\mathbf{w}^{A}  \neq \mathbf{w}^{A} $, then one defines
	\[ \mathbf{w}_{_1}^{A} \triangleq\mathbf{w}^{A}+ \widehat{R}(k)\mathbf{w}^{A},\quad \mathbf{w}_{_2}^{A} \triangleq\mathbf{w}^{A}- \widehat{R}(k)\mathbf{w}^{A}. \]
	For $ \mathbf{w}_{_1}^{A} $, the desired symmetry relation is immediate. For $ \mathbf{w}_{_2}^{A} $, one has $ \widehat{R}(k)[\iu\mathbf{w}_{_2}^{A}] = \iu\mathbf{w}_{_2}^{A}$. Thus, without loss of generality, one supposes the relations \eqref{eqn:Inversion} hold. From these relations, we have 
	\begin{gather*}
		\begin{split}
			w_{_{6,n}} = \overline{w_{_{4,n}}}\eu^{\iu nk},\quad w_{_{5,n}} = \overline{w_{_{5,n}}}\eu^{\iu nk};\\
			w_{_{3,n}} = \overline{w_{_{1,n}}}\eu^{\iu nk},\quad w_{_{2,n}} = \overline{w_{_{2,n}}}\eu^{\iu nk}.
		\end{split}
	\end{gather*} 
	Therefore, one can recast the solution $ \{(w_{_{j,n}})_{j=1}^6\}_{n\in\mathbb{Z}} $ into another set of solutions $ \{(v_{_{j,n}})_{j=1}^6\}_{n\in\mathbb{Z}}  $ 
	\begin{gather}
		\begin{split}
			w_{_{4,n}} = v_{_{4,n}}\ehalf,\quad w_{_{5,n}} = v_{_{5,n}}\ehalf,\quad w_{_{6,n}} = \overline{v_{_{4,n}}}\ehalf,\\
			w_{_{1,n}} = v_{_{1,n}}\ehalf,\quad w_{_{2,n}} = v_{_{2,n}}\ehalf,\quad w_{_{3,n}} = \overline{v_{_{1,n}}}\ehalf,
		\end{split}
	\end{gather}
	where $ v_{_{1,n}},v_{_{4,n}}\in \mathbb{C} $ and $ v_{_{2,n}},v_{_{5,n}}\in \mathbb{R} $ for all $ n\in \mathbb{Z} $. The original equations \eqref{eqn:ZeroIIEdge} can be equivalently rewritten for $ \{(v_{_{j,n}})_{j=1}^6\}_{n\in\mathbb{Z}}  $ as 
	\begin{gather}\label{eqn:ZeroIIEdgeR}
		\begin{split}
			&\left\{
			\begin{aligned}
				& -b_{_{n,\rmII}}v_{_{4,n}}-b_{_{n,\rmII}}v_{_{5,n}}-c_{_{n,\rmII}}\eu^{- \frac{\iu k}{2}} \overline{v_{_{4,n+1}}}=0,\\
				& -b_{_{n,\rmII}} v_{_{4,n}}-d_{_{n-2,\rmII}} v_{_{5,n-2}}-b_{_{n,\rmII}}\overline{v_{_{4,n}}}=0,\\
				&-c_{_{n,\rmII}}\eu^{\frac{\iu k}{2}}v_{_{4,n+1}} -b_{_{n,\rmII}}v_{_{5,n}}-b_{_{n,\rmII}} \overline{v_{_{4,n}}}=0,
			\end{aligned}\right.\\
			& \left\{
			\begin{aligned}
				& -b_{_{n,\rmII}} v_{_{1,n}}-b_{_{n,\rmII}}v_{_{2,n}}-c_{_{n-1,\rmII}}\eu^{- \frac{\iu k}{2}}\overline{v_{_{1,n-1}}}=0,\\
				& -b_{_{n,\rmII}} v_{_{1,n}}-d_{_{n,\rmII}} v_{_{2,n+2}}-b_{_{n,\rmII}}\overline{v_{_{1,n}}}=0,\\
				& -c_{_{n-1,\rmII}}\eu^{\frac{\iu k}{2}} v_{_{1,n-1}}-b_{_{n,\rmII}}v_{_{2,n}}-b_{_{n,\rmII}}\overline{v_{_{1,n}}}=0.
			\end{aligned}
			\right.	\end{split}
	\end{gather}
	Similarly, we denote the above solutions as $ \mathbf{v}^{A} = \{v^A(n)\}_{n\in \mathbb{Z}}$ and $\mathbf{v}^{B} = \{v^{B}(n) \}_{n\in \mathbb{Z}} $, which are given by 
	\begin{equation}
		v^{A}(n) = ( 0,0,0,v_{_{4,n}},v_{_{5,n}},v_{_{6,n}} )^{T},\quad v^{B}(n) = (v_{_{1,n}},v_{_{2,n}},v_{_{3,n}},0,0,0 )^{T}.
	\end{equation}
	We now turn to finding a nonzero solution to \eqref{eqn:ZeroIIEdgeR} in the following section.
	\subsection{Existence of Zero-energy Edge State}
	In this part, we intend to characterize the existence of a nonzero solution $ \mathbf{v} $ for $ k\in [-\pi,\pi) $. We first investigate the existence of a nonzero solution $ \mathbf{v}^{A} $. It satisfies the following equations for $ n\in\mathbb{Z} $. 
	\begin{gather}
		\left\{\begin{split}
			& -b_{_{n,\rmII}}v_{_{4,n}}-b_{_{n,\rmII}}v_{_{5,n}}-c_{_{n,\rmII}}\eu^{-\frac{\iu k}{2}} \overline{v_{_{4,n+1}}}=0,\\
			& -b_{_{n,\rmII}} v_{_{4,n}}-d_{_{n-2,\rmII}} v_{_{5,n-2}}-b_{_{n,\rmII}}\overline{v_{_{4,n}}}=0,\\
			&-c_{_{n,\rmII}}\eu^{\frac{\iu k}{2}}v_{_{4,n+1}} -b_{_{n,\rmII}}v_{_{5,n}}-b_{_{n,\rmII}} \overline{v_{_{4,n}}}=0.
		\end{split}\right.
	\end{gather}
	The above equations can be rewritten as 
	\begin{align}
		& -b_{_{n,\rmII}}(v_{_{4,n}}+v_{_{5,n}})-c_{_{n,\rmII}}\eu^{-\frac{\iu k}{2}} (\overline{v_{_{4,n+1}}}+v_{_{5,n+1}})+c_{_{n,\rmII}}\eu^{-\frac{\iu k}{2}}v_{_{5,n+1}}=0,\\
		& -b_{_{n+2,\rmII}} v_{_{4,n+2}}-d_{_{n,\rmII}} v_{_{5,n}}-b_{_{n+2,\rmII}}\overline{v_{_{4,n+2}}}=0,\label{eqn:eqnn+1}\\
		& -b_{_{n+1,\rmII}}(v_{_{4,n+1}}+v_{_{5,n+1}})-c_{_{n+1,\rmII}}\eu^{-\frac{\iu k}{2}} \overline{v_{_{4,n+2}}}=0.\label{eqn:eqnn+2}
	\end{align}
	Substituting \eqref{eqn:eqnn+2} into \eqref{eqn:eqnn+1}, one has
	\begin{gather*}
		\begin{split}
			& (v_{_{4,n}}+v_{_{5,n}})=-\frac{c_{_{n,\rmII}}}{b_{_{n,\rmII}}}\eu^{-\frac{\iu k}{2}} (\overline{v_{_{4,n+1}}}+v_{_{5,n+1}})+\frac{c_{_{n,\rmII}}}{b_{_{n,\rmII}}}\eu^{-\frac{\iu k}{2}}v_{_{5,n+1}},\\
			&v_{_{5,n}} = \frac{b_{_{n+2,\rmII}}b_{_{n+1,\rmII}}}{c_{_{n+1,\rmII}}d_{_{n,\rmII}}}\eu^{-\frac{\iu k}{2}}(\overline{v_{_{4,n+1}}}+v_{_{5,n+1}}) +\frac{b_{_{n+2,\rmII}}b_{_{n+1,\rmII}}}{c_{_{n+1,\rmII}}d_{_{n,\rmII}}}\eu^{\frac{\iu k}{2}}(v_{_{4,n+1}}+v_{_{5,n+1}}).\\
		\end{split}
	\end{gather*}
	Therefore, by taking complex conjugation, one has
	\begin{gather}
		\begin{split}
			\chi_{_{n+1}}&\triangleq((v_{_{4,n+1}}+v_{_{5,n+1}}),(\overline{v_{_{4,n+1}}}+v_{_{5,n+1}}),v_{_{5,n+1}})^{T} \\
			&= Q_{_{A,n}}^{-1}(k)\chi_{_{n}}.
		\end{split}
	\end{gather}
	The invertible matrix $ Q_{_{A,n}}(k)$ is given as follows:
	\begin{gather}
		Q_{_{A,n}}(k) \triangleq \begin{pmatrix}
			0 &-\frac{c_{_{n,\rmII}}}{b_{_{n,\rmII}}}\eu^{-\frac{\iu k}{2}}  & \frac{ c_{_{n,\rmII}} }{ b_{_{n,\rmII}}}\eu^{-\frac{\iu k}{2}}\\
			-\frac{c_{_{n,\rmII}}}{b_{_{n,\rmII}}}\eu^{\frac{\iu k}{2}} &0  & \frac{c_{_{n,\rmII}} }{b_{_{n,\rmII}}}\eu^{\frac{\iu k}{2}}\\
			\frac{b_{_{n+2,\rmII}}b_{_{n+1,\rmII}}}{c_{_{n+1,\rmII}}d_{_{n,\rmII}}}\eu^{\frac{\iu k}{2}} & \frac{b_{_{n+2,\rmII}}b_{_{n+1,\rmII}}}{c_{_{n+1,\rmII}}d_{_{n,\rmII}}}\eu^{-\frac{\iu k}{2}} & 0
		\end{pmatrix}.
	\end{gather}	
	Similarly for $ \mathbf{v}^B $, one has 
	\begin{gather}
		\begin{split}
			\xi_{_{n+1}}&\triangleq((v_{_{1,n+1}}+v_{_{2,n+1}}),(\overline{v_{_{1,n+1}}}+v_{_{2,n+1}}),v_{_{2,n+1}})^{T} \\
			&= Q_{_{B,n}}(k)\xi_{_{n}}.
		\end{split}
	\end{gather}
	And the matrix $ Q_{_{B,n}}(k) $ is given by
	\begin{gather}
		Q_{_{B,n}}(k) \triangleq \begin{pmatrix}
			0 &-\frac{c_{_{n,\rmII}}}{b_{_{n+1,\rmII}}}\eu^{-\frac{\iu k}{2}}  & \frac{ c_{_{n,\rmII}} }{ b_{_{n+1,\rmII}}}\eu^{-\frac{\iu k}{2}}\\
			-\frac{c_{_{n,\rmII}}}{b_{_{n+1,\rmII}}}\eu^{\frac{\iu k}{2}} &0  & \frac{c_{_{n,\rmII}} }{b_{_{n+1,\rmII}}}\eu^{\frac{\iu k}{2}}\\
			\frac{b_{_{n-1,\rmII}}b_{_{n,\rmII}}}{c_{_{n-1,\rmII}}d_{_{n-1,\rmII}}}\eu^{\frac{\iu k}{2}} & \frac{b_{_{n-1,\rmII}}b_{_{n,\rmII}}}{c_{_{n-1,\rmII}}d_{_{n-1,\rmII}}}\eu^{-\frac{\iu k}{2}} & 0
		\end{pmatrix}.
	\end{gather}
	We further define the invertible matrix $ Q(b,\varepsilon,k) $ for $ b,b+\varepsilon>0 $
	\begin{gather}
		Q(b,\varepsilon,k)\triangleq\begin{pmatrix}
			0 &-\frac{b+\varepsilon}{b}\eu^{-\frac{\iu k}{2}}  & \frac{b+\varepsilon}{b}\eu^{-\frac{\iu k}{2}}\\
			-\frac{b+\varepsilon}{b}\eu^{\frac{\iu k}{2}} &0  & \frac{b+\varepsilon}{b}\eu^{\frac{\iu k}{2}}\\
			\frac{b^2}{(b+\varepsilon)^2}\eu^{\frac{\iu k}{2}} & \frac{b^2}{(b+\varepsilon)^2}\eu^{-\frac{\iu k}{2}} & 0
		\end{pmatrix}.
	\end{gather}
	In summary, by the definition of the hopping coefficients \eqref{eqn:CoefII}, one has
	\begin{theorem}\label{thm:ZeroEigenvalueII}
		Fix $k\in[-\pi,\pi)$. The necessary and sufficient condition for the existence of nonzero solutions $ \mathbf{v}^A$ (respectively, $\mathbf{v}^B$) to the equations \eqref{eqn:ZeroIIEdgeR} is that the following relations \eqref{eqn:CondAII}\emph{(}respectively, \eqref{eqn:CondBII}\emph{)} hold for nonzero $ \{\chi_{_{n}}\}_{n\in \mathbb{Z}} (\text{respectively, } \{\xi_{_{n}}\}_{n\in\mathbb{Z}})$
		\begin{gather}
			\left\{\begin{split}\label{eqn:CondAII}
				\chi_{_{n+1}} &= [Q(b_+,\delta_{+},k)]^{-1}\chi_{_{n}},\quad n\ge 0,\\
				\chi_{_{0}} &= Q_{_{A,-1}}^{-1}(k)\chi_{_{-1}},\qquad \chi_{_{-1}} = Q_{_{A,-2}}^{-1}(k)\chi_{_{-2}}\\
				\chi_{_{n+1}} &= [Q(b_-,\delta_{-},k)]^{-1}\chi_{_{n}},\quad n\le -3.
			\end{split}\right.\\
			\left\{\begin{split}\label{eqn:CondBII}
				\xi_{_{n+1}} &= Q(b_+,\delta_{+},k)\xi_{_{n}},\quad n\ge 1,\\
				\xi_{_{1}} &= Q_{_{B,0}}(k)\xi_{_{0}},\qquad \xi_{_{0}} = Q_{_{B,-1}}(k)\xi_{_{-1}}\\
				\xi_{_{n+1}} &= Q(b_-,\delta_{-},k)\xi_{_{n}},\quad n\le -2.
			\end{split}\right.
		\end{gather}
		Here the matrices $ Q_{_{A,-1}},Q_{_{A,-2}} $ and $ Q_{_{B,0}},Q_{_{B,-1}} $ are given by 
		\begin{equation*}
			Q_{_{A,-1}} \triangleq 
			\begin{pmatrix}
				0 & -\frac{c}{b_{-}}\eu^{-\frac{\iu k}{2}} & \frac{c}{b_{-}}\eu^{-\frac{\iu k}{2}}\\
				-\frac{c}{b_{-}}\eu^{\frac{\iu k}{2}} & 0 & \frac{c}{b_{-}}\eu^{\frac{\iu k}{2}}\\
				\frac{b^{2}_{+}}{(b_{+}+\delta_{+})c}\eu^{\frac{\iu k}{2}} & \frac{b^{2}_{+}}{(b_{+}+\delta_{+})c} \eu^{-\frac{\iu k}{2}}& 0
			\end{pmatrix},
		\end{equation*}
		\begin{equation*}
			Q_{_{A,-2}}\triangleq
			\begin{pmatrix}
				0 & -\frac{b_{-}+\delta_{-}}{b_{-}}\eu^{-\frac{\iu k}{2}} & \frac{b_{-}+\delta_{-}}{b_{-}}\eu^{-\frac{\iu k}{2}}\\
				-\frac{b_{-}+\delta_{-}}{b_{-}}\eu^{\frac{\iu k}{2}} & 0 & \frac{b_{-}+\delta_{-}}{b_{-}}\eu^{\frac{\iu k}{2}}\\
				\frac{b_{+}b_{-}}{c^2}\eu^{\frac{\iu k}{2}} & \frac{b_{+}b_{-}}{c^2}\eu^{-\frac{\iu k}{2}} & 0
			\end{pmatrix},
		\end{equation*}
		\begin{equation*}
			Q_{_{B,0}} \triangleq 
			\begin{pmatrix}
				0 & -\frac{b_{+}+\delta_{+}}{b_{+}} \eu^{-\frac{\iu k}{2}}& \frac{b_{+}+\delta_{+}}{b_{+}}\eu^{-\frac{\iu k}{2}}\\
				-\frac{b_{+}+\delta_{+}}{b_{+}}\eu^{\frac{\iu k}{2}} & 0 & \frac{b_{+}+\delta_{+}}{b_{+}}\eu^{\frac{\iu k}{2}}\\
				\frac{b_{+}b_{-}}{c^2}\eu^{\frac{\iu k}{2}} & \frac{b_{+}b_{-}}{c^2}\eu^{-\frac{\iu k}{2}} & 0
			\end{pmatrix},
		\end{equation*}
		\begin{equation*}
			Q_{_{B,-1}}\triangleq 
			\begin{pmatrix}
				0 & -\frac{c}{b_{+}}\eu^{-\frac{\iu k}{2}} & \frac{c}{b_{+}}\eu^{-\frac{\iu k}{2}}\\
				-\frac{c}{b_{+}}\eu^{\frac{\iu k}{2}} & 0 & \frac{c}{b_{+}}\eu^{\frac{\iu k}{2}}\\
				\frac{b^{2}_{-}}{(b_{-}+\delta_{-})c}\eu^{\frac{\iu k}{2}} & \frac{b^{2}_{-}}{(b_{-}+\delta_{-})c}\eu^{-\frac{\iu k}{2}} & 0
			\end{pmatrix}.
		\end{equation*}
	\end{theorem}
	
	\subsection{Existence of Zero Eigenvalue at $ k=0 $}
	In this section, we analyze the existence of zero eigenvalues by taking advantage of \Cref{thm:ZeroEigenvalueII} at $ k=0 $. 
	To investigate the decay as $n\to  \pm \infty $, one can directly calculate the eigenvalues and eigenvectors of $ Q(b,\varepsilon,0) $.
	\begin{lemma}\label{lem:PropMatEigII}
		For $ b>0,\varepsilon\neq 0 $ satisfying $ b+\varepsilon>0 $, the matrix 
		\[ Q(b,\varepsilon,0)=\begin{pmatrix}
			0 &-\frac{b+\varepsilon}{b}  & \frac{b+\varepsilon}{b}\\
			-\frac{b+\varepsilon}{b} &0  & \frac{b+\varepsilon}{b}\\
			\frac{b^2}{(b+\varepsilon)^2} & \frac{b^2}{(b+\varepsilon)^2} & 0
		\end{pmatrix}, \] has three distinct eigenvalues $ \{ \mu_{j}(b,\varepsilon) \}_{j=1}^{3} $
		\begin{gather}\label{eqn:MatEigenvalueII}
			\left\{\begin{split}
				\mu_{1} &= -\frac{b+\varepsilon}{2b}-\frac{1}{2}\sqrt{\frac{(b+\varepsilon)^2}{b^2}+\frac{8b}{(b+\varepsilon)}}<-2, \\ 
				\mu_{2} &= -\frac{b+\varepsilon}{2b}+\frac{1}{2}\sqrt{\frac{(b+\varepsilon)^2}{b^2}+\frac{8b}{(b+\varepsilon)}}
				\left\{\begin{aligned}
					&>1,\quad \varepsilon<0,\\
					&\in (0,1),\quad \varepsilon>0,
				\end{aligned}\right.\\
				\quad\mu_{3} &= \frac{b+\varepsilon}{b}.
			\end{split}\right.
		\end{gather}
		The corresponding eigenvectors $ \{ v_{j}(b,\varepsilon) \}_{j=1}^{3} $ are 
		\begin{equation}\label{eqn:MatEigenvectorII}
			v_1 = (t_1,t_1,1)^{T},\quad v_{2} = (t_2,t_2,1)^{T},v_{3} = (\iu,-\iu,0)^{T},
		\end{equation}
		where $ t_1(b,\varepsilon),t_2(b,\varepsilon) $ are real numbers given by 
		\begin{equation}\label{eqn:MatEigvecCoefII}
			t_1 =\frac{b+\varepsilon}{\mu_{1}b+b+\varepsilon}= -\frac{b+\varepsilon}{\mu_{2} b},\quad t_2 = \frac{b+\varepsilon}{\mu_{2}b+b+\varepsilon} =-\frac{b+\varepsilon}{\mu_{1} b} .
		\end{equation}
	\end{lemma}
	\begin{remark}
		Although the matrix $ Q(b,\varepsilon,0) $ is real, we let $ v_3 \in \mathbb{C}^3$ to make sure that the original eigenvector $ \mathbf{w}\in l_{0}^{2}(\mathbb{Z};\mathbb{C}^6) $ satisfies $ \widehat{R}(0) w = w $.
	\end{remark}
	We will prove the following theorem characterizing the existence of zero eigenvalues of $ \hatHtwo(k) $ at $ k=0 $. 
	\begin{theorem}\label{thm:ZeroEigenvalueExistenceII}
		Suppose the hopping coefficients of the  Hamiltonian operator $ \hatHtwo(k) $ are given by \eqref{eqn:CoefII}. Then the following statements hold.
		\begin{itemize}
			\item[\textnormal{(i)}] When $ \delta_{+}\delta_{-}<0 $, a two-fold eigenvalue at $ E_{\rmII}(0) = 0 $ at $ k=0 $ always exists.
			\item[\textnormal{(ii)}] When $\delta_{+}\delta_{-}>0$, zero does not belong to the spectrum of $\hatHtwo(0)$.
		\end{itemize}
	\end{theorem}
	\begin{proof}
		We first begin with (i). By \eqref{eqn:CondAII} and \eqref{eqn:CondBII}, if $ \chi_{_{n}},\xi_{_{n}} $ decay as $ n\to \pm \infty $, it should hold that
		\begin{gather}\label{eqn:DecayCondII}
			\begin{split}
				\left\{\begin{aligned}
					\chi_{_{0}} &= h_1 v_1+h_2v_3,\; \chi_{_{-2}} = h_3 v_3,\\
					\xi_{_{1}}&= h_4v_2,\; \xi_{_{-1}} = h_5v_1 + h_6v_2,
				\end{aligned}\right.\quad \delta_{+}>0,\delta_{-}<0,\\
				\left\{\begin{aligned}
					\chi_{_{0}} &= h_1 v_1+h_2v_2,\; \chi_{_{-2}} = h_3 v_2,\\
					\xi_{_{1}}&= h_4v_3,\; \xi_{_{-1}} = h_5v_1 + h_6v_3,
				\end{aligned}\right.\quad \delta_{+}<0,\delta_{-}>0.
			\end{split}
		\end{gather}
		where $ \{h_{j}\}_{j=1}^{6} $ are real constants chosen so that the original eigenvector $ \mathbf{w}\in l_{0}^{2}(\mathbb{Z};\mathbb{C}^6) $ satisfies $ \widehat{R}(0) \mathbf{w} = \mathbf{w} $. Further, it should hold that 
		\begin{equation}\label{eqn:ConnectingCondII}
			Q_{_{A,-2}}Q_{_{A,-1}}(0)\chi_{_{0}} = \chi_{_{-2}},\quad \xi_{_{1}} = Q_{_{B,0}}Q_{_{B,-1}}(0)\xi_{-1}.
		\end{equation}
		By direct calculation, one can write each element of the real matrices. 
		\begin{gather}
			\begin{split}
				Q_{_{A,-2}}Q_{_{A,-1}}(0) = 
				\begin{pmatrix}
					\frac{ct_{-}}{b_{-}}+\frac{t_{-}b_{+}}{ct_{+}} & \frac{t_{-}b_{+}}{ct_{+}} & -\frac{ct_{-}}{b_{-}} \\
					\frac{t_{-}b_{+}}{ct_{+}} & \frac{ct_{-}}{b_{-}}+\frac{t_{-}b_{+}}{ct_{+}} & -\frac{ct_{-}}{b_{-}} \\
					-\frac{b_{+}}{c} & -\frac{b_{+}}{c} & \frac{2b_{+}}{c}
				\end{pmatrix},\\
				Q_{_{B,0}}Q_{_{B,-1}}(0) = 
				\begin{pmatrix}
					\frac{ct_{+}}{b_{+}}+\frac{t_{+}b_{-}}{ct_{-}} & \frac{t_{+ }b_{-}}{ct_{-}} & -\frac{ct_{+}}{b_{+}} \\
					\frac{t_{+}b_{-}}{ct_{-}} & \frac{ct_{+ }}{b_{+}}+\frac{t_{+}b_{-}}{ct_{-}} & -\frac{ct_{+}}{b_{+}} \\
					-\frac{b_{-}}{c} & -\frac{b_{-}}{c} & \frac{2b_{-}}{c}
				\end{pmatrix},
			\end{split}
		\end{gather}
		where $ t_{\pm} \triangleq (b_{\pm}+\delta_{\pm})/b_{\pm} $. For $ \delta_{+}>0,\delta_{-}<0 $, by \eqref{eqn:ConnectingCondII} and the fact that the matrix $ 	Q_{_{A,-2}}Q_{_{A,-1}}(0) $ is real, one has $ h_1 = 0 $. Therefore, supposing $ h_2 = 1 $, one has 
		\[ \chi_{_{-2}} = h_3 v_3(b_{-},\delta_{-}) = Q_{_{A,-2}}Q_{_{A,-1}}(0)v_3(b_{+},\delta_{+}). \]
		It can be easily seen that $ h_3 = ct_{-}/b_{-}>0 $.
		Next we proceed to deal with the equality $\xi_{_{1}} = Q_{_{B,0}}Q_{_{B,-1}}(0)\xi_{-1}  $. Setting $ h_4 = 1 $, a direct calculation gives
		\[ [Q_{_{B,0}}Q_{_{B,-1}}(0)]^{-1}v_2(b_{+},\delta_{+}) = (f_1,f_1,f_2)^{T} =  \xi_{-1} = h_5v_1(b_{-},\delta_{-})+h_6v_2(b_{-},\delta_{-}). \]
		Here $ f_2>f_1>0 $. Moreover, since we have $ t_1(b_{-},\delta_{-})\neq t_2(b_{-},\delta_{-})  $, elementary linear algebra shows that there exist unique real parameters $ h_5,h_6 $ such that 
		\[ h_5t_1(b_{-},\delta_{-})+h_6t_2(b_{-},\delta_{-}) = f_1,\quad h_5+h_6 = f_2. \]
		Thus, there exist unique nonzero vectors $ \{\chi_{_{n}}\}_{n\in\mathbb{Z}} $ and $ \{ \xi_{_{n}} \}_{n\in \mathbb{Z}} $, up to a constant. By \Cref{thm:ZeroEigenvalueII}, the Hamiltonian operator $ \hatHtwo(0) $ admits a two-fold eigenvalue at zero when $ \delta_{+}>0,\delta_{-}<0 $. \par 
		The other sign configuration is analogous, and hence we have proved the existence of nonzero vectors $ \{\chi_{_{n}}\}_{n\in\mathbb{Z}} $ and $ \{ \xi_{_{n}} \}_{n\in \mathbb{Z}} $ when $ \delta_{+}\delta_{-}<0$.\par 
		As for (ii), we argue by contradiction. Still by \eqref{eqn:CondAII} and \eqref{eqn:CondBII}, it should hold that 
		\begin{gather}\label{eqn:ContradictoryCond}
			\begin{split}
				\left\{\begin{aligned}
					\chi_{_{0}} &= h_1v_1+h_2v_3,\;\chi_{_{-2}} = h_3v_2,\\ 
					\xi_{_{1}} &= h_4 v_2,\; \xi_{_{-1}} = h_5v_1+h_6v_3,\\
				\end{aligned}\right.\quad \delta_{+}>0,\delta_{-}>0,\\	
				\left\{\begin{aligned}
					\chi_{_{0}} &= h_1v_1+h_2v_2,\;\chi_{_{-2}} = h_3v_3,\\ 
					\xi_{_{1}} &= h_4 v_3,\; \xi_{_{-1}} = h_5v_1+h_6v_2,
				\end{aligned}\right.\quad \delta_{+}<0,\delta_{-}<0.	
			\end{split}
		\end{gather}
		for real constants $ \{h_{j}\}_{j=1}^{6} $. It can be directly seen that \eqref{eqn:ConnectingCondII} does not hold since the left-hand sides of the equations are complex, while the right-hand sides are real. This contradiction proves the case when $ \delta_{+}<0$ and $ \delta_{-}<0 $. \par 
		When $ \delta_{+}>0 $ and $ \delta_{-}>0 $, one first observes that $ h_2 = h_6 = 0 $. Supposing $ h_1 = h_5 = 1 $, we have directly 
		\[ \chi_{_{0}} = v_1(b_{+},\delta_{+}),\;\chi_{_{-2}} = h_3v_2(b_{-},\delta_{-}),\; \xi_{_{1}} = h_4 v_2(b_{+},\delta_{+}),\; \xi_{_{-1}} = v_1(b_{-},\delta_{-}). \]
		By \eqref{eqn:ConnectingCondII}, one has 
		\[ -2\frac{b_{+}}{c}t_1(b_{+},\delta_{+}) + 2\frac{b_{+}}{c} = h_3,\quad -2\frac{b_{-}}{c}t_1(b_{-},\delta_{-}) + 2\frac{b_{-}}{c} = h_4,  \]
		where $ t_1 $ is given by \eqref{eqn:MatEigvecCoefII}. From \eqref{eqn:MatEigenvalueII}, one has by direct calculation
		\[ t_1(b_{\pm},\delta_{\pm}) = -\frac{t_{\pm}}{\mu_{2}(b_{\pm},\delta_{\pm})}<-1, \quad t_2(b_{\pm},\delta_{\pm}) = -\frac{t_{\pm}}{\mu_{1}(b_{\pm},\delta_{\pm})}\in(0,1).\]
		Thus $ h_3>0,h_4>0 $. However, still by \eqref{eqn:ConnectingCondII}, one has 
		\begin{gather*}
			\left\{\begin{split}
				\frac{t_{-}c}{b_{-}}t_1(b_{+},\delta_{+})+\frac{2t_{-}b_{+}}{t_{+}c}t_1(b_{+},\delta_{+}) -\frac{t_{-}c}{b_{-}} = h_3t_{2}(b_{-},\delta_{-}),\\
				\frac{t_{+}c}{b_{+}}t_1(b_{-},\delta_{-})+\frac{2t_{+}b_{-}}{t_{-}c}t_1(b_{-},\delta_{-}) -\frac{t_{+}c}{b_{+}} = h_4t_{2}(b_{+},\delta_{+}).
			\end{split}\right.
		\end{gather*}
		The left-hand sides of the above equations are negative, while the right-hand sides are positive. This contradiction proves the desired result. 
		
	\end{proof}
	The above theorem shows that when $ \delta_{+}\delta_{-}<0 $, the eigenvectors of $ \hatHtwo(0) $ that correspond to the zero eigenvalue at $ k=0 $ can be properly chosen so that
	\begin{equation}\label{eqn:TwoFoldII}
		\mathbf{w}^{A} = \{ (0,0,0,\iu x_{_n} ,0,-\iu x_{_n})\}_{n\in \mathbb{Z}},\quad \mathbf{w}^{B} = \{ (y_{_n},z_{_n},y_{_n},0 ,0,0) \}_{n\in \mathbb{Z}},
	\end{equation} 
	where $ \{ x_{_n},y_{_n},z_{_n} \}_{n\in \mathbb{Z}} $ are real parameters. Unlike the type-I interface, the two modes $\mathbf{w}^{A}$, $\mathbf{w}^{B}$ are not related by the antiunitary reflection symmetry $\widehat{R}(0)$. Therefore, no Kramers-type protection follows from this symmetry.

	\subsection{Local Behavior of $ E_{\rmII}(k) $ near $ k=0 $}
	In this section, we study the local behavior of $ E_{\rmII}(k) $ given the existence of zero eigenvalues. As noted in \Cref{rmk:spectral_isolation}, the zero eigenvalue is two-fold and isolated from the remaining spectrum of $\widehat{H}_{\rmII}(0)$, so first-order degenerate perturbation theory applies. 
	Similarly, one has
	\begin{theorem}\label{thm:CrossingCurveII}
		When $ \delta_{+}\delta_{-}<0 $, there exist two distinct point spectra $ E_{\rmII,-}(k)$, $E_{\rmII,+}(k) $ for sufficiently small but nonzero $ |k| $ that admit the following asymptotic expansion.
		\begin{equation}
			E_{\rmII,\pm}(k) = \pm E^{(1)}_{\rmII} | k | + \mathcal{O}(k^2).	
		\end{equation}
		Here $ E^{(1)}_{\rmII} $ is a nonzero real constant determined from the parameters $ ( b_{\pm},\delta_{\pm},c) $.
	\end{theorem} 
	To prove this theorem, we will first characterize the eigenfunctions $ \mathbf{w}^{A},\mathbf{w}^{B} $.
	\begin{lemma}\label{lem:EigenvecFormII}
		When $ \delta_{+}\delta_{-}<0 $, the eigenvectors of the Hamiltonian $ \hatHtwo $, $ \mathbf{w}^{A},\mathbf{w}^{B} $, that correspond to the zero eigenvalue can be properly chosen such that
		\begin{equation}\label{eqn:EigenvecFormII}
			\mathbf{w}^{A} = \{ (0,0,0,\iu x_{_n},0,-\iu x_{_n})\}_{n\in\mathbb{Z}},\quad \mathbf{w}^{B} = \{ (y_{_n},z_{_n},y_{_n},0 ,0,0) \}_{n\in \mathbb{Z}},
		\end{equation}
		such that $ x_{_n}>0,y_{_n}<0 $ for all $ n\in  \mathbb{Z} $.
	\end{lemma}
	\begin{proof}
		After multiplying by suitable complex constants, one can check that the eigenvectors take the form \eqref{eqn:EigenvecFormII}. It remains to prove that the coordinates $ x_n,y_n $ satisfy the positivity constraints. For simplicity, we suppose $ \delta_{+}>0,\delta_{-}<0 $. 
		Therefore, from \eqref{eqn:DecayCondII} and \eqref{eqn:ConnectingCondII}, one has  
		\begin{gather*}
			\begin{split}
				&Q_{_{A,-1}}(0)v_{3}(b_{+},\delta_{+}) = \chi_{_{-1}},\; Q_{_{A,-2}}(0)\chi_{_{-1}} = h_3v_3(b_{-},\delta_{-}),\\
				&[Q_{_{B,0}}(0)]^{-1}v_2(b_{+},\delta_{+}) = \xi_{_{0}},\; [Q_{_{B,-1}}(0)]^{-1}\xi_{_{0}}=   h_5v_1(b_{-},\delta_{-})+h_6v_2(b_{-},\delta_{-}).
			\end{split}
		\end{gather*}
		For $ \{\chi_{_{n}}\}_{n\in\mathbb{Z}} $, one has by direct calculation
		\[ \chi_{_{-1}} = ( c/b_{-},-c/b_{-},0)^{T},\quad \chi_{_{-2}} = (ct_{-}/b_{-},- ct_{-}/b_{-},0)^{T}. \]
		By definition, one has 
		\begin{gather*}
			\chi_{_{n}} = [Q(b_{+},\delta_{+},0)]^{-n}\chi_{_{0}} = \mu_{3}^{-n}(b_{+},\delta_{+})\chi_{_{0}},\\
			\chi_{_{-n}} = [Q(b_{-},\delta_{-},0)]^{n-2}\chi_{_{-2}} = \mu_{3}^{n-2}(b_{-},\delta_{-})\chi_{_{-2}} .
		\end{gather*}
		Therefore, $ \{ x_{_n} \}_{n\in \mathbb{Z}} $ are positive for all $ n\in\mathbb{Z} $. \par 
		As for $ \{ \xi_{_{n}} \}_{n\in\mathbb{Z}} $, still by calculation, 
		\begin{gather*}
			\begin{split}
				&\xi_{_{0}} = (\frac{c^2}{2b_{+}b_{-}},\frac{c^2}{2b_{+}b_{-}},\frac{t_{2}(b_{+},\delta_{+})}{t_{+}}+\frac{c^2}{2b_{+}b_{-}})^{T},\\
				&\xi_{_{-1}} = (f_1,f_1,f_2)^{T} = h_5v_1(b_{-},\delta_{-})+h_6v_2(b_{-},\delta_{-}),
			\end{split}
		\end{gather*}
		where
		\[
		\begin{aligned}
			f_1
			&=
			\frac{ct_-}{2b_-}
			\left(
			\frac{t_2(b_+,\delta_+)}{t_+}
			+\frac{c^2}{2b_+b_-}
			\right),\\
			f_2
			&=
			f_1+\frac{c}{2b_-}.
		\end{aligned}
		\]
		In particular, $f_2>f_1>0$. Therefore, one has 
		\begin{gather*}
			\begin{split}
				&w_{_{1,0}} = w_{_{3,0}}= -\frac{t_2(b_{+},\delta_{+})}{t_{+}}<0,\quad w_{_{2,0}} = \frac{t_{2}(b_{+},\delta_{+})}{t_{+}}+\frac{c^2}{2b_{+}b_{-}},\\
				&w_{_{1,-1}} = w_{_{3,-1}}= -\frac{c}{2b_{-}}<0,\quad w_{_{2,-1}} = f_2.
			\end{split}
		\end{gather*}
		To investigate the quantity $ w_{_{1,n}} $ for $ n\le-2 $, one first observes that 
		\[
			\xi_{_n}
			=
			h_5\mu_1^{n+1}(b_-,\delta_-)
			v_1(b_-,\delta_-)
			+
			h_6\mu_2^{n+1}(b_-,\delta_-)
			v_2(b_-,\delta_-),
			\quad n\leq -2.
			\]
		Thus, one needs to give a precise description of the parameters $ h_5,h_6 $, which satisfy the following relations 
		\[ h_5+h_6 = f_2,\quad h_5t_1(b_{-},\delta_{-})+h_6t_2(b_{-},\delta_{-})= f_1. \]
		Directly, one has 
		\begin{gather*}
			h_5 =  \frac{f_1[t_2(b_{-},\delta_{-})-1]+\frac{c}{2b_{-}}t_2(b_{-},\delta_{-})}{t_2(b_{-},\delta_{-})-t_1(b_{-},\delta_{-})},\\
			h_6 = \frac{f_1[1-t_1(b_{-},\delta_{-})]-\frac{c}{2b_{-}}t_1(b_{-},\delta_{-})}{t_2(b_{-},\delta_{-})-t_1(b_{-},\delta_{-})}.
		\end{gather*}
			
		Therefore,
		\[
		w_{_{1,n}}=h_5\mu_1^{n+1}(t_1-1)
		+h_6\mu_2^{n+1}(t_2-1)
		=\frac{h_5\mu_1^{n+2}}{\mu_2}
		+\frac{h_6\mu_2^{n+2}}{\mu_1},
		\quad n\leq-2.
		\]
		Here we have used
		\[
		t_1-1=\frac{\mu_1}{\mu_2},
		\quad
		t_2-1=\frac{\mu_2}{\mu_1}.
		\]
		Combining the formulae of $h_5$ and $h_6$, one has 
		\begin{gather*}
			\begin{aligned}
				w_{_{1,n}}
				&=
				\frac{f_1}{t_2-t_1}
				\left[
				\frac{(t_2-1)\mu_1^{n+2}}{\mu_2}-	
				\frac{(t_1-1)\mu_2^{n+2}}{\mu_1}
				\right]\\
				&\quad +
				\frac{c}{2b_-(t_2-t_1)}
				\left[
				\frac{t_2\mu_1^{n+2}}{\mu_2}-
				\frac{t_1\mu_2^{n+2}}{\mu_1}
				\right].
			\end{aligned}
		\end{gather*}
		
		When $n$ is even, tedious calculations show that 
		\begin{gather*}
			\begin{aligned}
				w_{_{1,n}}
				&<
				\frac{c}
				{2b_{-}(t_2-t_1)}
				\left[
				\frac{t_2\mu_1^{n+2}}{\mu_2}
				-
				\frac{t_1\mu_2^{n+2}}{\mu_1}
				\right]
				\\
				&=
				\frac{c}
				{2b_{-}(t_2
					-t_1)}
				\left[
				-\frac{t_{-}\mu_1^{n+1}}{\mu_2}
				+
				\frac{t_{-}\mu_2^{n+1}}{\mu_1}
				\right]\leq 0.
			\end{aligned}
		\end{gather*}
		The equality comes from \Cref{lem:PropMatEigII}. \par 
		When $n$ is odd, one similarly has, for $n\leq -3$, 
		\begin{gather*}
			\begin{aligned}
				w_{_{1,n}}
				&<\frac{f_1}{t_2-t_1}
				\left[
				\frac{(t_2-1)\mu_1^{n+2}
				}{\mu_2}
				-
				\frac{(t_1-1)\mu_2^{n+2}}{\mu_1}
				\right]
				\\
				&=
				\frac{f_1}
				{t_2-t_1}
				\left[\mu_1^{n+1}-\mu_2^{n+1}\right]<0.
			\end{aligned}
		\end{gather*}
		Therefore, we have proved that $w_{_{1,n}}< 0$ for all
		$n\in\mathbb{Z}$, given that
		$\xi_{_{1}}=v_2(b_{+},\delta_{+})$.
	\end{proof}
	Now we can give a full proof of \Cref{thm:CrossingCurveII}. We suppose $\mathbf{w}^A$ and $\mathbf{w}^B$ are normalized. 
	\begin{proof}[Proof of \Cref{thm:CrossingCurveII}]
		
		A Taylor expansion gives
		\[ \hatHtwo(k) = \hatHtwo(0)+k\hatHtwo^{(1)} + \mathcal{O}(k^2). \]
		The first-order part $ \hatHtwo^{(1)}: l_{k}^2(\mathbb{Z};\mathbb{C}^6)\to l^{2}_{k}(\mathbb{Z};\mathbb{C}^6) $ is given by 
		\begin{gather}
			\left\{
			\begin{split}
				& (\widehat{H}^{(1)}_{\rmII}\mathbf{w})_{_{1,n}} = \iu c_{_{n,\rmII}} w_{_{6,n+1}},\\
				& (\widehat{H}^{(1)}_{\rmII}\mathbf{w})_{_{2,n}} = -\iu d_{_{n-2,\rmII}} w_{_{5,n-2}},\\
				& (\widehat{H}^{(1)}_{\rmII}\mathbf{w})_{_{3,n}} = 0,
			\end{split}\right.\quad 
			\left\{
			\begin{split}
				& (\widehat{H}^{(1)}_{\rmII}\mathbf{w})_{_{4,n}} = 0,\\
				& (\widehat{H}^{(1)}_{\rmII}\mathbf{w})_{_{5,n}} = \iu d_{_{n,\rmII}} w_{_{2,n+2}},\\
				& (\widehat{H}^{(1)}_{\rmII}\mathbf{w})_{_{6,n}} = -\iu c_{_{n-1,\rmII}} w_{_{1,n-1}}.
			\end{split}
			\right.
		\end{gather}
		Following the aforementioned procedure, we calculate the quantities $ (\mathbf{w}^{A},\hatHtwo^{(1)}\mathbf{w}^{B}) $.
		\begin{gather*}
			\begin{split}
				(\mathbf{w}^{A},\hatHtwo^{(1)}\mathbf{w}^{A}) &= (\mathbf{w}^{B},\hatHtwo^{(1)}\mathbf{w}^{B}) = 0,\\
				(\mathbf{w}^{A},\hatHtwo^{(1)}\mathbf{w}^{B}) &= \sum_{n\in\mathbb{Z}}\iu c_{_{n-1,\rmII}}y_{_{n-1}} x_{_n}, \\
				(\mathbf{w}^{B},\hatHtwo^{(1)}\mathbf{w}^{A}) &= -\sum_{n\in\mathbb{Z}}\iu c_{_{n,\rmII}}y_{_{n}} x_{_{n+1}} =- (\mathbf{w}^{A},\hatHtwo^{(1)}\mathbf{w}^{B}).
			\end{split}
		\end{gather*}
		From \Cref{lem:EigenvecFormII}, the eigenvectors can be properly chosen such that $ x_{_n}>0,y_{_n}<0 $. Therefore, it is obvious that 
		\[ (\mathbf{w}^{A},\hatHtwo^{(1)}\mathbf{w}^{B})\neq 0. \]
		And the matrix 
		\begin{equation}
			\mathbf{M}_{0,\rmII} \triangleq 
			\begin{pmatrix}
				( \mathbf{w}^{A},\hatHtwo^{(1)}\mathbf{w}^{A})  &(\mathbf{w}^{A},\hatHtwo^{(1)}\mathbf{w}^{B}) \\ 
				(\mathbf{w}^{B},\hatHtwo^{(1)}\mathbf{w}^{A})& (\mathbf{w}^{B},\hatHtwo^{(1)}\mathbf{w}^{B} )
			\end{pmatrix},
		\end{equation}
		is nontrivial, with two distinct real eigenvalues 
		$ \pm\iu(\mathbf{w}^{A},\hatHtwo^{(1)}\mathbf{w}^{B}) $, which means that the crossing at $ k=0 $ is also non-tangential. And we also have 
			\[E_{\rmII,\pm} (k)= \pm |( \mathbf{w}^{A},\hatHtwo^{(1)}\mathbf{w}^{B})|  \,|k|+\mathcal{O}(k^2) \triangleq \pm E_{\rmII}^{(1)}|k|+\mathcal{O}(k^2).\]
	\end{proof}
	
	\begin{remark}
		The above results only concern the local behavior near $ k=0 $. In the next section, we numerically show that the energy curves $ E_{\rmI}(k),E_{\rmII}(k) $ do not vanish for $ k\neq 0 $. Therefore, there exist only two counter-propagating zero-energy edge states concentrated near $E_{\rmI}(0)=E_{\rmII}(0)=0 $.
	\end{remark}
	\section{Numerical Illustration}\label{sec:Numerics}
	In this section, we present a series of numerical results to verify the above analysis. First, to numerically investigate the existence of zero eigenvalues at $ k=0 $, we calculate the spectra of the Hamiltonian operators $\hatHone(k),\hatHtwo(k) $ for various choices of hopping coefficients and $ k\in [-\pi,\pi) $, with special emphasis on the point spectrum near $ 0 $. 
	Given the existence of zero eigenvalues at $ k=0 $, we also numerically simulate the time dynamics of wave packets concentrated near $ E_{\rmI}(0)=E_{\rmII}(0)=0 $. More specifically, we numerically solve the following systems:
	\begin{equation}\label{eqn:TimeSchrodinger}
		\iu \partial_{t}\varPhi = \hatHone\varPhi,\; \iu \partial_{t}\varPsi = \hatHtwo\varPsi, \quad t\in [0,T],
	\end{equation}
	with appropriately chosen initial values. Previous analysis shows that the localized wave packets will propagate along the interface with nonzero group velocity. Since the zero eigenvalues of $ \hatHone(0) $ and $ \hatHtwo(0) $ are two-fold, the propagation direction depends on the choice of initial value. 
	
	\subsection{The Spectra of $ \hatHone(k) $ and $ \hatHtwo(k) $}
	In this part, we investigate the spectrum of the Hamiltonians $ \hatHone(k),\hatHtwo(k) $. We adopt the supercell method to calculate the point spectrum. If eigenvectors localize on the artificial boundary, we remove the corresponding eigenvalues. \par
\begin{figure}[htbp] 
		\centering
		\begin{subfigure}{0.35\textwidth}
			\includegraphics[width=\textwidth]{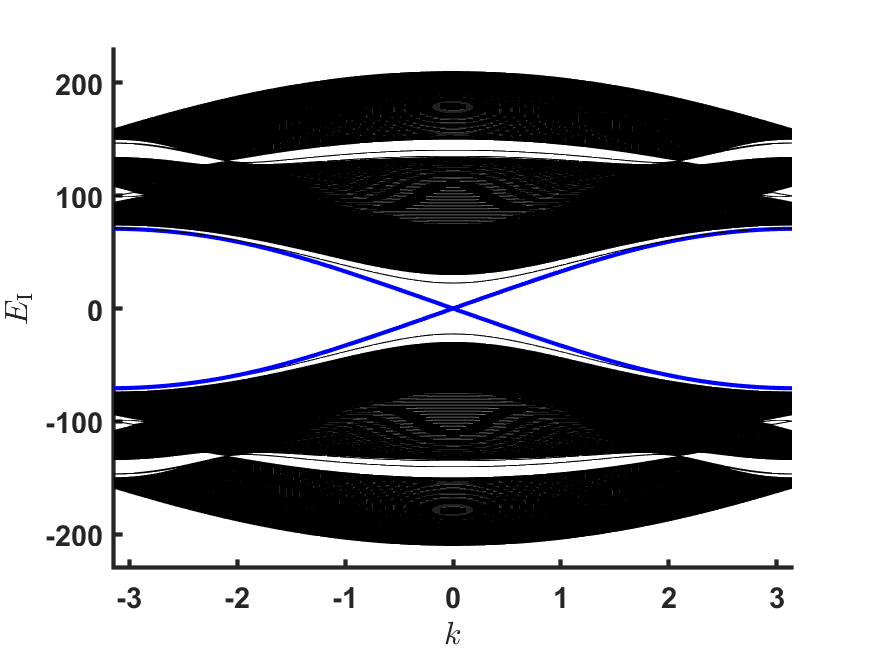}
			\caption{$ b_{\pm} = 60,\delta_{\pm}=30 $.}
			\label{subfig:ppzero}
		\end{subfigure}
		\hspace{-0.8cm}
		\begin{subfigure}{0.35\textwidth}
			\includegraphics[width=\textwidth]{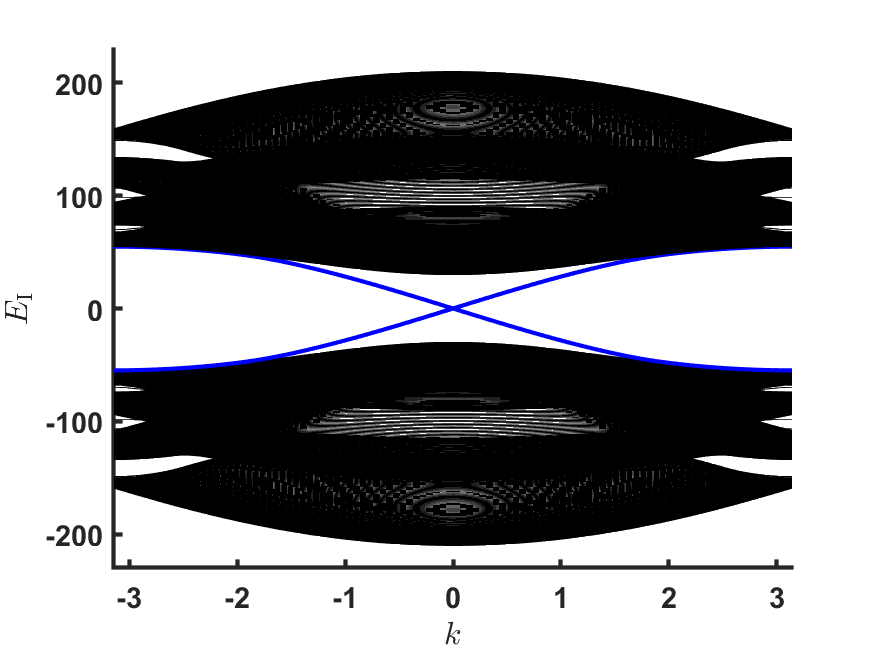}
			\caption{$ b_{\pm} = 60,\delta_{\pm}=\pm 30 $.}
			\label{subfig:pmzero}
		\end{subfigure}
		\hspace{-0.8cm}
		\begin{subfigure}{0.35\textwidth}
			\includegraphics[width=\textwidth]{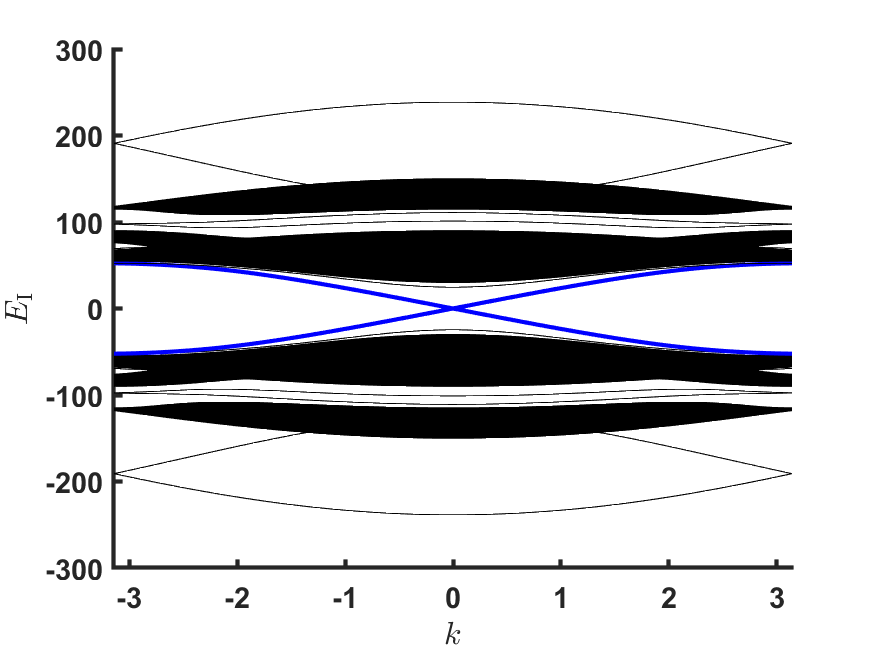}
			\caption{$ b_{\pm} = 60,\delta_{\pm}=- 30 $.}
			\label{subfig:mmzero}
		\end{subfigure}
		\caption{The spectra of the Hamiltonian operator $ \hatHone(k) $ for $ k\in [-\pi,\pi] $ when the hopping coefficient $ c $ across the boundary is chosen properly. The point spectrum near zero is shown in thick blue lines. There is a crossing at $ k=0 $ in each figure. Near $ k=0 $, the slope of $ E_{\rmI}(k) $ is nonzero. }
		\label{fig:tbedge_zero}
	\end{figure}

 
	First, for $ \hatHone(k) $, \Cref{thm:ZeroEigenvalueI} shows that given $ b_{\pm}\neq 0,\delta_{\pm}\neq 0 $, there exists a unique $ c $ such that there exists a two-fold zero eigenvalue at $ k=0 $. It can be easily seen from \Cref{fig:tbedge_zero} that the zero eigenvalue exists for various parameters $ b_{\pm},\delta_{\pm} $. It is also obvious that the two eigencurves cross at $ k=0 $. \par 
\begin{figure}[htbp] 
		\centering
		\begin{subfigure}{0.35\textwidth}
			\includegraphics[width=\textwidth]{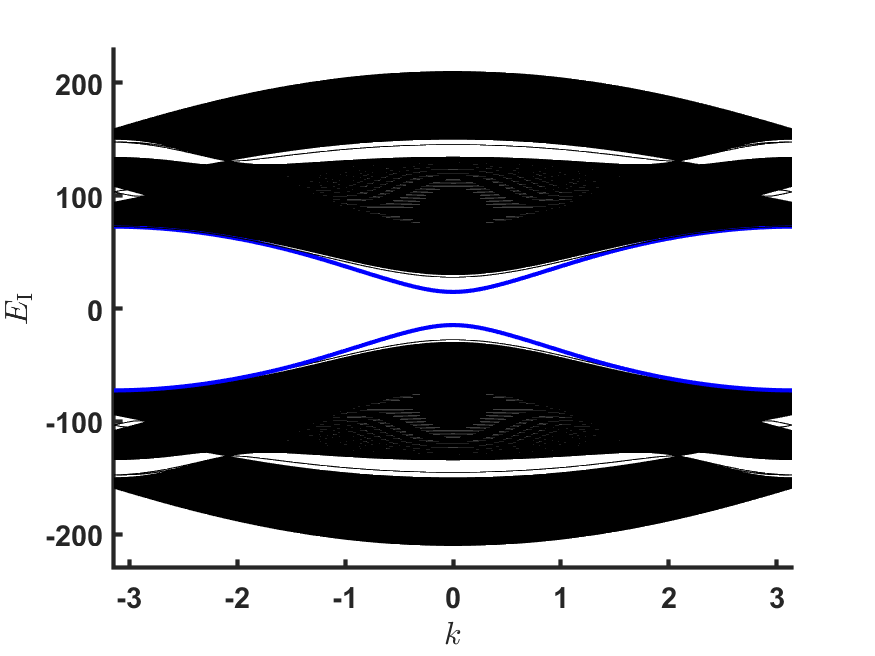}
			\caption{$ b_{\pm} = 60,\delta_{\pm}=30 $.}
			\label{subfig:ppnonzero}
		\end{subfigure}
		\hspace{-0.8cm}
		\begin{subfigure}{0.35\textwidth}
			\includegraphics[width=\textwidth]{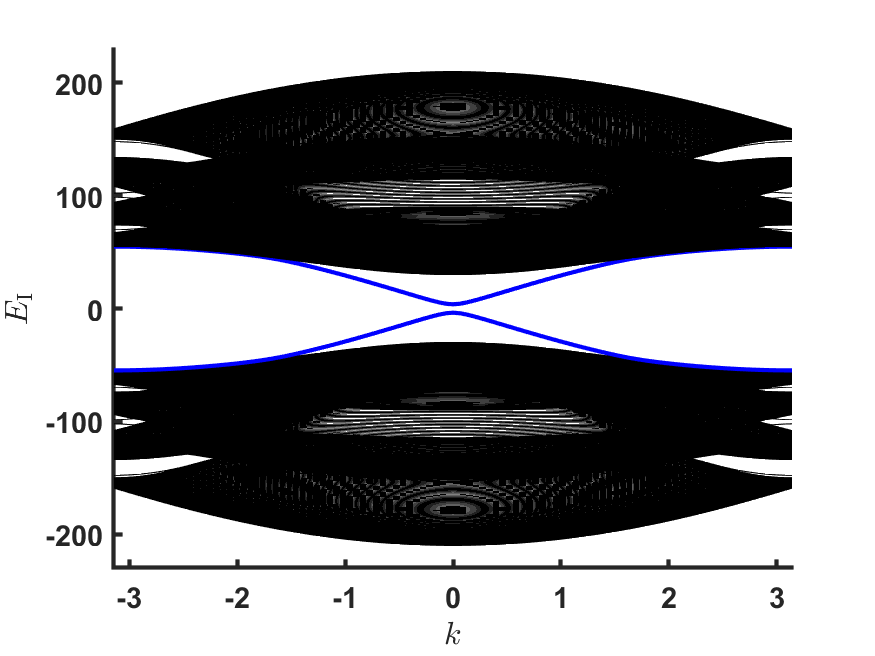}
			\caption{$ b_{\pm} = 60,\delta_{\pm}=\pm 30 $.}
			\label{subfig:pmnonzero}
		\end{subfigure}
		\hspace{-0.8cm}
		\begin{subfigure}{0.35\textwidth}
			\includegraphics[width=\textwidth]{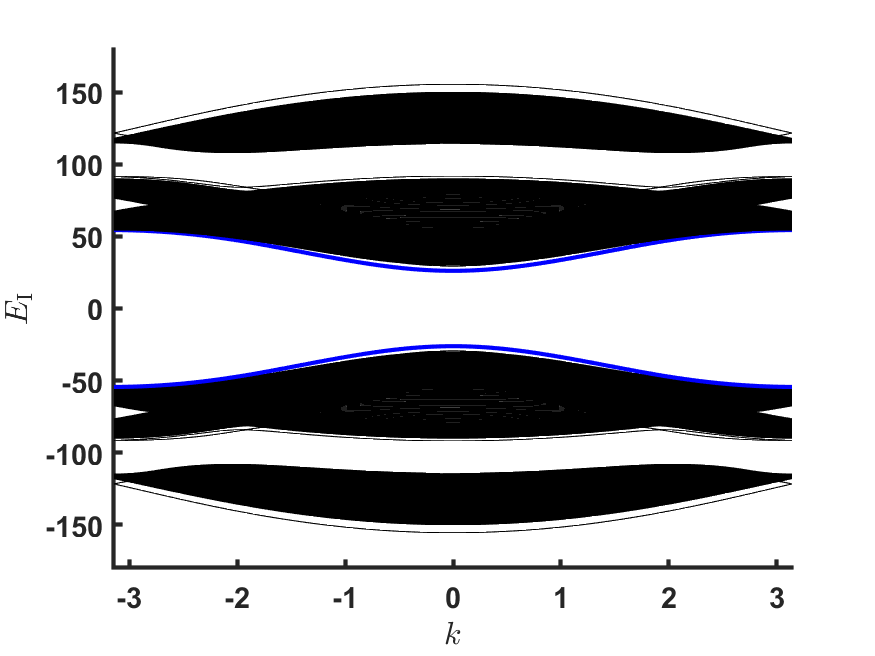}
			\caption{$ b_{\pm} = 60,\delta_{\pm}=- 30 $.}
			\label{subfig:mmnonzero}
		\end{subfigure}
		\caption{The spectra of the Hamiltonian operator $ \hatHone(k) $ for $ k\in [-\pi,\pi] $ when the hopping coefficient $ c=50 $. The point-spectrum branches near zero are shown by bold blue lines. There is no crossing at $ k=0 $ in any panel. }
		\label{fig:tbedge_nonzero}
	\end{figure}

	However, when the hopping coefficient $ c $ across the interface is chosen arbitrarily, zero eigenvalues do not exist in general, as shown in \Cref{fig:tbedge_nonzero}. It is evident that for the type-I interface, zero eigenvalues do not exist in general. \par 
	 \par
\begin{figure}[htbp] 
		\centering
		\begin{subfigure}{0.48\textwidth}
			\includegraphics[width=\textwidth]{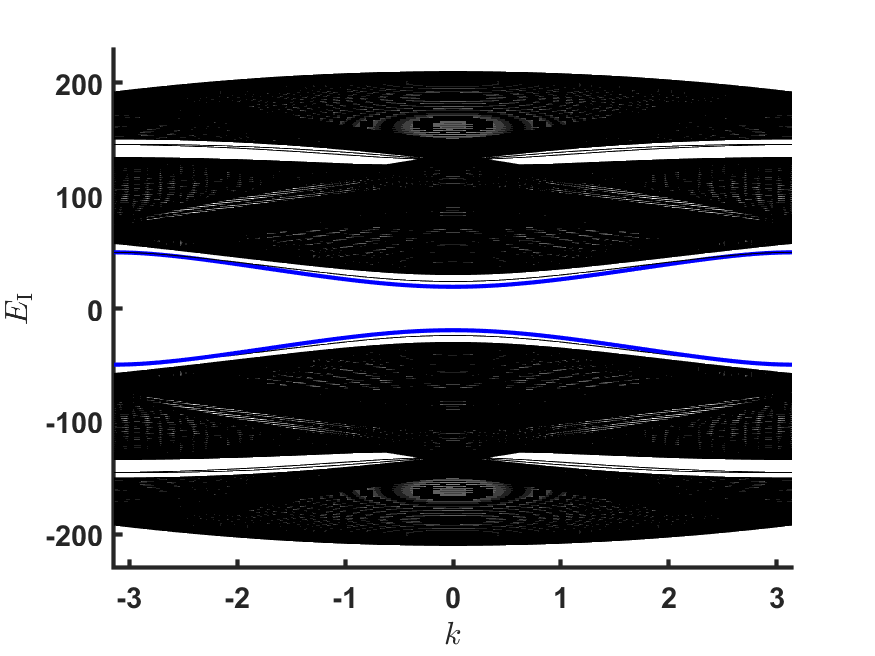}
			\caption{$ b_{\pm} = 60,\delta_{\pm}=30 $.}
			\label{subfig:ppnonzero2}
		\end{subfigure}
		\hspace{-0.8cm}
		\begin{subfigure}{0.48\textwidth}
			\includegraphics[width=\textwidth]{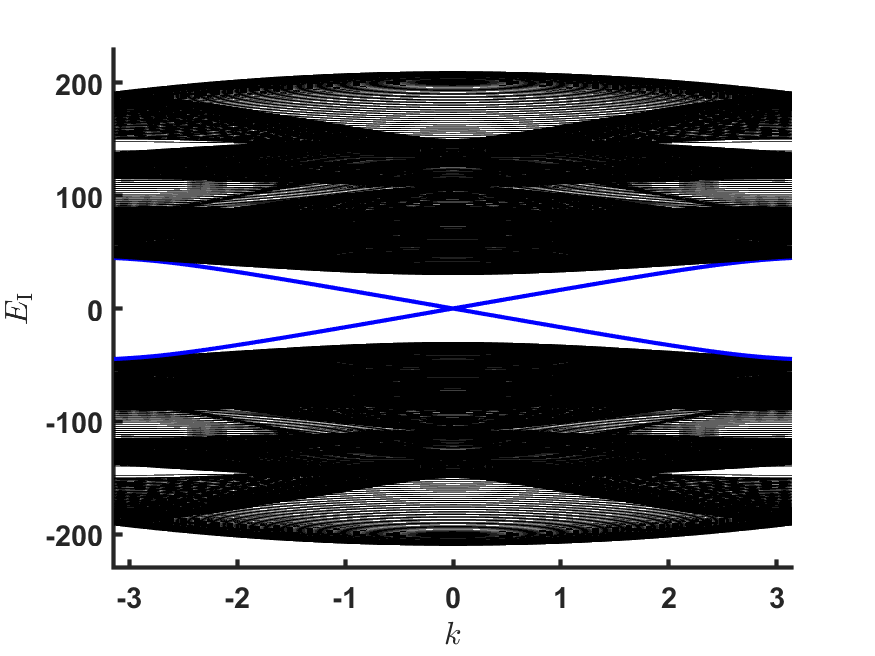}
			\caption{$ b_{\pm} = 60,\delta_{\pm}=\pm 30 $.}
			\label{subfig:pmzero2}
		\end{subfigure}
		\caption{The spectra of the Hamiltonian operator $ \hatHtwo(k) $ for $ k\in [-\pi,\pi] $ when the hopping coefficient $ c=50 $. The point spectrum near zero is shown in bold blue lines. When $ \delta_{+}\delta_{-}<0 $, two discrete eigenvalue branches cross at $ k=0 $. The zero-energy edge states do not exist when $ \delta_{+}\delta_{-}>0 $. }
		\label{fig:tbedge_2}
	\end{figure}

 
	Then, we calculate the spectrum of $ \hatHtwo(k) $. As shown in \Cref{lem:EigenvecFormII}, zero eigenvalues exist when $ \delta_{+}\delta_{-}<0 $. However, they never exist when $ \delta_{+}\delta_{-}>0 $. The corresponding numerical results are shown in \Cref{fig:tbedge_2}. The two eigenvalue curves also cross at $ k=0 $. \par 
	\clearpage
	\subsection{Dynamics of the Edge States}\label{subsec:dynamics}
	In this section, we numerically simulate the dynamics of wave packets for the Schr\"odinger-type equation \eqref{eqn:TimeSchrodinger}. Assuming the Hamiltonian operators $ \hatHone,\hatHtwo $ have zero eigenvalues at $ k=0 $, the wave packets localized nearby will propagate at nonzero group velocity by \Cref{thm:CrossingCurveI} and \Cref{thm:CrossingCurveII}. Moreover, the matrices $ \mathbf{M}_{0,\rmI},\mathbf{M}_{0,\rmII} $ have two distinct eigenvalues with opposite signs; therefore, wave packets can be chosen such that they propagate in either direction. For simplicity, we only choose the right-propagating state. 
	
	We adopt the supercell method to simulate the spatial operator. We use the fourth-order Runge-Kutta method to simulate the time dynamics. We plot the absolute value of the wave function $ \varPhi $ at each time. Later in this section, we label the interface configurations with red lines.\par 
		
\begin{figure}[htbp] 
			\centering
			\begin{subfigure}{0.48\textwidth}
				\includegraphics[width=\textwidth]{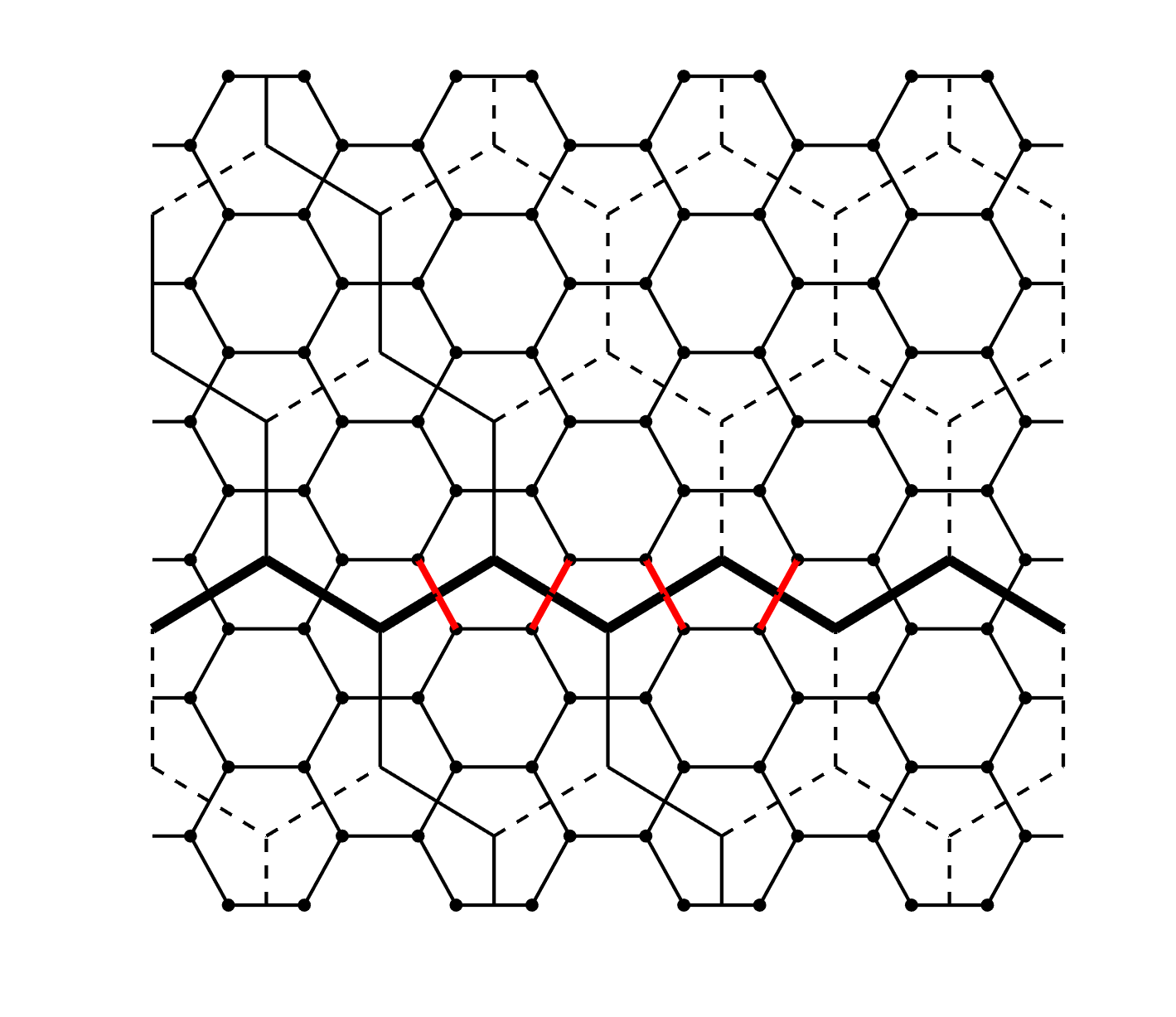}
				\caption{Perturbation at type-I interfaces.}
			\end{subfigure}
			\hspace{-0.8cm}
			\begin{subfigure}{0.48\textwidth}
				\includegraphics[width=\textwidth]{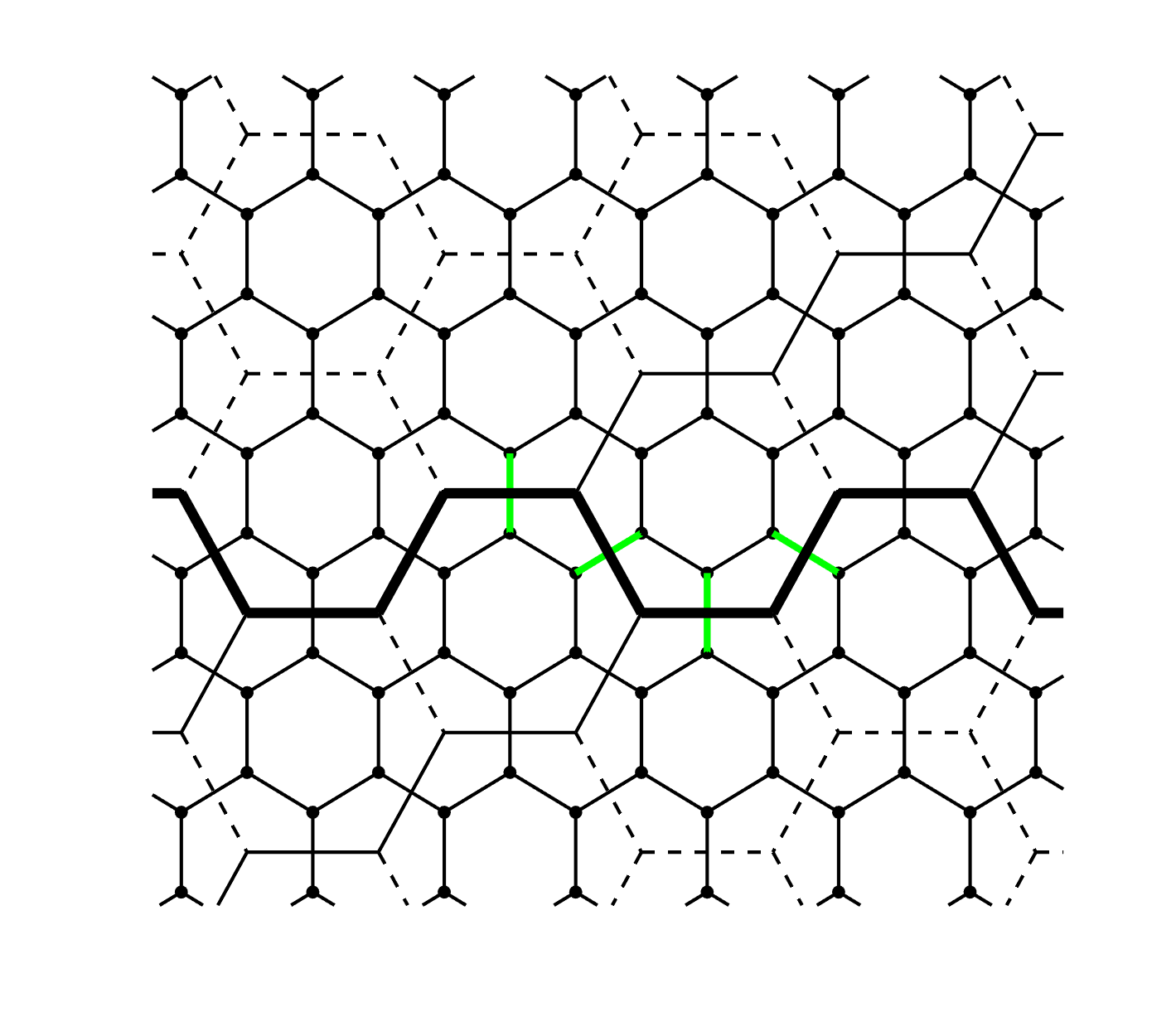}
				\caption{Perturbation at type-II interfaces.}
			\end{subfigure}
			\caption{Illustrations of local defects. We multiply the hopping coefficient $c$ across the interface by independent $\zeta \sim\mathrm{Uniform}[0,2] $. The perturbed coefficients are plotted with red and green lines for type-I and type-II interfaces, respectively.} %
			\label{fig:tbdefect}
		\end{figure}

	Since we do not define a topological invariant or prove topological protection against arbitrary perturbations, we simulate the behavior of wave packets passing through a local defect to illustrate the results and provide some intuitions. As shown in \Cref{fig:tbdefect}, we multiply each hopping coefficient $c$ by an independent random variable $\zeta\sim\mathrm{Uniform}[0,2]$. The perturbed coefficients are plotted with red and green lines for type-I and type-II interfaces, respectively. In the following, the interface and local defects are highlighted with red lines and dots. 
		
\begin{figure}[htbp] 
			\centering
			\begin{subfigure}{0.26\textwidth}
				\includegraphics[width=\textwidth]{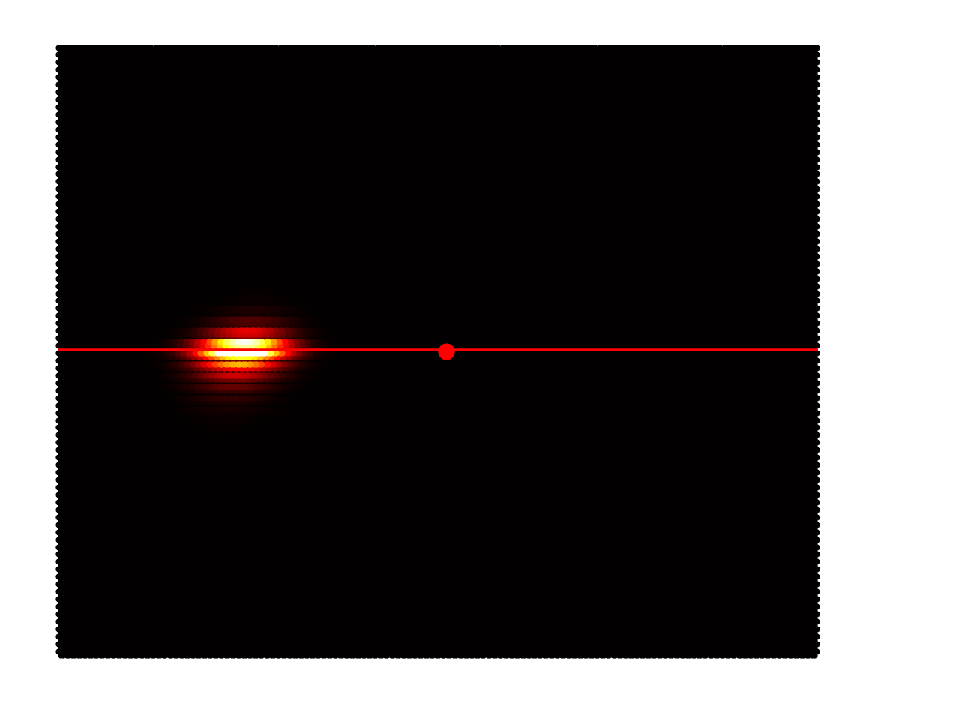}
				\caption{$ t=0 $.}
			\end{subfigure}
			\hspace{-0.6cm}
			\begin{subfigure}{0.26\textwidth}
				\includegraphics[width=\textwidth]{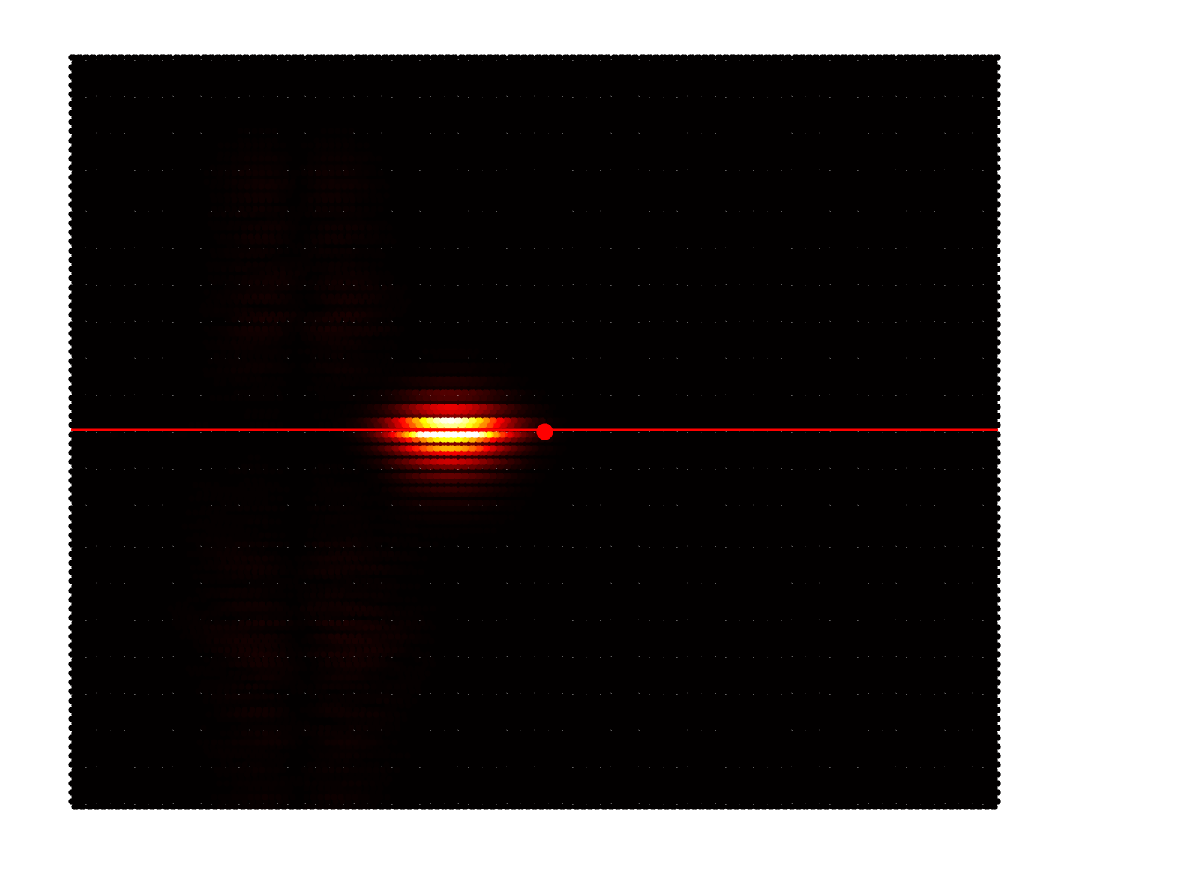}
				\caption{$ t=0.75 $.}
			\end{subfigure}
			\hspace{-0.6cm}
			\begin{subfigure}{0.26\textwidth}
				\includegraphics[width=\textwidth]{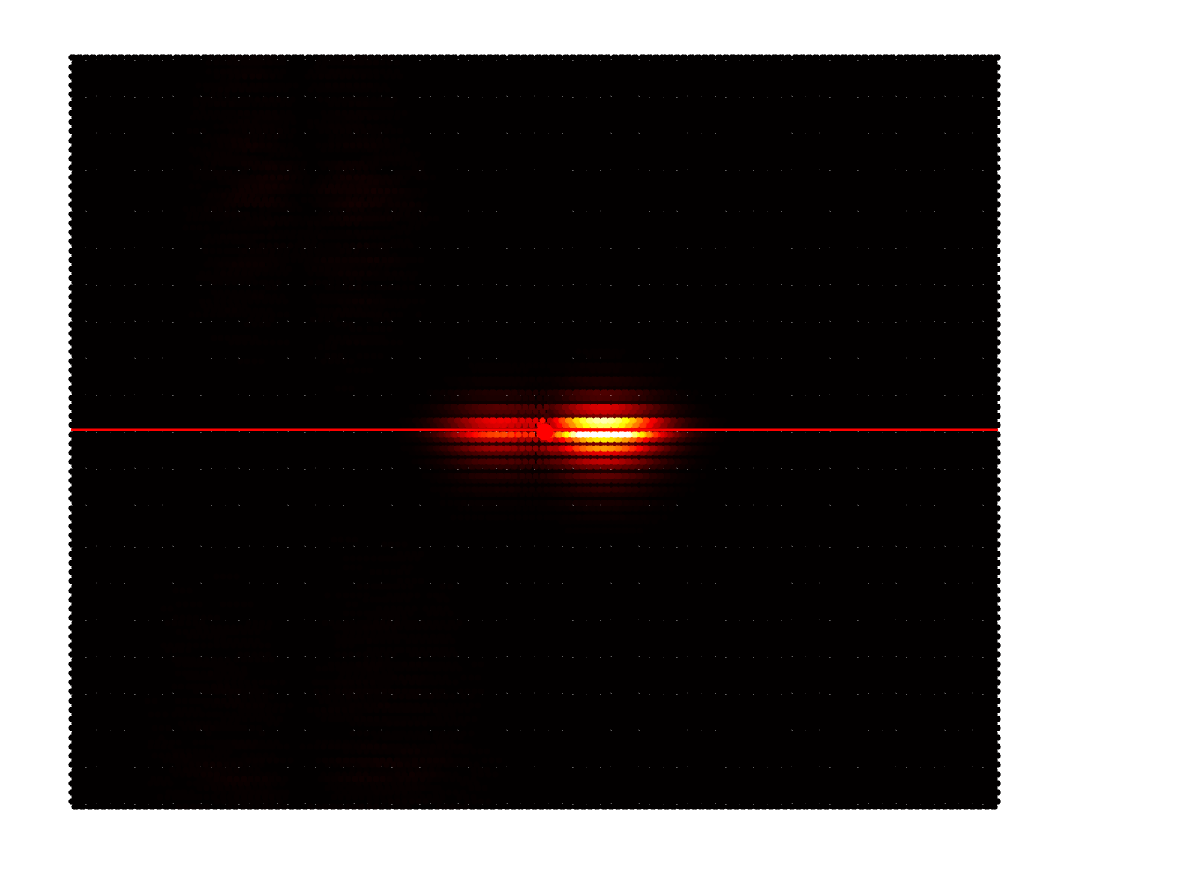}
				\caption{$ t=1.5 $.}
			\end{subfigure}
			\hspace{-0.6cm}
			\begin{subfigure}{0.26\textwidth}
				\includegraphics[width=\textwidth]{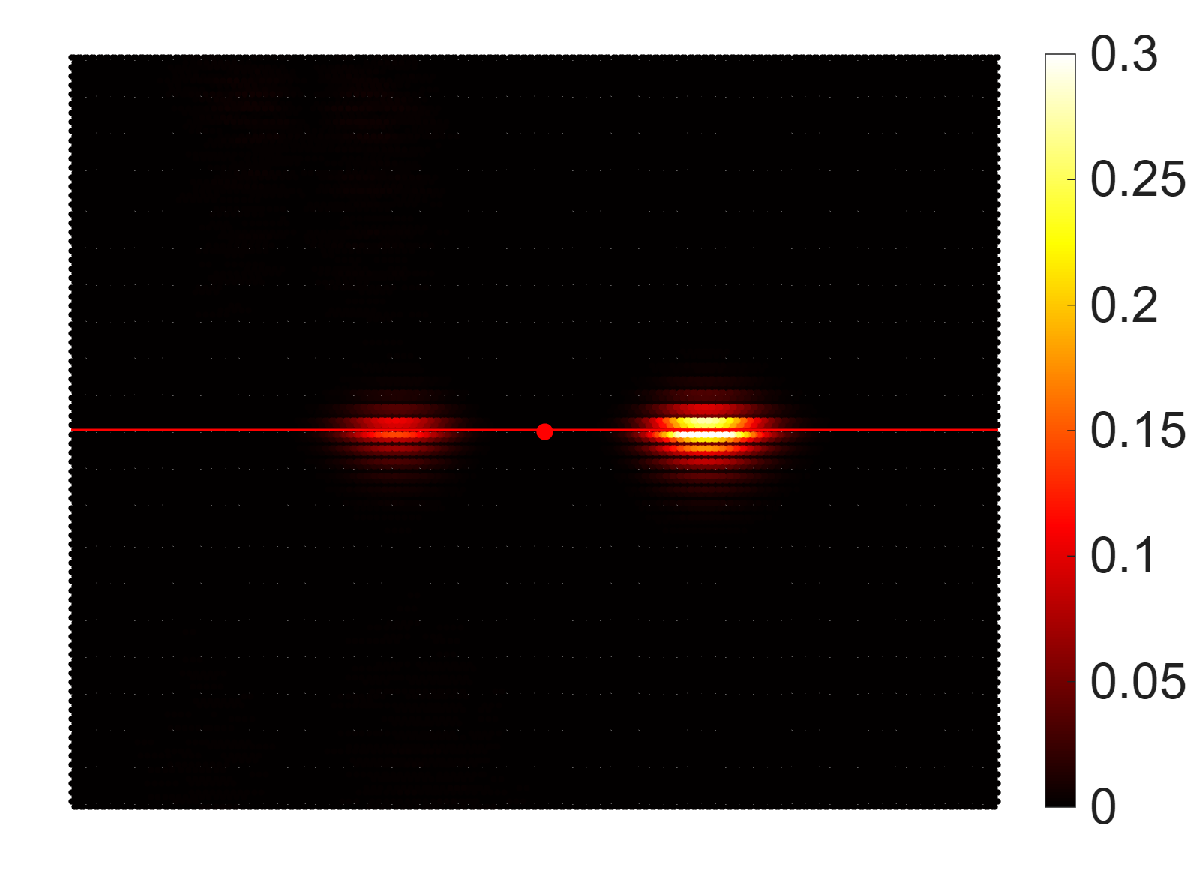}
				\caption{$ t=2 $.}
			\end{subfigure}
			\caption{ Dynamics of wave packets passing through a local defect at a type-I interface when $ \delta_{\pm} = \pm 30 $. In this case, there is a topological transition across the interface.}
			\label{fig:edgeIdef_diff}
		\end{figure}

\begin{figure}[htbp] 
			\centering
			\begin{subfigure}{0.26\textwidth}
				\includegraphics[width=\textwidth]{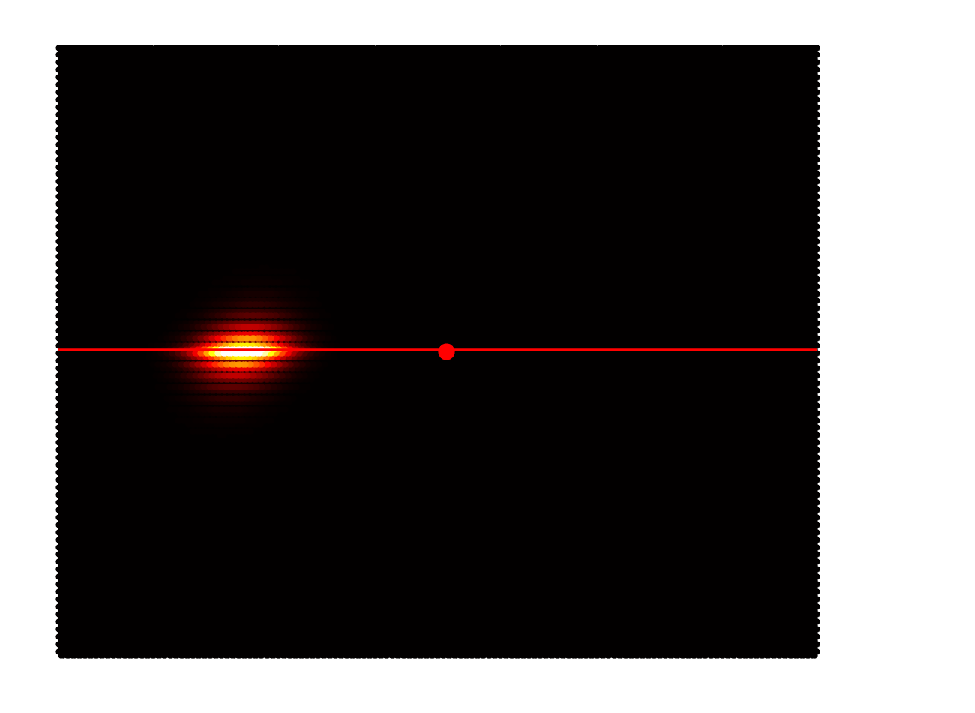}
				\caption{$ t=0 $.}
			\end{subfigure}
			\hspace{-0.6cm}
			\begin{subfigure}{0.26\textwidth}
				\includegraphics[width=\textwidth]{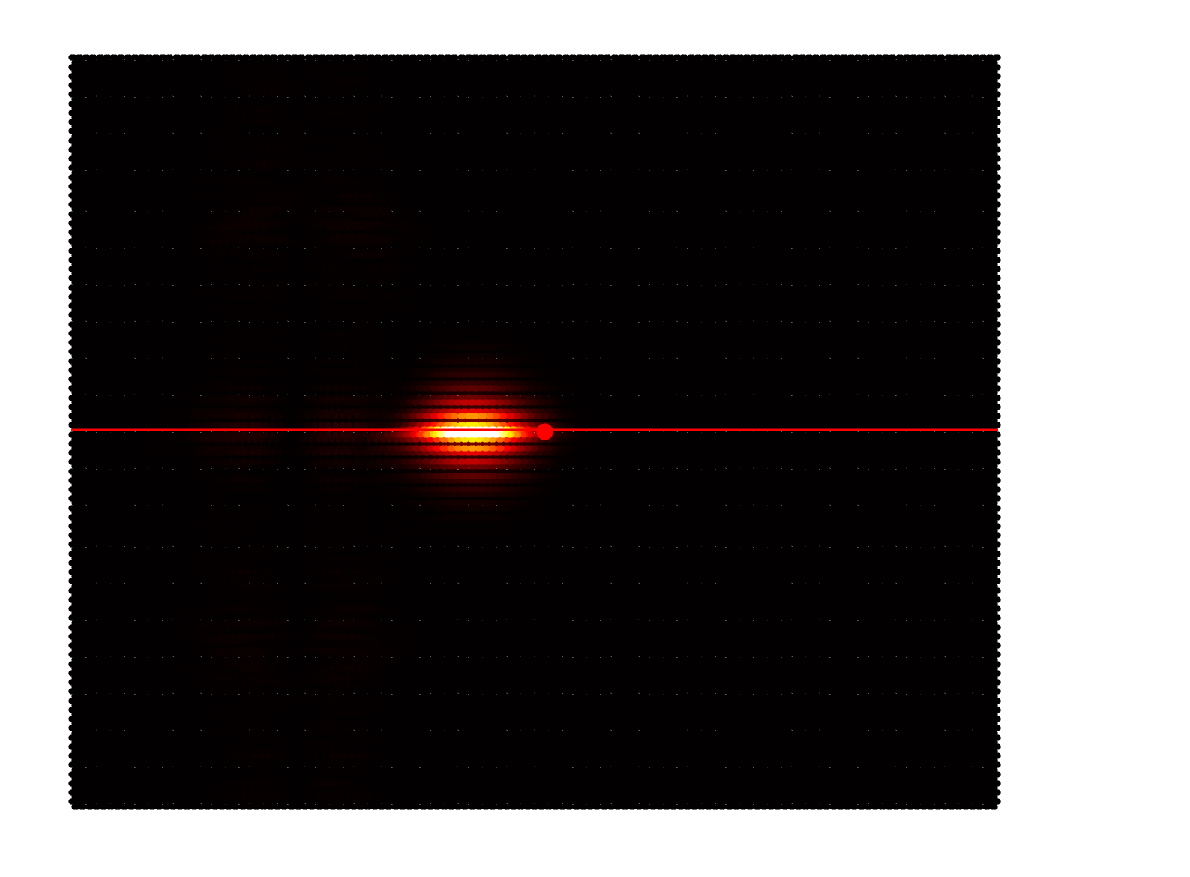}
				\caption{$ t=0.75 $.}
			\end{subfigure}
			\hspace{-0.6cm}
			\begin{subfigure}{0.26\textwidth}
				\includegraphics[width=\textwidth]{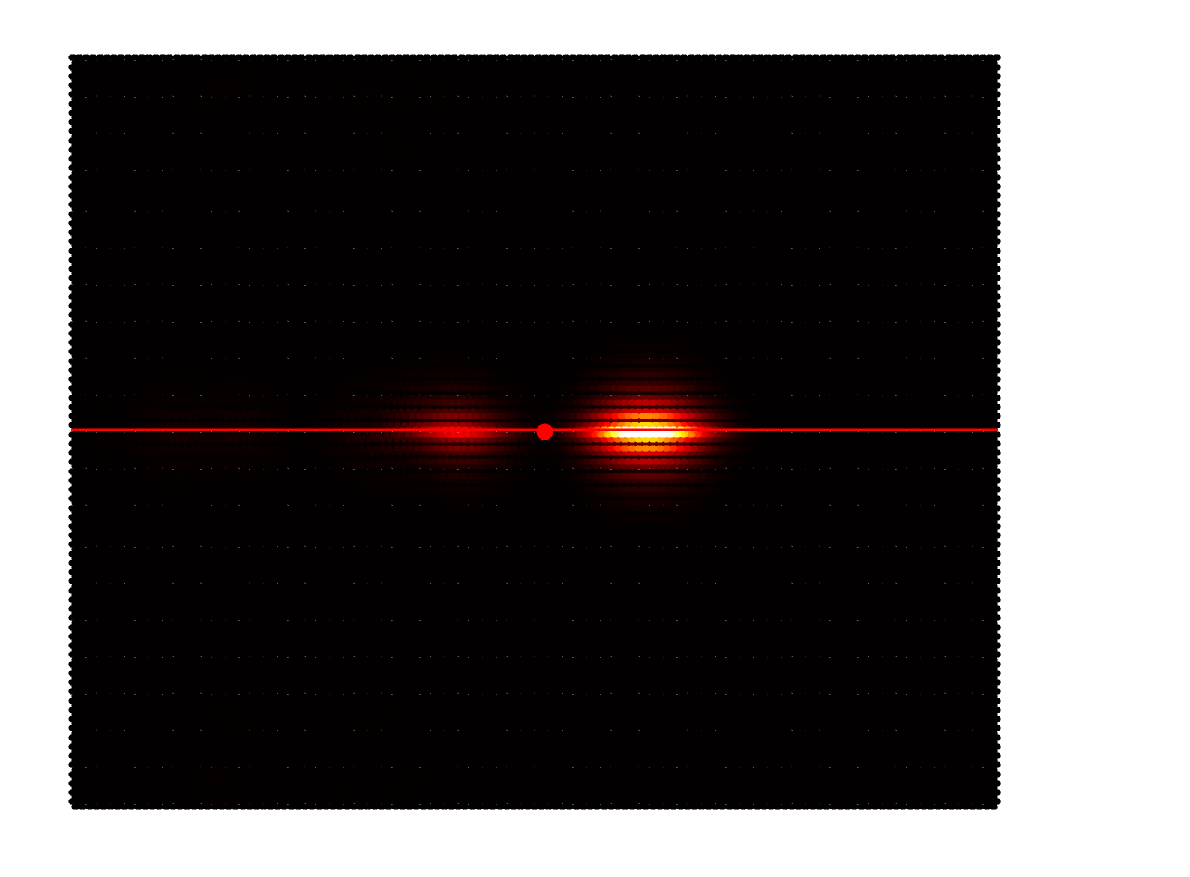}
				\caption{$ t=1.5 $.}
			\end{subfigure}
			\hspace{-0.6cm}
			\begin{subfigure}{0.26\textwidth}
				\includegraphics[width=\textwidth]{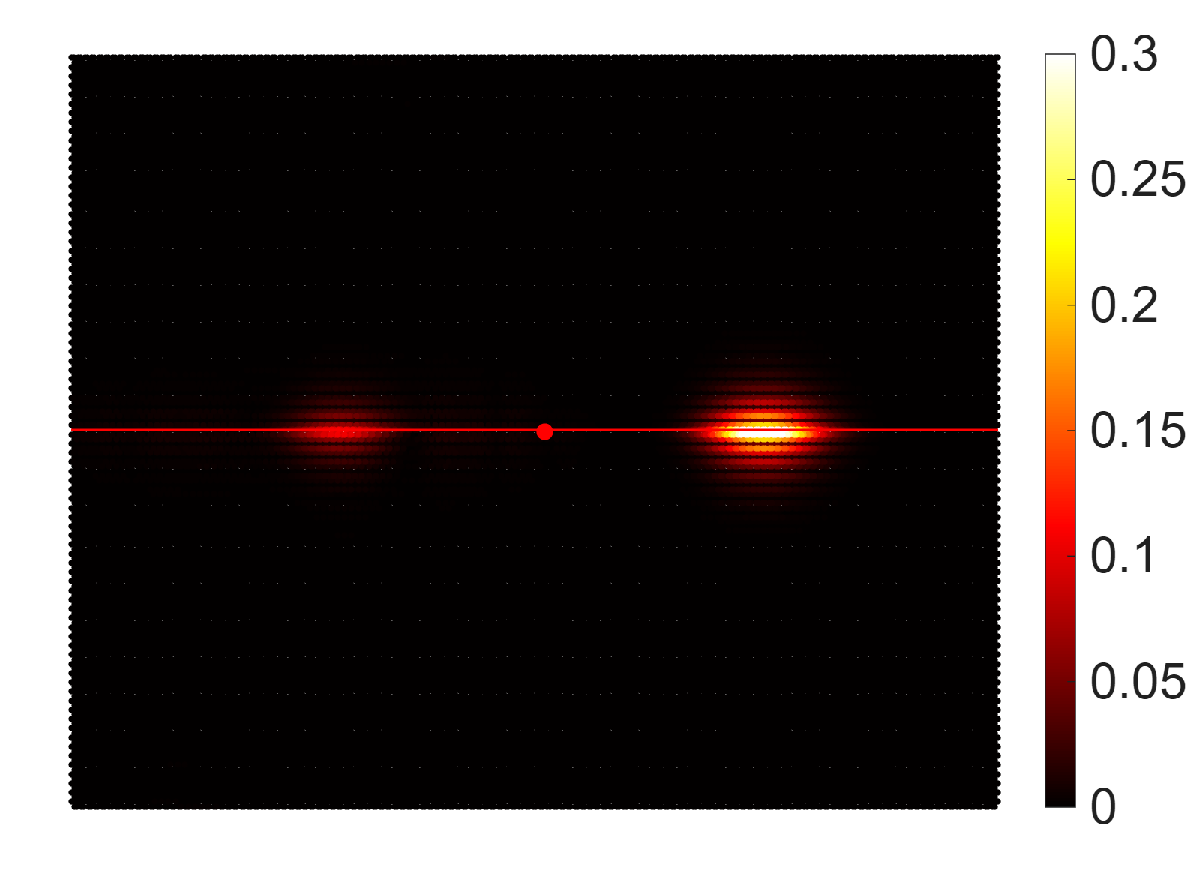}
				\caption{$ t=2 $.}
			\end{subfigure}\\
			\begin{subfigure}{0.26\textwidth}
				\includegraphics[width=\textwidth]{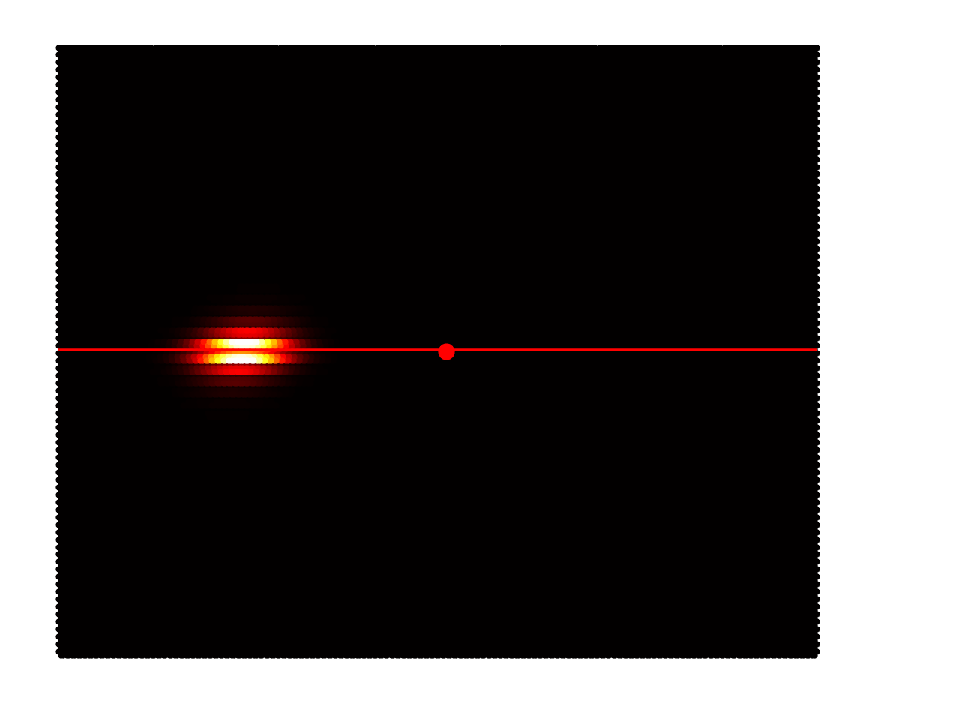}
				\caption{$ t=0 $.}
			\end{subfigure}
			\hspace{-0.6cm}
			\begin{subfigure}{0.26\textwidth}
				\includegraphics[width=\textwidth]{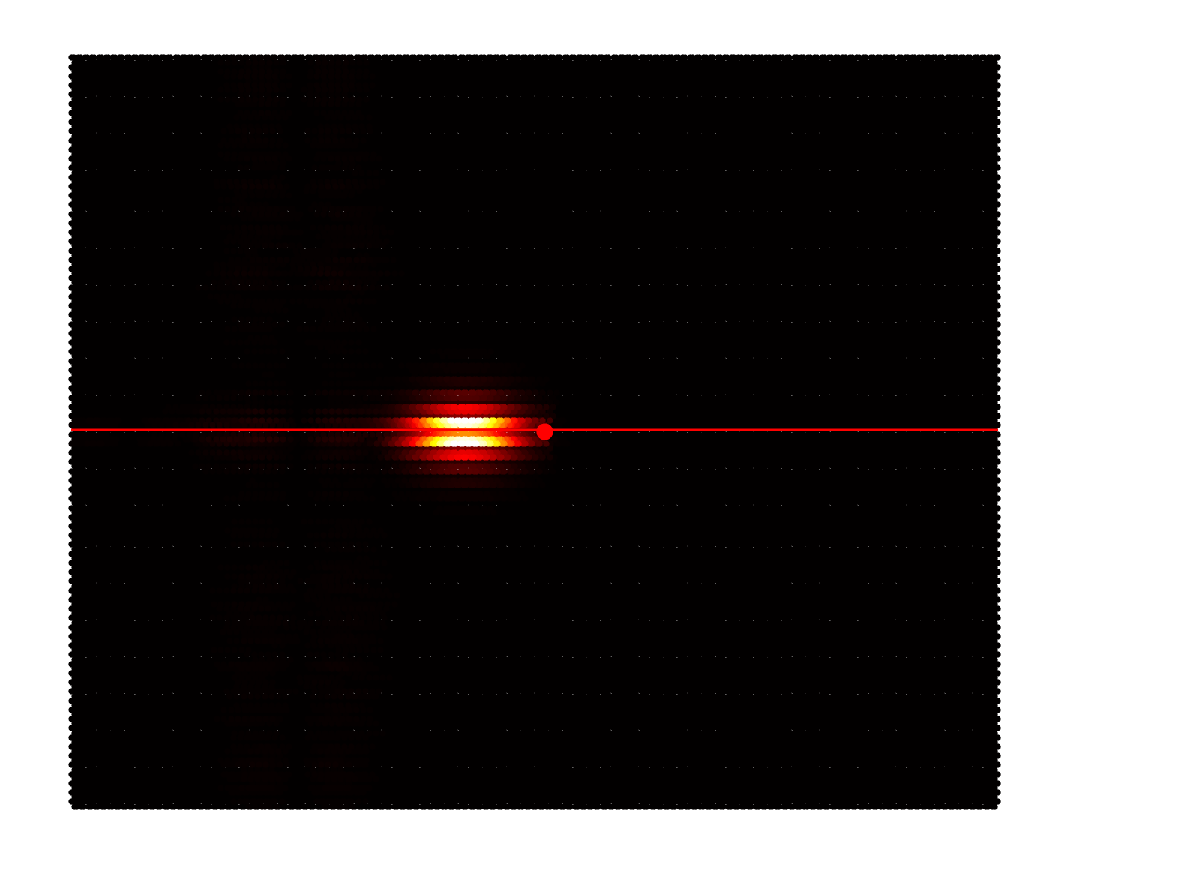}
				\caption{$ t=1 $.}
			\end{subfigure}
			\hspace{-0.6cm}
			\begin{subfigure}{0.26\textwidth}
				\includegraphics[width=\textwidth]{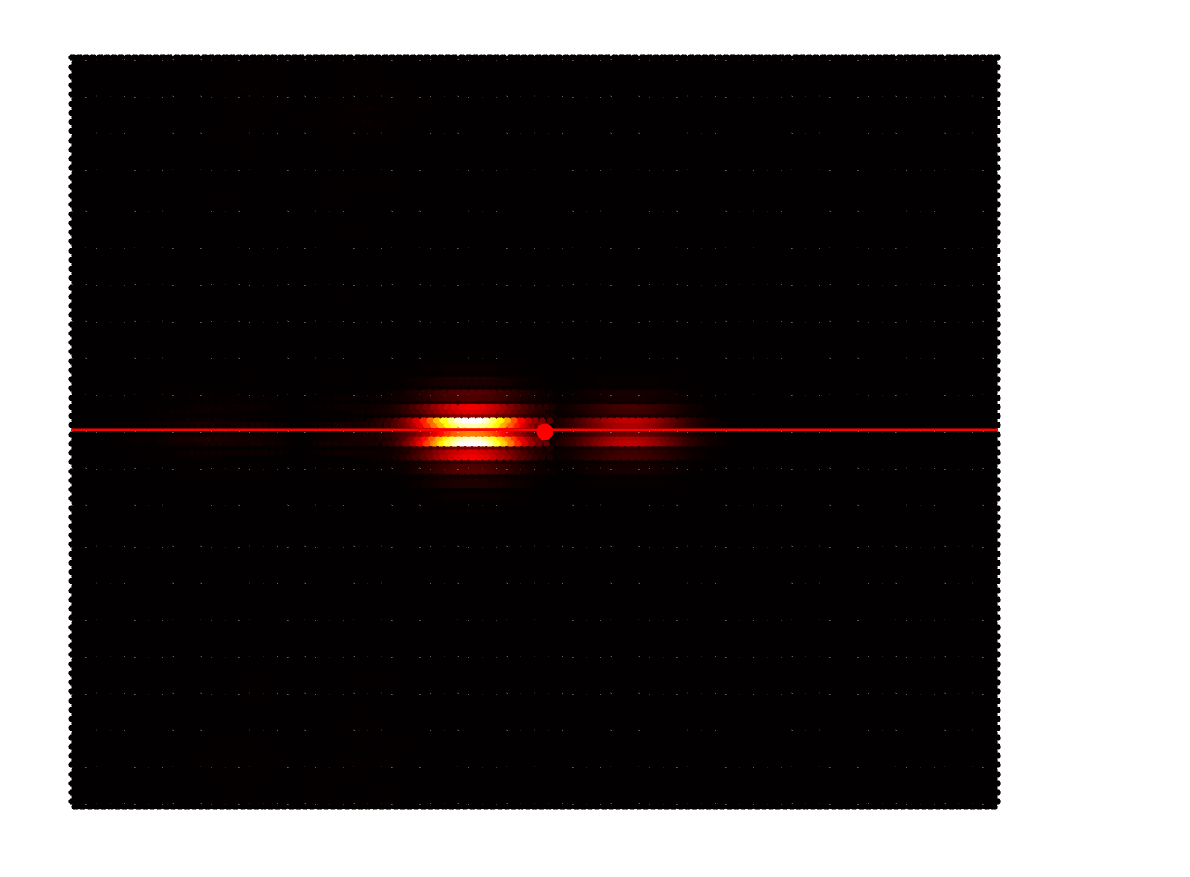}
				\caption{$ t=2 $.}
			\end{subfigure}
			\hspace{-0.6cm}
			\begin{subfigure}{0.26\textwidth}
				\includegraphics[width=\textwidth]{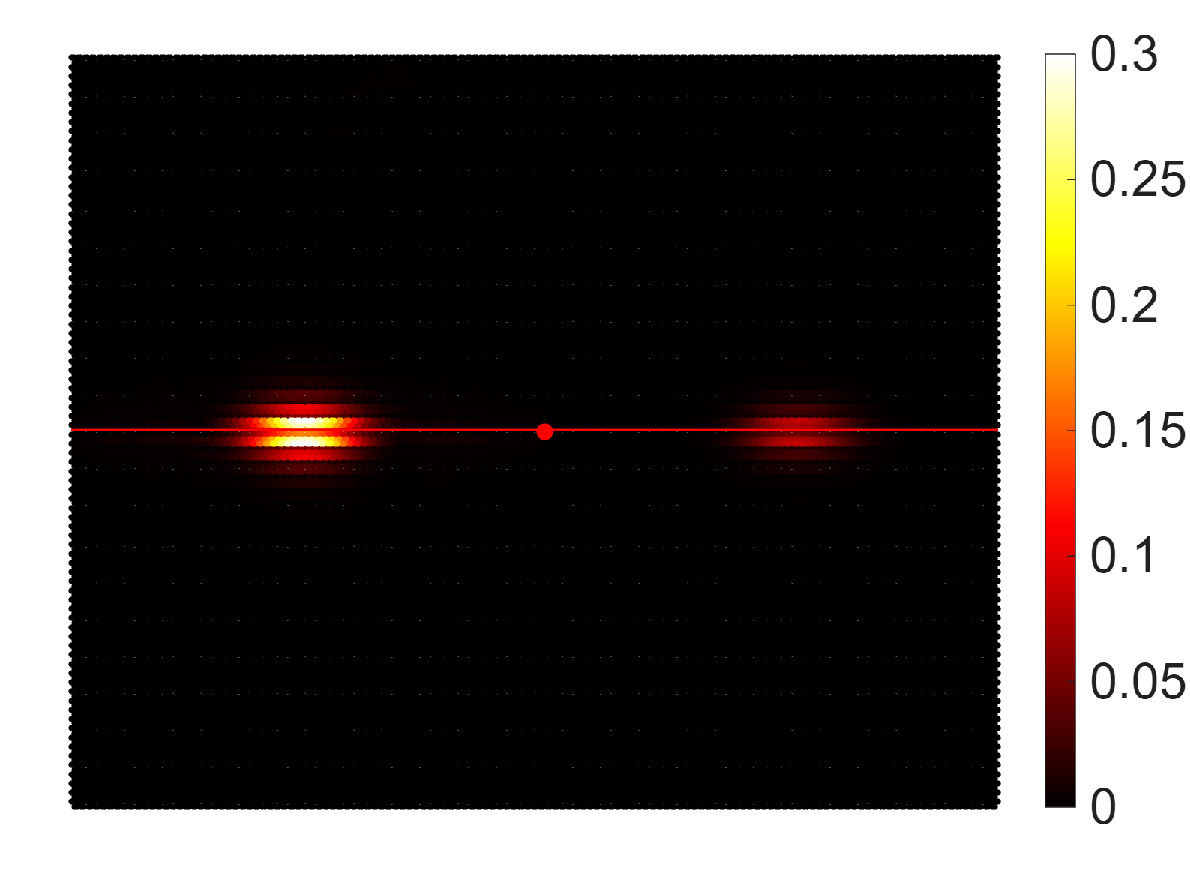}
				\caption{$ t=3 $.}
			\end{subfigure}
			\caption{Dynamics of wave packets passing through a local defect at a type-I interface when $\delta_{+}\delta_{-}>0$. In these two cases, there is no topological transition across the interface. Upper panels: $ \delta_{\pm} = 30$. Lower panels: $ \delta_{\pm} = -30$.}
			\label{fig:edgeIdef_same}
		\end{figure}

\begin{figure}[htbp] 
			\centering
			\begin{subfigure}{0.25\textwidth}
				\includegraphics[width=\textwidth]{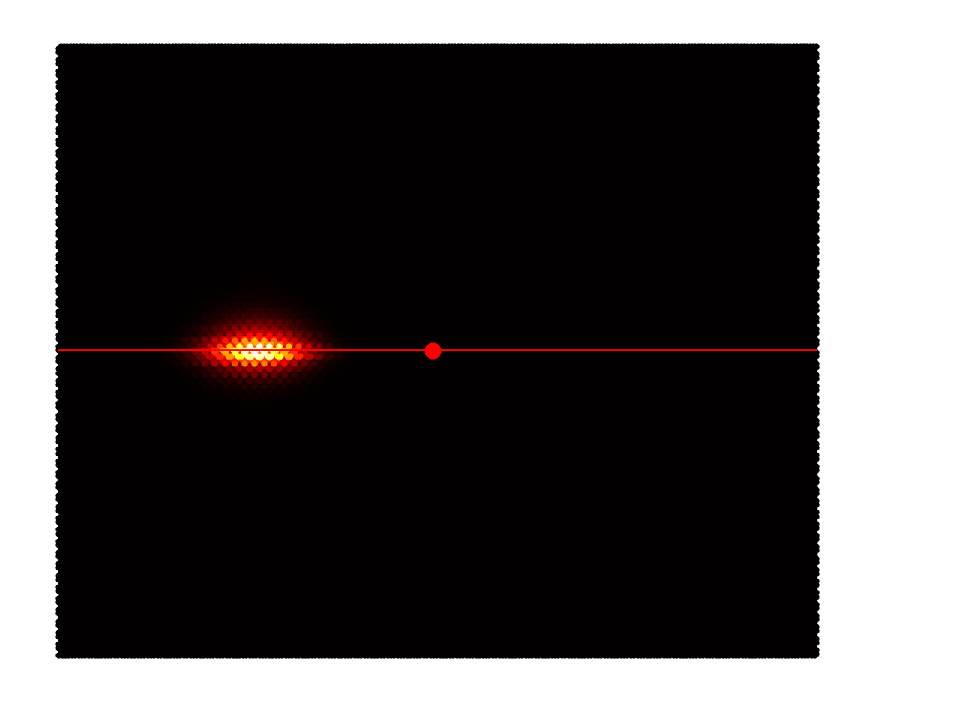}
				\caption{$ t=0 $.}
			\end{subfigure}
			\hspace{-0.6cm}
			\begin{subfigure}{0.25\textwidth}
				\includegraphics[width=\textwidth]{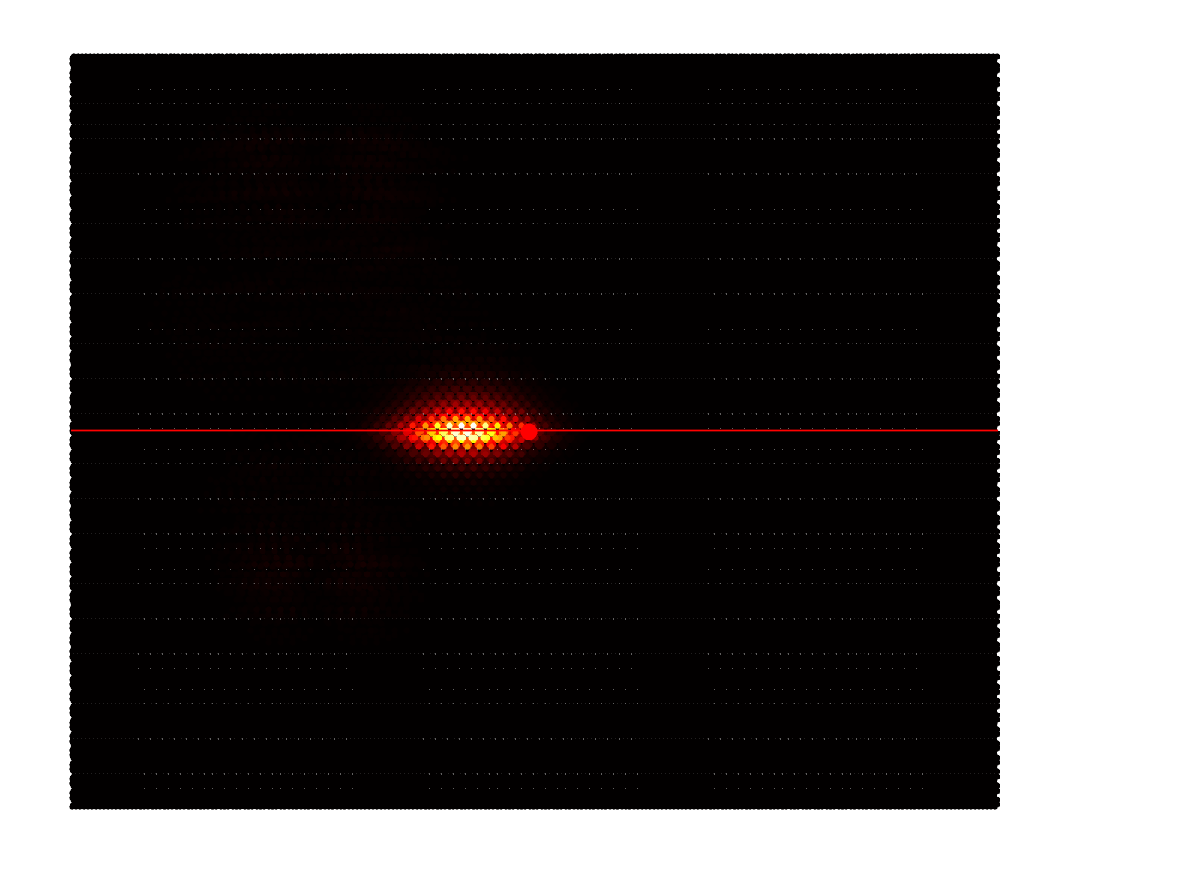}
				\caption{$ t=1 $.}
			\end{subfigure}
			\hspace{-0.6cm}
			\begin{subfigure}{0.25\textwidth}
				\includegraphics[width=\textwidth]{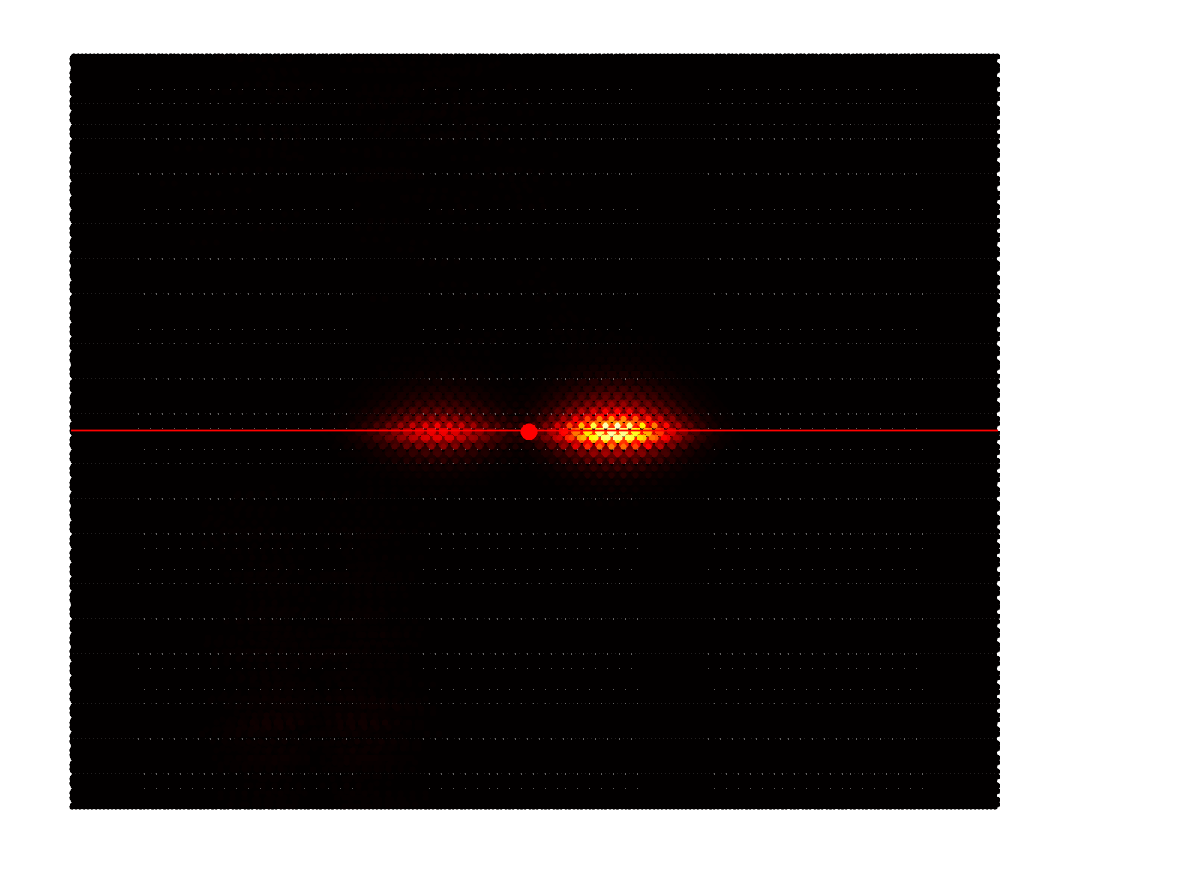}
				\caption{$ t=1.5 $.}
			\end{subfigure}	
			\hspace{-0.6cm}
			\begin{subfigure}{0.25\textwidth}
				\includegraphics[width=\textwidth]{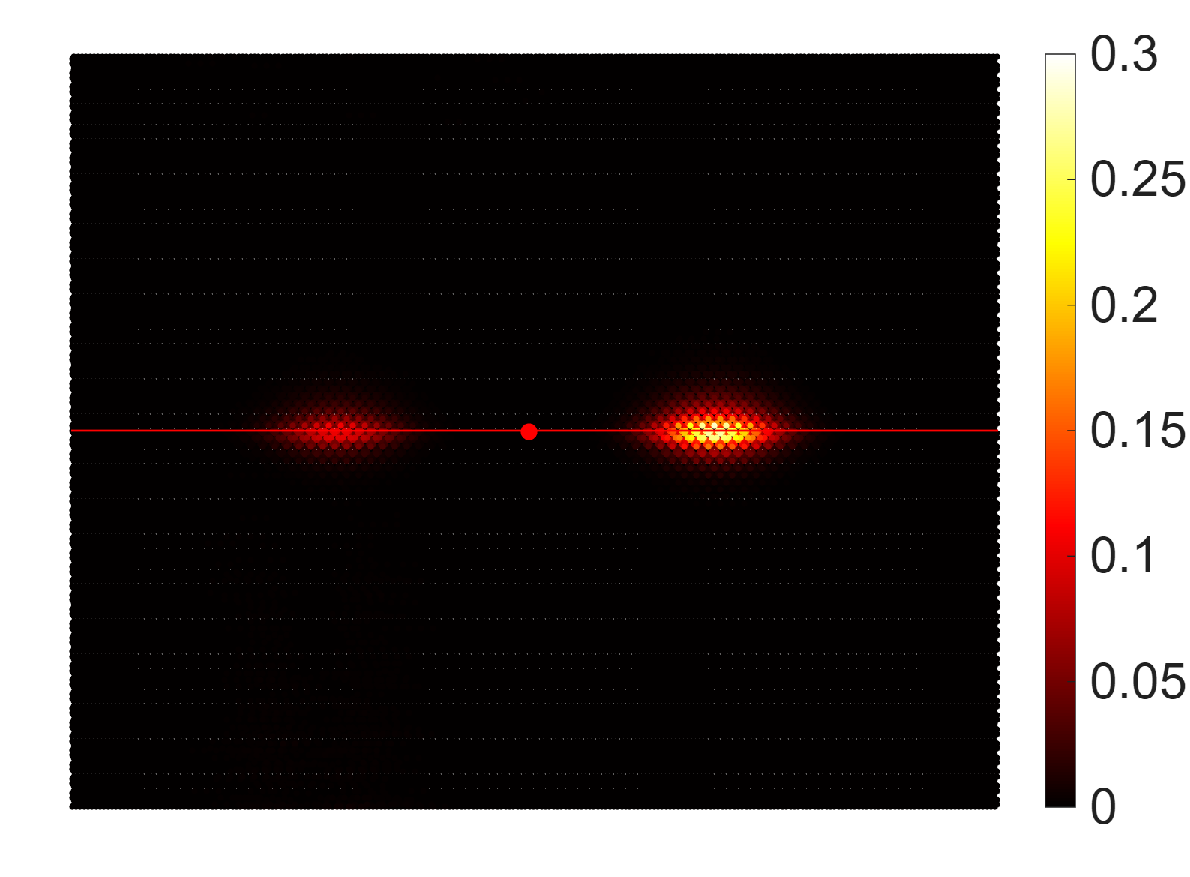}
				\caption{$ t=2 $.}
			\end{subfigure}
			\caption{Dynamics of wave packets passing through a local defect at a type-II interface.  }
			\label{fig:edgeIIdef}
		\end{figure}

		First, we investigate the dynamics of wave packets at type-I interfaces, when there is a topological transition across the interface. The time snapshots in \Cref{fig:edgeIdef_diff} show that the wave packet exhibits backscattering when passing through a local defect.
		We also numerically calculate the dynamics at a type-I interface when $\delta_{\pm} = 30$ and $\delta_{\pm}=-30$. The time snapshots are given in \Cref{fig:edgeIdef_same}. In these two cases, there is no topological transition across the interface.
		For the type-II interface, the snapshots are shown in \Cref{fig:edgeIIdef}.
		
\begin{figure}[htbp] 
			\centering
			\begin{subfigure}{0.48\textwidth}
				\includegraphics[width=\textwidth]{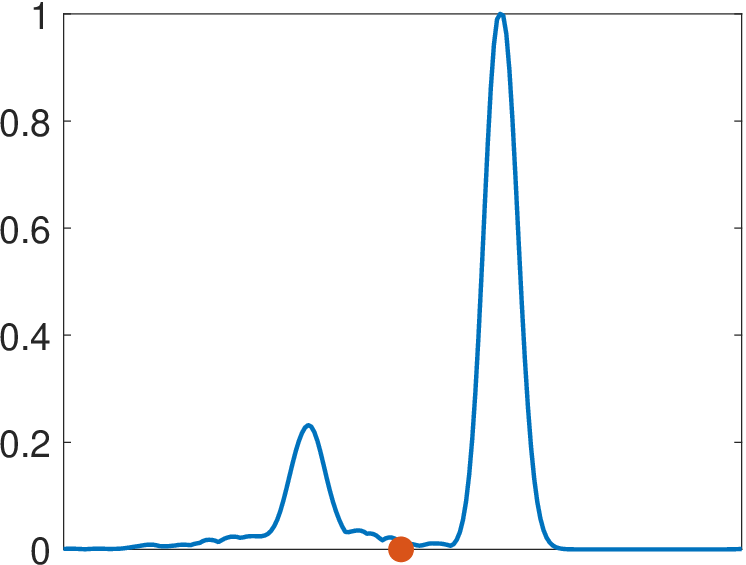}
				\caption{Type-I interface, $\delta_{\pm}=30$.}
			\end{subfigure}
			\begin{subfigure}{0.48\textwidth}
				\includegraphics[width=\textwidth]{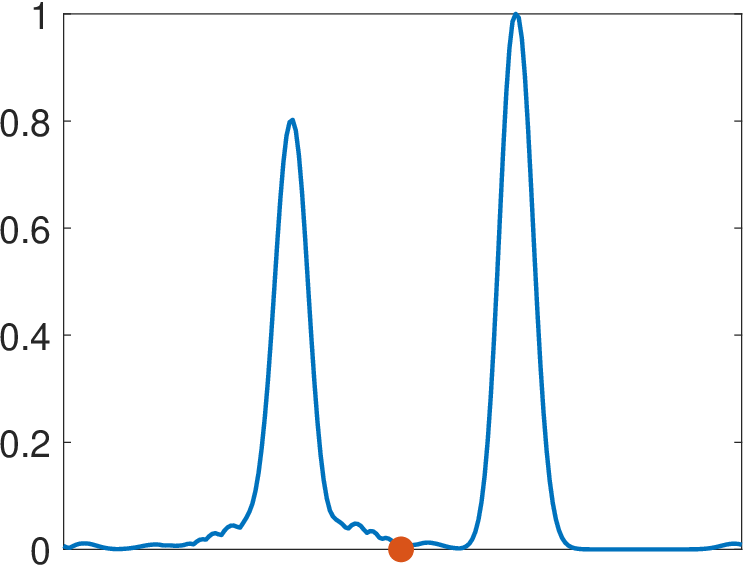}
				\caption{Type-I interface, $\delta_{\pm}=-30$.}
			\end{subfigure}\\
			\begin{subfigure}{0.48\textwidth}
				\includegraphics[width=\textwidth]{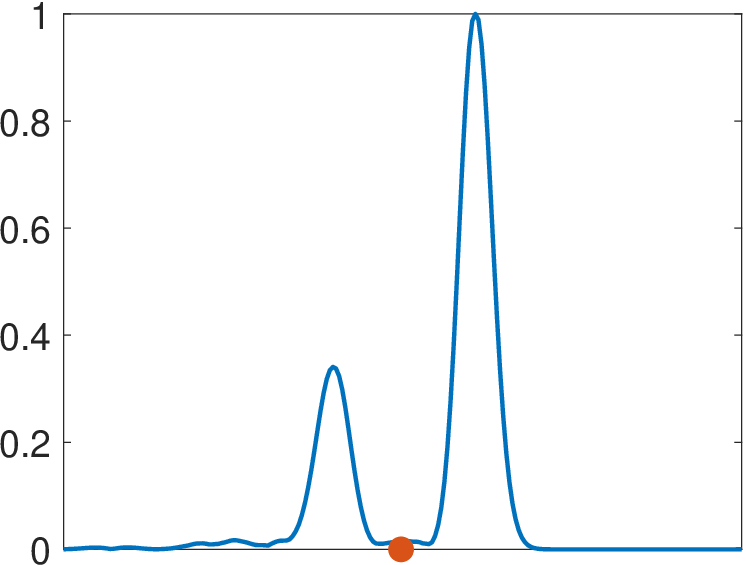}
				\caption{Type-I interface, $\delta_{\pm}=\pm30$.}
			\end{subfigure}	
			\begin{subfigure}{0.48\textwidth}
				\includegraphics[width=\textwidth]{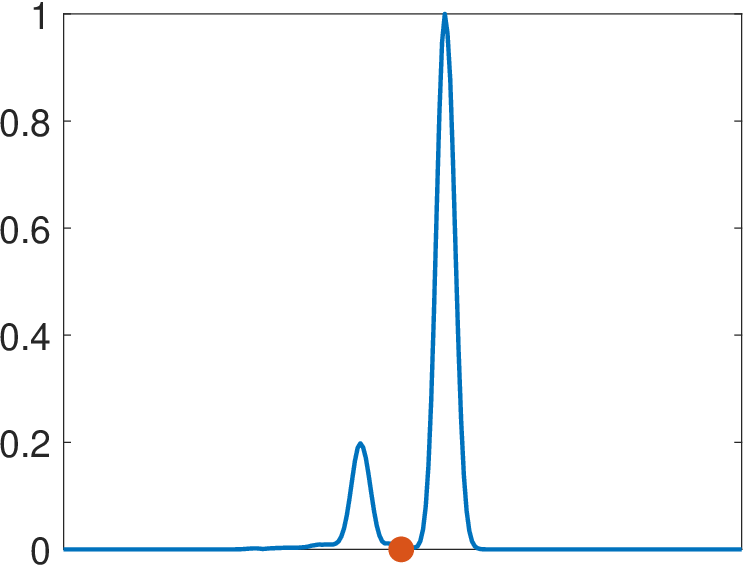}
				\caption{Type-II interface, $\delta_{\pm}=\pm30$.}
			\end{subfigure}
			\caption{The average maximum amplitude of 32 samples of wave packet dynamics.  }
			\label{fig:edgetransmission}
		\end{figure}

		To investigate the robustness of each interface, we simulate 32 random realizations and compute the maximum absolute value of the wave function $\boldsymbol{\Phi}$ at the final time in each vertical cell. We then average these quantities over the 32 realizations. After normalizing the maximum averaged amplitude to $1$, we obtain \Cref{fig:edgetransmission}. For the particular settings, the type-II interface exhibits a smaller reflected component than the corresponding opposite-sign type-I simulation. We wish to characterize this phenomenon in future work.

	\clearpage

	\begin{appendices}
		\section{Bulk Properties of the Operator}\label{apsec:bulk}
\begin{figure}[htbp] 
			\centering
			\begin{subfigure}{0.34\textwidth}
				\includegraphics[width=\textwidth]{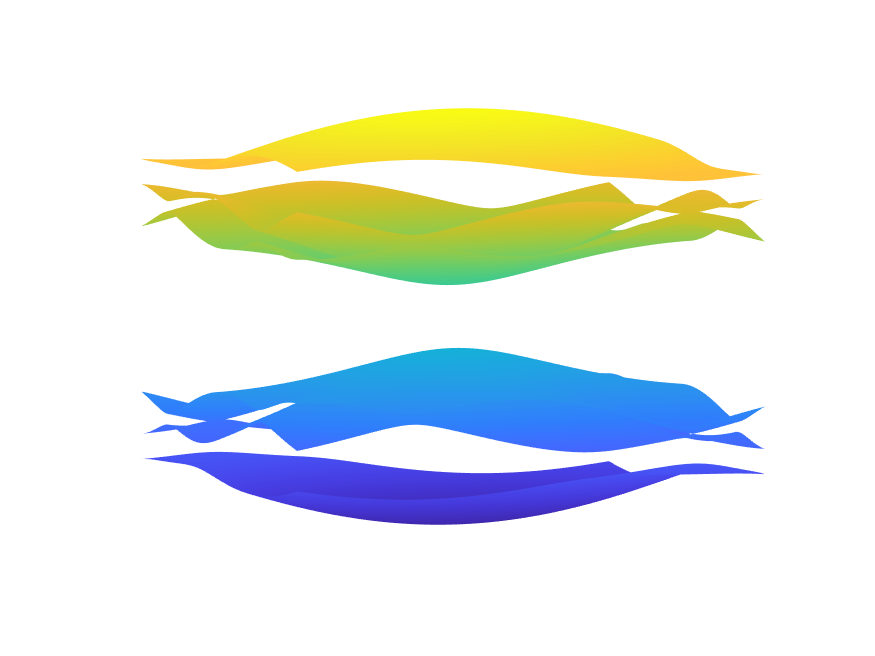}
				\caption{$ b = 5,\varepsilon=-2 $.}
			\end{subfigure}
			\hspace{-0.5cm}
			\begin{subfigure}{0.34\textwidth}
				\includegraphics[width=\textwidth]{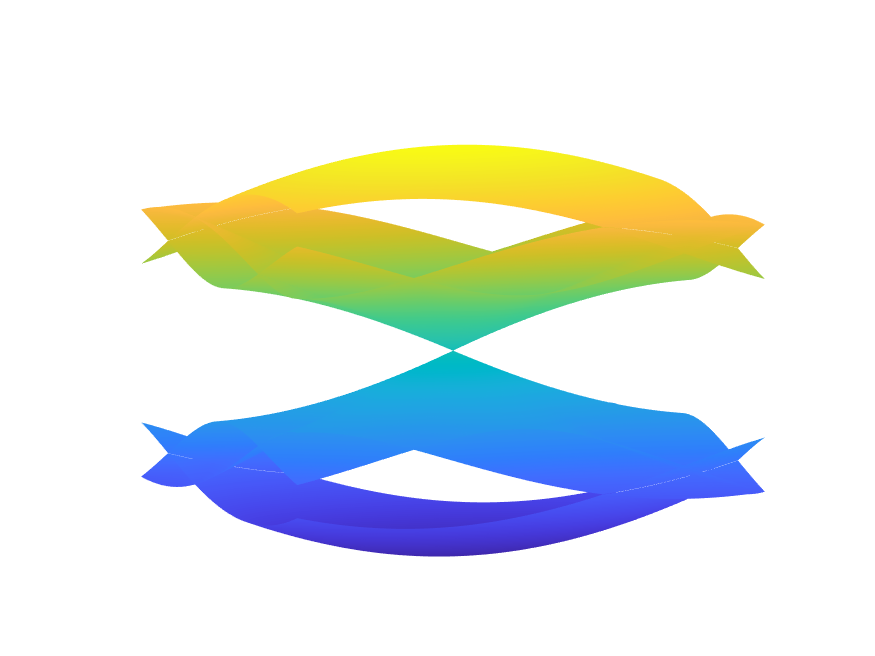}
				\caption{$ b = 5,\varepsilon=0 $.}
			\end{subfigure}
			\hspace{-0.5cm}
			\begin{subfigure}{0.34\textwidth}
				\includegraphics[width=\textwidth]{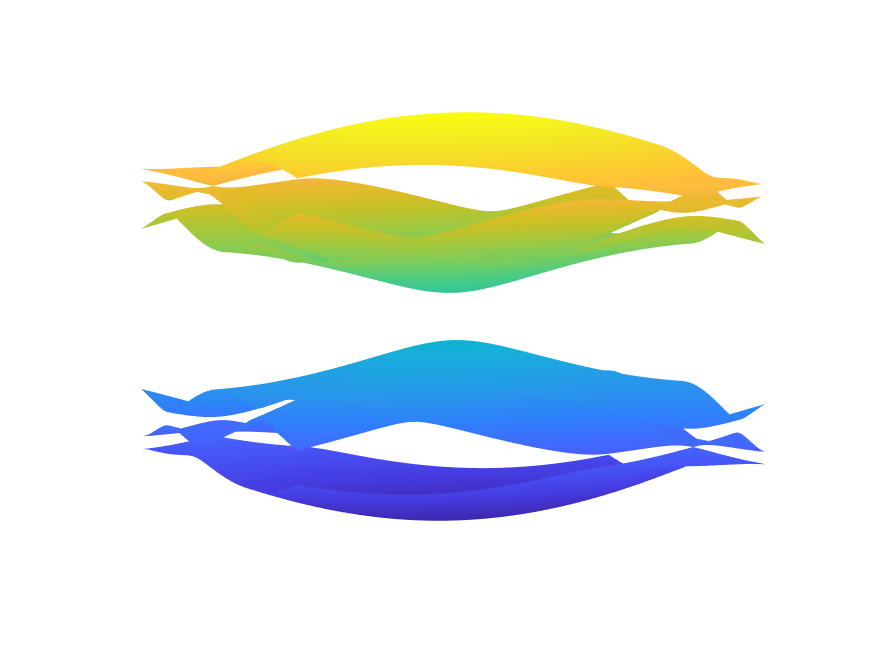}
				\caption{$ b = 5,\varepsilon=2 $.}
			\end{subfigure}
			\caption{The energy surfaces $ \protect\{ \lambda_{j}(\mathbf{k})\}_{j=1}^{6}$ for different values of $ \protect\varepsilon $. A double Dirac cone occurs when $ \protect\varepsilon=0 $, whereas a band gap opens when $ \protect\varepsilon\neq 0 $.}
			\label{fig:BandStructureTB}
		\end{figure}

		In this section, we suppose that the hopping coefficients $ \{ a_{_\lambda},b_{_\lambda},c_{_\lambda},d_{_\lambda} \}_{\lambda\in \Lambda} $ are given by 
		\begin{gather*}
			b_{_{\lambda}} = b,\quad a_{_\lambda}=c_{_\lambda}=d_{_\lambda} = b+\varepsilon.
		\end{gather*}
		The assumptions \eqref{eqn:HoppingCoef} ensure that $ b>0 $ and $ \varepsilon>-b $. We denote the operator by $ \hatHbulk = \hatHbulk(b,\varepsilon) $.
		\begin{remark}
			From a physical point of view, if the intercell interaction is larger than the intracell interaction, then the corresponding hopping coefficient $ b+\varepsilon $ is larger than $ b $, and vice versa. To realize it practically, one can dilate or contract the inclusions in super-honeycomb structures; see, for example, \cite{Miao2024}.
		\end{remark}
		Since the operator $\hatHbulk$ is periodic in $ \circv_{\alpha} $ and $ \circv_{\beta} $, one can decompose the Hilbert space $ l^2(\Lambda;\mathbb{C}^6) $ into a direct integral over the space $ \mathbb{C}^6_{\mathbf{k}} $, with quasi-momentum $ \mathbf{k}\in \{ s\mathring{\mathbf{k}}_{\alpha}+t\mathring{\mathbf{k}}_{\beta}:s,t\in[0,1) \} \triangleq Y^{\ast}$, where $ \mathring{\mathbf{k}}_{\alpha}, \mathring{\mathbf{k}}_{\beta} $ are dual basis vectors defined by 
		\[ \mathring{\mathbf{k}}_{\alpha}\cdot \circv_{\alpha} = \mathring{\mathbf{k}}_{\beta}\cdot \circv_{\beta} = 2\pi,\quad \mathring{\mathbf{k}}_{\alpha}\cdot \circv_{\beta} = \mathring{\mathbf{k}}_{\beta}\cdot \circv_{\alpha} = 0.  \]
		The Hamiltonian operator $ \hatHbulk $ can also be decomposed into a direct integral of operators $ \hatHbulk(\mathbf{k}) $. For $ \nu\in \mathbb{C}^6_{\mathbf{k}}  $, the Hamiltonian is given by 
		\begin{gather}
			\begin{split}
				(\widehat{H}_{\mathrm{bulk}}(\mathbf{k})\nu)_{_{1}} &\triangleq -b \nu_{_{4}}  -b\nu_{_{5}} -(b+\varepsilon)\eu^{-\iu \mathbf{k}\cdot \circv_{\alpha}}\nu_{_{6}},\\
				(\widehat{H}_{\mathrm{bulk}}(\mathbf{k})\nu)_{_{2}} &\triangleq -b \nu_{_{4}}  -(b+\varepsilon)\eu^{\iu \mathbf{k}\cdot(\circv_{\alpha}-\circv_{\beta})} \nu_{_{5}} -b\nu_{_{6}},\\
				(\widehat{H}_{\mathrm{bulk}}(\mathbf{k})\nu)_{_{3}} &\triangleq -(b+\varepsilon) \eu^{\iu \mathbf{k}\cdot \circv_{\beta}}\nu_{_{4}} - b \nu_{_{5}} -b\nu_{_{6}},\\
				(\widehat{H}_{\mathrm{bulk}}(\mathbf{k})\nu)_{_{4}} &\triangleq -b \nu_{_{1}}-b\nu_{_{2}} -(b+\varepsilon)\eu^{-\iu \mathbf{k}\cdot \circv_{\beta}} \nu_{_{3}},\\
				(\widehat{H}_{\mathrm{bulk}}(\mathbf{k})\nu)_{_{5}} &\triangleq -b \nu_{_{1}}-(b+\varepsilon)\eu^{-\iu \mathbf{k}\cdot(\circv_{\alpha}-\circv_{\beta})}\nu_{_{2}} -b \nu_{_{3}},\\
				(\widehat{H}_{\mathrm{bulk}}(\mathbf{k})\nu)_{_{6}} &\triangleq -(b+\varepsilon)\eu^{\iu \mathbf{k}\cdot \circv_{\alpha}} \nu_{_{1}}-b\nu_{_{2}} -b \nu_{_{3}}.
			\end{split}
		\end{gather}
		For $ \mathbf{k}\in Y^{\ast} $, one can verify the following fact about the eigenvalues of the Hamiltonian operator $ \hatHbulk(\mathbf{k}) $.
		\begin{proposition}\label{approp:bulk_gap}
			For $ \mathbf{k}\in Y^{\ast} $, the Hamiltonian operator $ \hatHbulk(\mathbf{k}) $ has six real eigenvalues $ \{ \lambda_{j}(\mathbf{k}) \}_{j=1}^6 $ in ascending order. Moreover, the eigenvalues satisfy 
			\[ \lambda_{j}(\mathbf{k}) \le -|\varepsilon|,\quad \lambda_{l}(\mathbf{k}) \ge |\varepsilon|,\quad j = 1,2,3,\;l = 4,5,6. \] 
			Specifically, when $ \mathbf{k}=\mathbf{0} $, one has
			\[ \lambda_{1}(\mathbf{0}) = -(3b+\varepsilon),\;\lambda_{2}(0) = \lambda_3(\mathbf{0}) = -|\varepsilon|\;,\lambda_{4}(\mathbf{0}) = \lambda_{5}(\mathbf{0}) = |\varepsilon|\;,\lambda_{6}(\mathbf{0}) =3b+\varepsilon.  \]
		\end{proposition}
		For $ \varepsilon=0 $, it can be proved easily that there exists a double Dirac point at $ \mathbf{k}=\mathbf{0} $.
		\begin{theorem}
			For $ \varepsilon=0 $, there exists a double Dirac point at $ \mathbf{k}=\mathbf{0} $. More specifically, there exists a positive constant $ \lambda^{\sharp} $ such that the following relations hold
			\begin{gather}
				\begin{split}
					&\lambda_{j}(\mathbf{k}) = -\lambda^{\sharp}|\mathbf{k}| + \mathcal{O}(|\mathbf{k}|^2),\quad j = 2,3,\\
					&\lambda_{l}(\mathbf{k}) = \lambda^{\sharp}|\mathbf{k}| + \mathcal{O}(|\mathbf{k}|^2),\quad l = 4,5.
				\end{split}
			\end{gather} 
		\end{theorem}
		We plot the eigenvalues $ \{ \lambda_{j}(\mathbf{k}) \}_{j=1}^6  $ as a function of $ \mathbf{k}\in Y^{\ast} $. We set $ b = 5 $ and $\varepsilon = -2,0,2$, as shown in \Cref{fig:BandStructureTB}. We can see that a double Dirac point occurs at $ \mathbf{k}=\mathbf{0} $ when $ \varepsilon =0$. For $ \varepsilon\neq 0 $, a band gap opens at $ \mathbf{k} = \mathbf{0} $.
		\par 
		Furthermore, one can prove the band inversion by direct calculation.
		\begin{proposition}
			For $ \mathbf{k} = 0 $, the Hamiltonian operator $ \hatHbulk(0) $ has six eigenvectors $ \{ \nu_{_{j}} \}_{j=1}^{6} $ that correspond to the eigenvalues $ \{ \lambda_{j}(\mathbf{0})\}_{j=1}^{6} $. Moreover, when $ \varepsilon<0 $, the second to fifth eigenvectors are given by 
			\begin{gather*}
				\begin{split}
					\nu_{_2} = (1,0,-1,1,0,-1)^{T},&\quad \nu_{_3} = (1,-2,1,-1,2,-1)^T,\\
					\nu_{_4} = (1,0,-1,-1,0,1)^{T},&\quad \nu_{_5} = (1, -2, 1, 1 , -2,1 )^T.
				\end{split}
			\end{gather*}
			When $ \varepsilon>0 $, the second to fifth eigenvectors are given by
			\begin{gather*}
				\begin{split}
					\nu_{_2} = (1,0,-1,-1,0,1)^{T},&\quad \nu_{_3} = (1, -2, 1, 1 , -2,1 )^T,\\
					\nu_{_4} = (1,0,-1,1,0,-1)^{T},&\quad \nu_{_5} = (1,-2,1,-1,2,-1)^T.
				\end{split}
			\end{gather*}
		\end{proposition}
		
		This eigenvector inversion motivates us to treat contraction and dilation as topologically distinct configurations in this work.
		\section*{Acknowledgments}
		This work was supported by the National Key R\&D Program of China (Grant No. 2021YFA0719200). 
		\section*{Data availability}
		The data that support the findings of this study are available from the authors upon reasonable request.
		\section*{Conflict of Interest}
		The authors declare no conflicts of interest.

	\end{appendices}
	
	\bibliography{references}
	
\end{document}